\documentclass{amsart}

\usepackage{amsmath}
\usepackage{amsfonts}
\usepackage{amssymb}
\usepackage{hyperref}
\usepackage{subfig}
\usepackage{graphicx}
\usepackage{array}
\usepackage{indentfirst}
\usepackage{cite}
\usepackage{esint}
\usepackage{enumerate}
\usepackage{bbm}

\newtheorem{theorem}{Theorem}[section]
\newtheorem{lemma}[theorem]{Lemma}
\newtheorem{corollary}[theorem]{Corollary}
\newtheorem{proposition}[theorem]{Proposition}

\theoremstyle{definition}
\newtheorem{definition}[theorem]{Definition}

\theoremstyle{remark}
\newtheorem{remark}[theorem]{Remark}

\numberwithin{equation}{section}

\newcommand{\rest}{
  \,\raisebox{-.127ex}{\reflectbox{\rotatebox[origin=br]{-90}{$\lnot$}}}\,}

\begin{document}

\title[On the singular set of free interface]{On the singular set of free interface \\ in an optimal partition problem}

\author{Onur Alper}
\address{Department of Mathematics, Purdue University \newline \indent 150 N.University Street, West Lafayatte, IN 47907 USA}
\email{oalper@purdue.edu}

\begin{abstract}
We study the singular set of free interface in an optimal partition problem for the Dirichlet eigenvalues.
We prove that its upper $(n-2)$-dimensional Minkowski content, and consequently, its $(n-2)$-dimensional Hausdorff measure are locally finite.
We also show that the singular set is countably $(n-2)$-rectifiable, 
namely it can be covered by countably many $C^1$-manifolds of dimension $(n-2)$, up to a set of $(n-2)$-dimensional Hausdorff measure zero.
Our results hold for optimal partitions on Riemannian manifolds and harmonic maps into homogeneous trees as well.
\end{abstract}


%\date{...}

%\dedicatory{...}

%\keywords{... }

\maketitle

\section{Introduction} \label{IntroSection}

\subsection{Optimal partition problem}
We consider the following optimal partition problem,
\begin{equation}
\inf_{\mathcal{P}_N ( \Omega)} \sum_{k=1}^N \lambda_1 \left ( \Omega_k \right ),
\label{EigenvalueMinimization}
\end{equation}
where $\lambda_1 \left ( \Omega_k \right )$ denotes the first Dirichlet eigenvalue of $\Omega_k$ with zero boundary data for $k= 1$,...,$N$, and
$\mathcal{P}_N (\Omega)$ denotes the partitions of a connected, open, bounded domain $\Omega \subset \mathbf{R}^n$ into $N$ pairwise disjoint, {\emph{admissible}} sets $\Omega_k$.
That is, $\Omega_k \cap \Omega_{k'} = \emptyset$ for $k \neq k'$, $\Omega_k \subset \Omega$ and $\Omega_k$ belongs to the {\emph{admissible class}} for every $k$, $k'= 1$,...,$N$.

The definition of {\emph{admissible}} sets is a subtle issue.
The existence of an optimal partition into pairwise disjoint, quasi-open sets $\Omega_k$ follows from a more general result by Bucur, {\emph{et al.} in \cite{BBH}.
We refer to this seminal work for the definition of quasi-open sets, as well as the existence results for various optimal partition problems.
In \cite{CTV05F} Conti, {\emph{et al.}} considered a similar yet slightly more general minimization problem, 
\begin{equation}
\inf_{\mathcal{P}_N ( \Omega )} \sum_{k=1}^N \left ( \lambda_1 \left ( \Omega_k \right ) \right )^p, \quad p > 0,
\label{EigenvalueMinimization2}
\end{equation}
where $\mathcal{P}_N (\Omega)$ denotes the partitions of $\Omega$ into pairwise disjoint, measurable sets. 
Among the various interesting results they proved are that the optimal partition in this class is in fact a collection of $N$ pairwise disjoint, open and connected sets $\Omega_k$,
and that the eigenfunctions satisfy useful extremality conditions.
They also combined these extremality conditions with their general regularity result in \cite{CTV05} to deduce the local Lipschitz regularity of eigenfunctions $u_1$, ..., $u_N$ for the optimal partition.
Below we state these results for $p=1$, which is the case we consider in this article.

\begin{theorem}[Conti, {\emph{et al.}}, \cite{CTV05, CTV05F}] \label{ExtremalityTheorem}
There exists a partition $ \{\Omega_1$, ..., $\Omega_N \} $ achieving \eqref{EigenvalueMinimization} over all open, connected, pairwise disjoint partitions of $\Omega$.
Furthermore, the associated $L^2$-normalized eigenfunctions $u_1$, $u_2$, ..., $u_N$ are Lipschitz continuous, and 
for $k = 1, 2, ..., N$, the following inequalities hold in $\mathcal{D}'(\Omega)$:
\begin{align}
- \Delta u_k & \leq \lambda_k u_k, \label{ExtremalityConditionA} \\
- \Delta \left ( u_k - \sum_{j \neq k} u_j \right ) & \geq \lambda_k u_k - \sum_{j \neq k} \lambda_j u_j. \label{ExtremalityConditionB}
\end{align} 
\end{theorem}

The proof of \eqref{ExtremalityConditionA} and \eqref{ExtremalityConditionB} involves considering the eigenfunctions corresponding to the solutions of \eqref{EigenvalueMinimization2}
for $p > 1$ and passing to a limit as $p \downarrow 1$.
(Note that by \cite[Eqn. 20 and Lemma 3.1]{CTV05F}, the positive coefficients $a_1$, ...., $a_N$ in the original statement \cite[Theorem 1.1]{CTV05F} are equal to $1$ in the limiting case $p=1$.)
The local Lipschitz regularity of eigenfunctions $u_1$, ..., $u_N$, follows from 
combining \eqref{ExtremalityConditionA} and \eqref{ExtremalityConditionB} with the general regularity result due to the same authors in \cite{CTV05}.
Moreover, the Lipschitz continuity of eigenfunctions is true up to the boundary, if we assume that $\partial \Omega$ is of class $C^{1}$, cf. \cite[Theorem 8.4]{CTV05}.
However, in the absence of this assumption, the best result is global H\"older bounds for every $\alpha \in (0,1)$, cf. \cite{N2010}.

Furthermore, in \cite{CTV05, CTV05F}, the study of interfaces between subdomains $\Omega_1$, $...$, $\Omega_N$, was initiated in the special case $\Omega \subset \mathbf{R}^2$.
In this article, we will focus on this free interface (and its junction points in particular) in arbitrary dimensions. 
Hence, we refer to \cite{Helffer, BNH} and the references therein for numerous fascinating problems and developments related to \eqref{EigenvalueMinimization} and its variants.
Moreover, optimal partition problems such as \eqref{EigenvalueMinimization} and its variants arise in the context of spatial segregation in reaction-diffusion systems,
as well as population models of competitive type. We refer to \cite{CTV05F, CTV05A} and the references therein for numerous interesting applications. 

\subsection{Free interface} \label{CLSection}
In order to analyze the free interface in higher dimensions, in \cite{CL07} Caffarelli and Lin considered the solutions to the constrained energy minimization problem,
\begin{equation}
\min_{v \in \mathcal{K}} \int_{\Omega} | \nabla v |^2 \mathrm{d}y,
\label{EquivalentProblem}
\end{equation}
where $\mathcal{K} = \left \{ v \in H_0^1 \left ( \Omega, \Sigma_N \right ) \; : \; \int_{\Omega} \left | v_k \right |^2 \mathrm{d}y =1, \; k = 1,2,...,N \right \}$, and
$$
\Sigma_N = \left \{ Y \in \mathbf{R}^N \, : \, \sum_{k \neq k'} Y_{k}^2 Y_{k'}^2 = 0, \; \mathrm{and} \; Y_k \geq 0, \; k=1,2,...,N \right \}. 
$$
We remark that $\Sigma_N$ is a singular target with nonpositive curvature in the sense of Alexandrov, cf. \cite{GS92}, and any homogeneous tree with $N$ branches,
that is a simplical tree with one vertex and $N$ edges attached to it, is clearly isometric to $\Sigma_N$.

We observe that $\pi_1 \left ( \Sigma_N \right ) = 0$ implies the existence of a Lipschitz retraction from a neighborhood of $\Sigma_N$ in $\mathbf{R}^N$ onto $\Sigma_N$. 
Consequently, the set $H^1_0 \left ( \Sigma_N \right )$ is nonempty, and Lipschitz continuous maps are dense in it. Therefore, problem \eqref{EquivalentProblem} has a solution.
We also note that if $u \in \mathcal{K}$ is a minimizer, there exist Lagrange multipliers $\lambda_1, ..., \lambda_N$ such that $u$ is a stationary point of the functional
\begin{equation}
F (u) = \int_{\Omega} | \nabla u |^2 \, \mathrm{d}x - \sum_{k=1}^N \lambda_k \int_{\Omega} \left | u_k \right  |^2 \, \mathrm{d}y,
\label{Functional}
\end{equation}
with respect to admissible domain and target variations. That is,
\begin{equation}
\begin{aligned}
\frac{d}{dt}  \vert_{t=0} F \left ( u \circ \Phi_t \right ) = 0,  \\
\frac{d}{dt}  \vert_{t=0} F \left ( \Psi_t \circ u \right ) = 0,
\end{aligned}
\label{CombinedVariations}
\end{equation}
where $\left \{ \Phi_t \right \}_{t \in ( - \epsilon, \epsilon )}$ is $C^1$-family of diffeomorphisms of $\Omega$ such that $\Phi_0 = id$ and $\Phi_t \vert_{\partial \Omega} = id$,
and $\left \{ \Psi_t \right \}_{t \in ( - \epsilon, \epsilon )}$ is a $C^1$- family of maps in $C_0^\infty \left ( \Sigma_N, \Sigma_N \right )$.
In particular, for every $k = 1, 2,..., N$, we have
\begin{equation}
- \Delta u_k = \lambda_k u_k \quad \mathrm{in} \; \left \{ x \in \Omega \, : \, u_k(x) > 0 \right \}.
\label{EulerLagrange}
\end{equation}

In \cite[Section 2]{CL07} Caffarelli and Lin observed the following:
Firstly, the minimal value in \eqref{EquivalentProblem} is equal to \eqref{EigenvalueMinimization} over quasi-open, pairwise disjoint sets.
Secondly, once extended to $\Omega$ as $u \equiv 0$ in $\Omega \backslash \Omega_k$, the eigenfunction corresponding to the eigenvalue $\lambda_1 \left ( \Omega_k \right )$ coincides with the $k^{th}$-component of a quasi-continuous representative of map $u \, : \, \Omega \to \Sigma_N$ solving \eqref{EquivalentProblem}.
Lastly, the Lagrange multiplier $\lambda_k$ in \eqref{Functional} is equal to $\lambda_1 \left ( \Omega_k \right )$, and $\Omega_k = \left \{ x \in \Omega \, : \, u_k(x) > 0 \right \}$ up to a set of capacity zero.
In short, the problems \eqref{EigenvalueMinimization} (over quasi-open sets first and also over open sets by the regularity theory in \cite{CTV05, CTV05F, CL07}) and \eqref{EquivalentProblem} are equivalent. 
In particular, $u \, : \, \Omega \to \Sigma_N$, where $u = \left ( u_1, ..., u_N \right )$ for $u_1, ..., u_N$, the eigenfunctions in Theorem \ref{ExtremalityTheorem} extended to $\Omega$ as above,
is a solution to the problem \eqref{EquivalentProblem}. We refer to \cite[Section 2]{CL07} for further details on this equivalence, and state the main consequence of \cite{CL07} below.

\begin{theorem}[Caffarelli-Lin, \cite{CL07}] \label{CaffarelliLinTheorem}
Suppose $u \, : \, \Omega \to \Sigma_N$ is the solution to the constrained minimization problem \eqref{EquivalentProblem}, corresponding to the collection of eigenfunctions $u_1$, ..., $u_N$ in Theorem \ref{ExtremalityTheorem}.
Then the following hold:
\begin{enumerate}[(i)]
 \item $u^{-1} \{0\}$ has Hausdorff dimension $(n-1)$. \label{TopDimEstimate}
 \item There exists a set $\mathcal{S}(u) \subset u^{-1} \{0\}$, which is relatively closed in $u^{-1} \{0\}$ and is of Hausdorff dimension $(n-2)$ at most, such that locally $u^{-1} \{0\} \backslash \mathcal{S}(u)$ is
 a $C^{1, \alpha}$ hypersurface for some $\alpha > 0$. \label{DecompositionRule}
\end{enumerate}
\end{theorem}

Firstly, we remark that for Theorem \ref{CaffarelliLinTheorem} ensures the exhaustion property $\overline{\Omega} = \cup_{k=1}^N \overline{\Omega_k}$, 
which is not {\emph{a priori}} clear for optimal partitions in Theorem \ref{ExtremalityTheorem}. Optimal partitions satisfying \eqref{TopDimEstimate} and \eqref{DecompositionRule} are called {\emph{regular partitions}}.
Secondly, the set $\mathcal{S}(u)$ in the statement of Theorem \ref{CaffarelliLinTheorem} corresponds to the junctions, where three or more subdomains come together.
Therefore, it is evident that the Hausdorff dimension estimate on $\mathcal{S}(u)$ is sharp, as long as $N \geq 3$.
Since $\mathcal{S}(u)$ is the central object of study in this article, we give a precise definition below.

\begin{definition} \label{SingularSetDefined}
For any map $u \, : \, \Omega \to \Sigma_N$, we define $\mathcal{S}(u)$, the {\emph{singular set}} of $u$, as the subset of $u^{-1} \{0\}$ such that 
for every $x \in \mathcal{S}(u)$ and every $r \in \left ( 0, \mathrm{dist}\left ( x,\partial \Omega \right ) \right )$, there exist at least three indices $k \in \left \{ 1, ..., N \right \}$ satisfying 
\begin{equation*}
B_r(x) \cap \left \{ x \in \Omega \, : \, u_k(x) > 0 \right \} \neq \emptyset. 
\end{equation*}
\end{definition}

Consequently, the set $u^{-1} \{0\} \backslash \mathcal{S}(u)$ corresponds to {\emph{domain walls}} separating precisely two subdomains. 
We will recall in Corollary \ref{DimensionReduction} one of the key ingredients of Theorem \ref{CaffarelliLinTheorem}, which is that the map $u$ has vanishing order 1 at $z \in u^{-1} \{0\} \backslash \mathcal{S}(u)$,
while it has vanishing order greater than or equal to $1 + \delta_n$ at any $z \in \mathcal{S}(u)$ for some dimensional constant $\delta_n > 0$.
We remark that $\delta_2 = 1/2$ was already proved in \cite{CTV05F}, and Helffer, {\emph{et al.}} \cite{HHT} proved that $\delta_3 = 1/2$. Determining whether $\delta_n = 1/2$ for $n \geq 4$ is an open problem.

There are several extensions of Theorem \ref{CaffarelliLinTheorem}.
Building upon \cite{CL09}, Snelson proved in \cite{Snelson} an analogue of Theorem \ref{CaffarelliLinTheorem} for nonlocal heat flows corresponding to the constrained minimization problem \eqref{EquivalentProblem}.
Ramos, {\emph{et al.}} proved in \cite{RTT} the analogue of Theorem \ref{CaffarelliLinTheorem} for functionals depending on higher eigenvalues, 
despite the absence of extremal conditions \eqref{ExtremalityConditionA} and \eqref{ExtremalityConditionB}.
On the other hand, an extension that is very useful for our purposes is due to Tavares and Terracini \cite{TT12}, proving the analogue of Theorem \ref{CaffarelliLinTheorem} for segregated critical configurations
that are assumed to satisfy a weak reflection law, that is the Pohozaev identity \eqref{RenormalizedEnergy}.
Combining this result with the Lipschitz regularity result in \cite[Section 4]{N2010}, implies the conclusion of Theorem \ref{CaffarelliLinTheorem}, 
in particular the Hausdorff dimension estimate for the singular set, for strong limits solutions to \eqref{EquivalentProblem}.

\subsection{Main results} \label{MainResultsStated}
In the special case, $n=2$, the fact that $\mathcal{S}(u)$ is a locally discrete set was established in \cite{CTV05, CTV05F}. Hence, for the rest of this article we assume that $n \geq 3$.
In general, Theorem \ref{CaffarelliLinTheorem} does not contain information on the structure or size of the singular set $\mathcal{S}(u)$, other than the conclusion that $\mathcal{S}(u)$ has Hausdorff dimension $(n-2)$ .
In particular, it leaves open the question of whether it has locally finite $(n-2)$-dimensional Hausdorff measure.
This is the first of the two questions we address in this article, while the rectifiability of singular set $\mathcal{S}(u)$ is the second.
Before stating our main results, we recall the definition of upper $m$-dimensional Minkowski content.

\begin{definition} \label{UMCdefined}
For $A \subset \mathbf{R}^n$, $0 \leq m \leq n$, the upper $m$-dimensional Minkowski content of $A$ is defined as
\begin{equation}
M^{*,n-m} (A) = \limsup_{r \downarrow 0} \frac{\mathcal{L}^n \left ( B_r \left ( A \right ) \right )}{ \omega(n-m) r^{n-m}},
\label{UMC}
\end{equation}
where $B_r(A) = \left \{ x \in \mathbf{R}^n \, : \, \mathrm{dist}(x,A) < r \right \}$ for $\mathrm{dist}(x,A) = \inf_{y \in A} |x-y|$,
$\mathcal{L}^n$ is the $n$-dimensional Lebesgue measure, and $\omega(n-m)$ is the volume of unit ball in $\mathbf{R}^{n-m}$.
\end{definition}

Now we are ready to state the main results of this article.

\begin{theorem} \label{MainTheorem1}
Suppose $u \, : \, \Omega \to \Sigma_N$ is a solution to the constrained minimization problem \eqref{EquivalentProblem}, 
corresponding to the collection of $L^2$-normalized eigenfunctions in Theorem \ref{ExtremalityTheorem}.
For any compact $K \subset \Omega$, there exist positive $C \left (u,K, n, \lambda_M \right )$ and $r_0 \left ( n, N, \lambda_M, u \right ) < \mathrm{dist} \left ( K, \partial \Omega \right )$, 
where $\lambda_M = \max_{1 \leq k \leq N} \lambda_k$,
such that for every $r \in \left ( 0, r_0 \right)$,
\begin{equation}
\mathcal{L}^n \left ( B_r \left ( \mathcal{S}(u) \right ) \cap K \right ) \leq C(u,K) r^2,
\label{StrongestSizeEstimate}
\end{equation}
In particular, $\mathcal{S}(u)$ has locally finite upper $(n-2)$-dimensional Minkowski content and $(n-2)$-dimensional Hausdorff measure.
\end{theorem}

The estimate \eqref{StrongestSizeEstimate} clearly yields a bound on the upper $(n-2)$-dimensional Minkowski content of $\mathcal{S}(u) \cap K$.
We also remark that the upper $k$-dimensional Minkowski content of a set controls its $k$-dimensional Hausdorff measure from above, cf. \cite[Proposition 3.3.3]{Mattila}.
In fact, this essentially amounts to bounding the $(n-2)$-dimensional Hausdorff measure of $\mathcal{S}(u) \cap K$ from above by using \eqref{StrongestSizeEstimate} directly.

Our second result is on the structure of singular set $\mathcal{S}(u)$.

\begin{theorem} \label{MainTheorem2}
Suppose $u \, : \, \Omega \to \Sigma_N$ is a solution to the constrained minimization problem \eqref{EquivalentProblem}, 
corresponding to the collection of $L^2$-normalized eigenfunctions in Theorem \ref{ExtremalityTheorem}.
Then $\mathcal{S}(u)$ is countably $(n-2)$-rectifiable, namely it can be covered by countably many $C^1$-manifolds of dimension $(n-2)$, up to a set of $(n-2)$-dimensional Hausdorff measure zero.
\end{theorem}

In addition, up to minor modifications, the proofs of Theorems \ref{MainTheorem1} and \ref{MainTheorem2} also yield the following result.

\begin{theorem} \label{MainTheorem3}
For $(M,g)$, a Riemannian manifold and its metric, and $\Omega \subseteq M$, a smooth and bounded domain,
the conclusions of Theorem \ref{MainTheorem1} (with a constant $C = C(g,u,K)$ in \eqref{StrongestSizeEstimate} instead) and Theorem \ref{MainTheorem2} apply to $\mathcal{S}(u)$ in the following cases as well:
\begin{enumerate}[(A)]
 \item $ \{\Omega_1$, ..., $\Omega_N \}$ is a partition of $\Omega$ into open, connected and pairwise disjoint subsets,
  that minimizes \eqref{EigenvalueMinimization}, where $\lambda_1 \left ( \Omega_k \right )$ denotes the first Dirichlet eigenvalue -with respect to the Laplace-Beltrami operator $\Delta_g$ on $(M,g)$-
  of $\Omega_k$ with zero boundary data, for $k= 1$,...,$N$, and $u \, : \, \Omega \to \Sigma_N$ is the corresponding $N$-tuple of $L^2$-normalized eigenfunctions. \label{Manifolds}
 \item $u \, : \, \Omega \to \Sigma_N$ is a harmonic map that minimizes the Dirichlet energy in $\Omega$
 with respect to the boundary condition $u = g$ on $\partial \Omega$ for a map $g \in H^{1/2} \left ( \partial \Omega, \Sigma_N \right )$. \label{Trees}
\end{enumerate}
\end{theorem}

\subsection{Strategy of the proof} \label{StrategyStated}
We adopt the constrained mapping approach of Caffarelli and Lin to the optimal partition problem as in \cite{CL07,CL08,CL10}.
We emphasize that $\Sigma_N$ is a singular target, and the key dimensionless quantities in our analysis are the classical Almgren frequency function \eqref{ClassicalAlmgren} 
and its variants \eqref{ClassicalGenFrequency}, \eqref{SmoothedFrequency}, \eqref{SmoothedGenFrequency}, \eqref{SmoothedFrequencyAdditive}.
The classical Almgren frequency function was originally introduced in \cite{Almgren83} in the analysis of a special class of maps into singular targets, namely harmonic $\mathcal{Q}$-valued functions.
Since then the frequency function and its variants have been heavily utilized in the study of maps into various other singular targets.
A few examples are round cones arising from liquid crystal theory \cite{L89,L91}, spaces of nonpositive curvature in the sense of Alexandrov \cite{GS92},
and piecewise uniformly regular Lipschitz manifolds \cite{Wang}. We remark that $\Sigma_N$ also has nonpositive curvature in the sense of Alexandrov.

In \cite{NV, NV2} Naber and Valtorta introduced several powerful techniques for studying the size and structure of stationary maps into Riemannian manifolds.
These techniques were extended by De Lellis, {\emph{et al.}} in \cite{DLMSV} to the context of energy-minimizing harmonic $\mathcal{Q}$-valued functions.
As motivated by the Ericksen model in liquid crystal theory, the singular set of energy-minimizing maps into cones over the real projective plane
was studied in \cite{A18} by combining the approach of \cite{DLMSV} with the blowup analysis in \cite{AHL}. 
In this article, we adapt these techniques to the context of maps that are stationary points of the functional \eqref{Functional} under the variations \eqref{CombinedVariations}.

The analysis in both \cite{DLMSV} and \cite{A18} essentially concerns the preimage of vertex of a conical target, while the specific features differ due to different properties of respective targets and blowup maps.
Likewise, Focardi and Spadaro adapted the techniques of \cite{DLMSV} to thin obstacle problems in \cite{FS, FS2}, where they study the free boundary in its entirety. 
In contrast, we do not consider the preimage of $0 \in \Sigma_N$ as a whole in this article.
As the top dimensional part of nodal set has already been addressed in Theorem \ref{CaffarelliLinTheorem} and enjoys relatively nice local regularity,
instead we focus on the $(n-2)$-dimensional part of nodal set, which corresponds to the junctions where three or more {\emph{domain walls}} meet.

Note that in \cite{EE} Edelen and Engelstein also studied the lower dimensional strata of free boundary in various classes of free boundary problems by using the techniques of \cite{NV}.
Their work is closer in spirit to \cite{NV, NV2}, as they resort to the quantitative stratification of free boundary.
We remark that the problems they consider enjoy the property that their blowups are $1$-homogeneous functions.
In contrast, the problem we consider gives rise to cylindrical blowups with homogeneity $k/2$ for $k=3,4,5...$. See Section \ref{HomSubsection} ahead.
On the other hand, little is known about the spectrum of admissible vanishing orders for blowups of $3$ or more variables at the points in $\mathcal{S}(u)$,
and such points comprise the lower strata of $\mathcal{S}(u)$. 
Therefore, we do not implement a quantitative stratification of $\mathcal{S}(u)$ in this article, and instead make use of a gap condition for the frequency function (see Section \ref{DecompositionSection} ahead) 
in order to distinguish points in $\mathcal{S}(u)$ and prove their persistence within $u^{-1} \{0\}$ in compactness arguments.

We remark that the most natural class of maps $u \, : \, \Omega \to \Sigma_N$ to consider in our analysis are those, which are $C^{0,\alpha}$-regular for every $\alpha \in \left ( 0, 1 \right )$
and stationary with respect to the variations \eqref{CombinedVariations} of functional \eqref{Functional}. 
In fact, such maps do have local Lipschitz regularity as well. See Section \ref{LipschitzSub} ahead. 
Hence, without loss of generality, we will work with Lipschitz continuous maps that are stationary with respect to the variations \eqref{CombinedVariations} of functional \eqref{Functional}. 
Such a class of maps was also discussed in a more general context in \cite{TT12}. They have good compactness properties and enjoy the variational properties that are indispensable to our analysis.
And most importantly, the globally constrained minimizers \eqref{EquivalentProblem} and their weak limits are contained in this class. See Section \ref{CompactnessSubsection} ahead.

The proof combines first variation formulas, compactness arguments and spectral analysis for the homogeneous blowup maps with techniques from geometric measure theory.
Due to the presence of lower order term arising from the Lagrange multipliers in \ref{Functional}, we will restrict our attention to suitably small scales, 
as the (almost) monotonicity formulas we derive are valid at sufficiently small scales essentially depending on $n$ and $\lambda_M = \max_{1 \leq k \leq N} \lambda_k$.
We note that if we consider stationary maps with respect to \eqref{Functional} with $\lambda_1 = ... = \lambda_N = 0$, then these formulas are valid at all scales.
Consequently, the estimates we obtain through the first variation formulas will also be applied at sufficienty small scales.
While we can scale our maps to normalize the scale at which the variational formulas hold, this operation translates to a smallness requirement for $\lambda_M$.
Since the eigenvalues of subdomains in \eqref{EigenvalueMinimization} are at the heart of our analysis, we refrain from such a scaling argument.
Instead we adopt an explicit approach to tracking the spatial scales throughout our analysis.

A second restriction due to the lower order term in \eqref{Functional} concerns the logarithmic derivatives of various quantities that are central to our analysis.
In order to control such derivatives in variational formulas, we will often restrict our attention to the free interface of optimal partitions,
as the blowup analysis yields lower bounds for the Almgren frequency function and its variants at the points on the free interface. See Sections \ref{PrelimSection} and \ref{SmoothSection} ahead.
In the case $\lambda_1 = ... = \lambda_N = 0$ considered in Section \ref{FinalRemarksSection}, we do not need such a restriction either.

\subsection{Future directions} \label{FutureDirections}
We plan to extend in future works our analysis of stationary maps into the singular target $\Sigma_N$ to more general settings with the following applications in mind:
\begin{enumerate}[(1)]
 \item spatial segregation in reaction-diffusion systems, as well as competitive systems in population dynamics as in \cite{CTV05, TT12},
 \item Bose-Einstein condensates in multiple hyperfine spin states as in \cite{CLLL, TT12},
 \item optimal partition problems involving higher eigenvalues as in \cite{RTT},
 \item harmonic maps into singular targets of curvature bounded from above in the sense of Alexandrov and their geometric applications as in \cite{BFHMSZ} and the references therein.
\end{enumerate}

\subsection{Plan of the article} \label{Plan}
\begin{itemize}
 \item In Section \ref{PrelimSection} we collect various useful results such as: 
regularity and compactness properties, first variation formulas and (almost) monotone quantities, homogeneous blowups and the gap condition for vanishing order at the singular points.
 \item In Section \ref{SmoothSection} we introduce the smoothed versions of quantities from the previous section, as well as additional first variation formulas, and prove several basic estimates.
 \item In Section \ref{MainFreqEstimatesSection} we prove the generalization of a crucial {\emph{frequency pinching estimate}} that originated in \cite{DLMSV} to the setting of functional \eqref{Functional}.
 \item In Section \ref{JonesNumberSection} we prove a pointwise bound on the Jones' $\beta_2$-number, which controls the deviation from flatness for a set in the Euclidean space, in the spirit of \cite{NV,NV2,DLMSV}.
 \item In Section \ref{SpineSection} we prove important technical lemmas concerning almost homogeneous maps, their singular sets and the oscillations of their frequencies on the singular set.
 \item In Section \ref{MinkowskiSection} we prove Theorem \ref{MainTheorem1} through inductive covering arguments in the spirit of \cite{NV2, DLMSV}.
 \item In Section \ref{RectifiabilitySection} we prove Theorem \ref{MainTheorem2} by combining the estimates from the previous section with a characterization of rectifiable measures by Azzam and Tolsa \cite{AT15}.
 \item Finally, in Section \ref{FinalRemarksSection} we briefly discuss how our proofs can be adapted to conclude Theorem \ref{MainTheorem3} as well.
\end{itemize}

\section{Preliminaries} \label{PrelimSection}

In this section for the reader's convenience we collect several results from \cite{CTV05, CTV05F, CL07, TT12}, which we will routinely invoke in the rest of the article.
We include the proofs whenever the ideas they involve are useful on their own for the following sections, while referring to these works otherwise. 

\subsection{Lipschitz regularity of eigenfunctions} \label{LipschitzSub}

We briefly review the regularity theory for eigenfunctions corresponding to optimal partitions as in \cite{CTV05} and \cite{CTV05F}.
It relies heavily on two main ingredients: the observation that the eigenfunctions in optimal partitions satisfy extremality conditions \eqref{ExtremalityConditionA} and \eqref{ExtremalityConditionB},
and a version of the well-known Alt-Caffarelli-Friedman monotonicity formula, cf. \cite{ACF}, proved in \cite{CJK}.

\begin{lemma} \label{LocalLipschitzOptimal}
Suppose $u \, : \, \Omega \to \Sigma_N$ is a solution to the constrained minimization problem \eqref{EquivalentProblem}. 
Then $u$ is a locally Lipschitz continuous map in $\Omega$. 
Furthermore, if $\partial \Omega$ is of class $C^1$, then $u \in W^{1, \infty} \left ( \overline{\Omega} \right )$.
\end{lemma}

\begin{proof}
By \cite[Section 2]{CL07}, each component of $u$ is an eigenfunction corresponding to the partition realizing \eqref{EigenvalueMinimization}.
Hence, by \cite[Theorem 1.1]{CTV05F}, $u$ satisfies the extremality conditions \eqref{ExtremalityConditionA} and \eqref{ExtremalityConditionB}.
Therefore, \cite[Theorem 8.3]{CTV05} applies, and consequently, $u$ is a locally Lipschitz continuous map in $\Omega$.
Likewise, the Lipschitz continuity up to the boundary, under the additional assumption $\partial \Omega$ is of class $C^1$, follows from \cite[Theorem 8.4]{CTV05}.
\end{proof}

For the compactness result in Section \ref{CompactnessSubsection}, we also need the following observation due to Noris, {\emph{et al.}} \cite{N2010}. 

\begin{lemma} \label{LocalLipschitzStationary}
Suppose $u \, : \, \Omega \to \Sigma_N$ is a stationary map with respect to the variations \eqref{CombinedVariations} of functional \eqref{Functional}.
If $u \in C^{0, \alpha}_{loc} \left ( \Omega \right )$ for every $\alpha \in (0,1)$, then $u$ is locally Lipschitz continuous in $\Omega$.
\end{lemma}

\begin{proof}
The proof relies on a combination of the variational formulas \eqref{RenormalizedEnergy}, \eqref{ClassicalFreqScalarDeriv} and \eqref{ClassicalLogDerivHeight}, 
which are stated and reviewed in Section \ref{FirstVariationsSub},
and the proof of Lemma \ref{LocalLipschitzOptimal} in \cite[Section 8]{CTV05}.  We refer to \cite[Section 4]{N2010} for details.
\end{proof}

\begin{remark} \label{TT12Remark}
In fact, one can further relax the assumptions in Lemma \ref{LocalLipschitzStationary}, 
and obtain the conclusions of Theorems \ref{ExtremalityTheorem} and \eqref{CaffarelliLinTheorem}, for any collection of functions
$u_1$, ..., $u_K$ with disjoint supports in $\Omega$, satisfying \eqref{ExtremalityConditionA}, \eqref{ExtremalityConditionB} and \eqref{RenormalizedEnergy}, 
by combining the main results of \cite{CTV05} and \cite{TT12}.
\end{remark}

\subsection{Variational formulas} \label{FirstVariationsSub}

Below we introduce the quantities that play a central role in our analysis.

\begin{definition} \label{ClassicalQuantities}
The classical localized key functional is
\begin{equation}
F(x,r) = \int_{B_r(x)} | \nabla u |^2 \, \mathrm{d}x - \sum_{k=1}^N \lambda_k \int_{B_r(x)} \left | u_k \right  |^2 \, \mathrm{d}y.
\label{LocalFunctional}
\end{equation}
The classical Dirichlet energy is
\begin{equation}
D(x,r) = \int_{B_r(x)} | \nabla u |^2 \, \mathrm{d}y.
\label{Classical Dirichlet}
\end{equation}
The classical height function is
\begin{equation}
H(x,r) = \int_{\partial B_r(x)} |u|^2 \, \mathrm{d}\sigma(y).
\label{Height}
\end{equation}
The classical Almgren frequency is
\begin{equation}
I(x,r) = \frac{r D(x,r)}{H(x,r)}.
\label{ClassicalAlmgren}
\end{equation}
The classical generalized frequency is
\begin{equation}
G(x,r) = \frac{r F(x,r) + H(x,r)}{H(x,r)}.
\label{ClassicalGenFrequency}
\end{equation}
\end{definition}

\begin{lemma} \label{FirstVariations}
Suppose $u \, : \, \Omega \to \Sigma_N$ is stationary with respect to the variations \eqref{CombinedVariations} of functional \eqref{Functional}, $\Omega \subset \mathbf{R}^n$, 
$\varphi \in C_c^\infty \left ( \Omega; \mathbf{R}^n \right )$, and  $\psi \in C_c^\infty \left ( \Omega \right )$.
Then we have the following distributional identities:
\begin{equation}
\int_{\Omega} \left ( 2 \left \langle \nabla u \otimes \nabla u , \nabla \varphi \right \rangle - |\nabla u|^2 \mathrm{div}\varphi 
+ \sum_{k=1}^N \lambda_k \left | u_k \right |^2 \mathrm{div}\varphi \right ) \,
\mathrm{d}y = 0,
\label{DomainVariation}
\end{equation}
\begin{equation}
\int_{\Omega} \left ( \left \langle \nabla u , \nabla \psi u  \right \rangle + | \nabla u |^2 \psi 
- \sum_{k=1}^N \lambda_k \left | u_k \right |^2 \psi \right ) \,
\mathrm{d}y = 0.
\label{TargetVariation}
\end{equation}
\end{lemma}

\begin{proof}
The identity \eqref{DomainVariation} follows from the admissible domain variations $u(y) \mapsto u \left ( y+ \epsilon \varphi(y) \right )$,
whereas \eqref{TargetVariation} follows from the admissible target variations $u(y) \mapsto \left ( 1 + \epsilon \psi(y) \right ) u(y)$.
\end{proof}

\begin{lemma} \label{ClassicalPseudoMonotonicity}
Suppose $u \, : \, \Omega \to \Sigma_N$ is stationary with respect to the variations \eqref{CombinedVariations} of functional \eqref{Functional}, 
$\Omega \subset \mathbf{R}^n$ and $x \in \Omega$ satisfies $\mathrm{dist} (x, \partial \Omega ) > r$.
Then the following identities hold: \begin{equation}
D(x,r) 
= \int_{\partial B_r(x)} \left \langle \frac{\partial u}{\partial \nu} , u \right \rangle \mathrm{d}\sigma 
+ \int_{B_r(x)} \sum_{k=1}^N \lambda_k \left | u_k \right |^2 \mathrm{d}y.
\label{DirichletAlternative}
\end{equation}
\begin{equation}
\begin{aligned}
\frac{d}{dr} \left ( \frac{D(x,r)}{r^{n-2}} \right ) = 
\frac{2}{r^{n-2}} \int_{\partial B_r(x)} \left | \frac{\partial u}{\partial \nu} \right |^2 \mathrm{d}\sigma & +
\frac{1}{r^{n-2}} \int_{\partial B_r(x)} \sum_{k=1}^N \lambda_k \left | u_k \right |^2 \mathrm{d}\sigma \\ & -
\frac{n}{r^{n-1}} \int_{B_r(x)} \sum_{k=1}^N \lambda_k \left | u_k \right |^2 \mathrm{d}y.
\end{aligned}
\label{RenormalizedEnergy}
\end{equation}
\begin{equation}
\frac{d}{dr} \left ( \frac{H(x,r)}{r^{n-1}} \right ) 
= \frac{2 D(x,r)}{r^{n-1}}  
- \frac{2}{r^{n-1}} \int_{B_r(x)} \sum_{k=1}^N \lambda_k \left | u_k \right |^2 \mathrm{d}y = \frac{2 F(x,r)}{r^{n-1}}.
\label{RenormalizedHeight}
\end{equation}
\end{lemma}

\begin{proof}
Note that the absolute continuity of $D(x,r)$ and $H(x,r)$ with respect to $r$ on $\left ( 0, \infty \right )$ follows by the argument of \cite[A.1.2(5)]{A2000}.
Using \eqref{TargetVariation} with $\psi_\epsilon$ smoothly approximating the characteristic function of $B_r(x)$ and passing to the limit $\epsilon \downarrow 0$ gives the identity \eqref{DirichletAlternative}.
The identity \eqref{RenormalizedEnergy} follows from considering \eqref{DomainVariation} with the vector field $\varphi_\epsilon (y) = \phi_\epsilon (y) (y-x)$,
where $\phi_\epsilon$ are smooth approximations of the characteristic function of $B_r(x)$, and passing to the limit $\epsilon \downarrow 0$.
Finally, the formula \eqref{RenormalizedHeight} follows from a direct calculation using \eqref{DirichletAlternative}.
\end{proof}

\begin{lemma} \label{BoundaryPoincareCorollary}
Suppose $u \, : \, \Omega \to \Sigma_N$ is stationary with respect to the variations \eqref{CombinedVariations} of functional \eqref{Functional}, 
$\overline{r} \leq \sqrt{\frac{n-1}{2 \lambda_M}}$ is fixed, where $\lambda_M = \max_{1 \leq k \leq N} \lambda_k$, and $x \in \Omega$ satisfies $\mathrm{dist} (x, \partial \Omega ) > \overline{r}$.
Then for every $r \in \left ( 0, \overline{r} \right )$, we have
\begin{equation}
\frac{F(x,r)}{r^{n-2}} + \frac{H(x,r)}{r^{n-1}} \geq \frac{1}{2} \left ( \frac{D(x,r)}{r^{n-2}} + \frac{H(x,r)}{r^{n-1}} \right ) \geq 0.
\label{BoundaryPoincareComparison}
\end{equation}
\end{lemma}
\begin{proof}
For every $v \in H^1_{\mathrm{loc}} \left ( \mathbf{R}^n \right )$ and every $r > 0$, using Fubini's theorem, integration by parts and Cauchy-Schwarz inequality, it is easy to verify that 
the following Poincar\'{e}-type inequality holds
\begin{equation}
\int_{B_r(x)} |v|^2 \mathrm{d}y \leq \frac{1}{n-1} \left [ r^2 \int_{B_r(x)} | \nabla v |^2 \mathrm{d}y + r \int_{\partial B_r(x)} |v|^2 \mathrm{d}\sigma \right ].
\label{BoundaryPoincare}
\end{equation}
Consequently, we have the lower bound
\begin{equation*}
\frac{F(x,r)}{r^{n-2}} \geq \left ( 1 -  \frac{\lambda_M r^2 }{n-1} \right ) \frac{D(x,r)}{r^{n-2}} - \left ( \frac{\lambda_M r^2}{n-1} \right ) \frac{H(x,r)}{r^{n-1}}.
\end{equation*}
\eqref{BoundaryPoincareCorollary} follows immediately from this lower bound for any $r \in \left ( 0, \min \left \{ R_x, \sqrt{\frac{n-1}{2 \lambda_M}} \right \} \right )$.
\end{proof}

\begin{lemma} \label{FirstAlmostMonotonicityFormulas}
Suppose $u \, : \, \Omega \to \Sigma_N$ is as in Lemma \ref{LocalLipschitzStationary}, 
$\overline{r} \leq \sqrt{\frac{n-1}{2 \lambda_M}}$ is fixed, and $x \in \Omega$ satisfies $\mathrm{dist} (x, \partial \Omega ) > \overline{r}$.
Then for every $r \in \left ( 0, \overline{r} \right )$, $H(x,r) \neq 0$.
Furthermore, for almost every $r \in \left ( 0, \overline{r} \right )$, we have
\begin{equation}
\frac{d}{dr} \log G (x,r) \geq - \frac{4 \lambda_M}{n-1} r,
\label{ClassicalFreqScalarDeriv}
\end{equation}
where $\lambda_M = \max_{1 \leq k \leq N} \lambda_k$. 
In particular, $e^{\frac{2 \lambda_M}{n-1} r^2} G(x,r)$ is monotone nondecreasing in $r$ on the interval $\left ( 0, \overline{r} \right )$, 
and for almost every $r \in \left ( 0, \overline{r} \right )$, we have
\begin{equation}
\frac{d}{dr} \left ( \log \frac{H(x,r)}{r^{n-1}} \right ) = \frac{2}{r} \left ( G(x,r) - 1 \right ).
\label{ClassicalLogDerivHeight}
\end{equation}
\end{lemma}

\begin{proof}
Firstly, assume $H(x,r) \neq 0$. Using \eqref{RenormalizedHeight} we immediately obtain
\begin{equation*}
\begin{aligned}
\frac{d}{dr} \left ( \log \frac{H(x,r)}{r^{n-1}} \right ) 
& = \frac{r^{n-1}}{H(x,r)} \frac{d}{dr} \left ( \frac{H(x,r)}{r^{n-1}} \right ) \\
& = \frac{2}{r} \left ( \frac{rF(x,r)}{H(x,r)} \right ) = \frac{2}{r} \left ( G(x,r) - 1 \right ).
\end{aligned}
\end{equation*}

Next, we calculate
\begin{equation}
\frac{d}{dr} \log G (x,r) = 
\frac{\frac{d}{dr} \left ( \frac{F(x,r)}{r^{n-2}} \right ) \frac{H(x,r)}{r^{n-1}} 
- \frac{F(x,r)}{r^{n-2}} \frac{d}{dr} \left ( \frac{H(x,r)}{r^{n-1}} \right )}{\left ( \frac{F(x,r)}{r^{n-2}} + \frac{H(x,r)}{r^{n-1}} \right ) \frac{H(x,r)}{r^{n-1}}}.
\label{BasicStep}
\end{equation}
We note that by \eqref{RenormalizedEnergy}, 
\begin{equation}
\frac{d}{dr} \left ( \frac{ F(x,r) }{r^{n-2}} \right ) = 
\frac{2}{r^{n-2}} \int_{\partial B_r(x)} \left | \frac{\partial u}{\partial \nu} \right |^2 \mathrm{d}\sigma - \frac{2}{r^{n-1}} \int_{B_r(x)} \sum_{k=1}^N \lambda_k \left | u_k \right |^2 \mathrm{d}y,
\label{LocalFunctionalDeriv}
\end{equation}
while from \eqref{RenormalizedHeight} and \eqref{DirichletAlternative} we get
\begin{equation*}
\begin{aligned}
& \frac{F(x,r)}{r^{n-2}} = \frac{1}{r^{n-2}} \int_{\partial B_r(x)} \left \langle \frac{\partial u }{\partial \nu} , u \right \rangle \mathrm{d}\sigma, \quad \mathrm{and} \\
& \frac{d}{dr} \left ( \frac{H(x,r)}{r^{n-1}} \right ) = \frac{2 F(x,r)}{r^{n-1}} = \frac{2}{r^{n-1}} \int_{\partial B_r(x)} \left \langle \frac{\partial u }{\partial \nu} , u \right \rangle \mathrm{d}\sigma.
\end{aligned}
\end{equation*}
Therefore, applying the Cauchy-Schwarz inequality to the second term in the numerator of \eqref{BasicStep}, we obtain
\begin{equation*}
\frac{d}{dr} \log G (x,r) \geq - \frac{\frac{2}{r^{n-1}} \int_{B_r(x)} \sum_{k=1}^N \lambda_k \left | u_k \right |^2 \mathrm{d}y}{\frac{F(x,r)}{r^{n-2}} + \frac{H(x,r)}{r^{n-1}}},
\end{equation*}
where the equality is realized, if and only if $u$ is a homogeneous map.

Using \eqref{BoundaryPoincare} and \eqref{BoundaryPoincareComparison}, for $r \in \left ( 0, \overline{r} \right )$, we estimate
\begin{equation}
\begin{aligned}
\frac{2}{r^{n-1}} \int_{B_r(x)} \sum_{k=1}^N \lambda_k \left | u_k \right |^2 \mathrm{d}y 
& \leq \frac{2 \lambda_M r}{n-1} \left [ \frac{D(x,r)}{r^{n-2}} + \frac{H(x,r)}{r^{n-1}} \right ] \\
& \leq \frac{4 \lambda_M r }{n-1} \left [ \frac{F(x,r)}{r^{n-2}} + \frac{H(x,r)}{r^{n-1}} \right ].
\label{BallBound}
\end{aligned}
\end{equation}
Consequently, we have
\begin{equation*}
\frac{d}{dr} \log G (x,r) \geq - \frac{4 \lambda_M}{n-1} r,
\end{equation*}
and $e^{\frac{2 \lambda_M}{n-1} r^2} G(x,r)$ is monotone nondecreasing in $r$ on the interval $ \left ( 0, \overline{r} \right )$.

Finally, we justify the assumption $H(x,r) \neq 0$ for $r \in \left ( 0, \overline{r} \right )$.
Firstly, we show that the interior of $u^{-1} \{0\}$ is empty.
Suppose that the interior of $u^{-1} \{0\}$ is nonempty. Then there exists an $x \in u^{-1} \{0\}$ and an $\overline{r} > 0$ such that $B_r(x) \subset u^{-1} \{0\}$.
By Lemma \ref{LocalLipschitzStationary}, $d = \mathrm{dist} \left (x, \partial u^{-1} \{0\} \right )$ is well defined, 
and there exists a sufficiently small $\epsilon > 0$ such that $H(x,r) > 0$ for $\left ( d, d + \epsilon \right )$.
Without loss of generality we can assume that $d < r_0$.
Hence, by \eqref{ClassicalLogDerivHeight}, $\frac{H(x,r)}{r^{n-1}}$ solves the initial value problem,
\begin{equation*}
\begin{aligned}
\frac{d}{dr} \left ( \frac{H(x,r)}{r^{n-1}} \right ) & = \frac{2}{r} \left ( G(x,r) - 1 \right ) \frac{H(x,r)}{r^{n-1}}, \quad r \in \left ( d, d  + \epsilon \right ), \\
\frac{H(x,d)}{d^{n-1}} & = 0.
\end{aligned}
\end{equation*}
Note that $\frac{2}{r} \left ( G(x,r) - 1 \right )$ is absolutely continuous on $\left [ d, d  + \epsilon \right )$,  as $e^{\frac{2 \lambda_M}{n-1} r^2} G(x,r)$ is monotone on this interval.
Hence, the solution of the above initial value problem is unique, and therefore $H(x,r) \equiv 0$ on $\left ( d, d  + \epsilon \right )$, which is a contradiction.

Now suppose that $H \left (x,r_0 \right ) = 0$ for some $r_0 \in \left ( 0, \overline{r} \right )$.
Using \eqref{RenormalizedHeight} and \eqref{BoundaryPoincare}, we see that for almost every $r \in \left ( 0, \overline{r} \right )$ we have
\begin{equation*}
\frac{d}{dr} \left ( \frac{H(x,r)}{r^{n-1}} \right ) + \left ( \frac{\lambda_M}{n-1} \right ) 2r \frac{ H(x,r)}{r^{n-1}} \geq \frac{D(x,r)}{r^{n-1}},
\end{equation*}
and therefore,
\begin{equation}
\frac{d}{dr} \left ( e^{\frac{\lambda_M}{n-1} r^2} \frac{H(x,r)}{r^{n-1}} \right ) \geq e^{\frac{\lambda_M}{n-1} r^2 } \frac{D(x,r)}{r^{n-1}} \geq 0. 
\label{AlmostMonotoneHeight}
\end{equation}
As $H(x,r)$ is absolutely continuous, $H \left ( x, r_0 \right ) = 0$ and \eqref{AlmostMonotoneHeight} together imply $u \equiv 0$ in $B_{r_0}(x)$, contradicting that $u^{-1} \{0\}$ has empty interior.
\end{proof}

\subsection{Homogeneous blowups and the gap condition} \label{DecompositionSection}

We review the existence of homogeneous blowups and the frequency gap condition in a slightly more general context than it was originally stated in \cite{CL07}.
Our motivation in doing so is that we will need Hausdorff dimension estimates on the nodal and critical sets of limits of solutions to \eqref{EquivalentProblem} for our compactness arguments.

\begin{proposition} \label{HomogeneousBlowupsGeneral}
Suppose $u \, : \, \Omega \to \Sigma_N$ is a map that is stationary with respect to the variations \eqref{CombinedVariations} of functional \eqref{Functional},
belongs to $C^{0, \alpha} \left ( \Omega \right )$ for every $\alpha \in (0,1)$, and satisfies the Pohozaev identity \eqref{RenormalizedEnergy}.
Then up to a subsequence, $u_{x, \rho_i} (y) = c \left ( \rho_i \right ) u \left ( x + \rho_i y \right )$, where $ c ( \rho_i ) =  \rho_i^{-n/2} \left \| u \right \|_{L^2 \left ( B_{\rho_i} (x) \right )}$, 
converge strongly in $H^1_{loc} \cap C^{0, \alpha}_{loc} \left ( \mathbf{R}^n \right )$.  
Moreover, there exists a $\delta_n > 0$ such that for every $x \in u^{-1} \{ 0 \}$, 
either $I \left ( x, 0^+ \right ) = 1$, or $I \left ( x, 0^+ \right ) \geq 1 + \delta_n$.
\end{proposition}

\begin{proof}
When $u$ is locally Lipschitz continuous in $\Omega$, this is a special case of the result in \cite[Sections 4-6]{TT12}, 
which involves deriving uniform Lipschitz bounds for blowup sequences and analyzing the homogeneous blowup maps.
Moreover, taking Lemma \ref{LocalLipschitzStationary} into account, the result extends to stationary maps in $C^{0, \alpha} \left ( \Omega \right )$ for every $\alpha \in (0,1)$,
satisfying the Pohozaev identity \eqref{RenormalizedEnergy}.
\end{proof}

A first consequence of Proposition \ref{HomogeneousBlowupsGeneral} is the following Hausdorff dimension estimate on the nodal and critical sets of H\"older continuous stationary maps.

\begin{corollary} \label{DimensionReduction}
Suppose $u \, : \, \Omega \to \Sigma_N$ is a stationary map with respect to the variations \eqref{CombinedVariations} of functional \eqref{Functional},
belongs to $C^{0, \alpha} \left ( \Omega \right )$ for every $\alpha \in (0,1)$, and satisfies the Pohozaev identity \eqref{RenormalizedEnergy}.
Then the following hold:
\begin{enumerate}[(i)]
 \item $u^{-1} \{0\}$ has Hausdorff dimension $(n-1)$.
 \item $I \left ( x, 0^+ \right ) = 1$ for $x \in u^{-1} \{0\} \backslash \mathcal{S}(u)$, and
 $$ \mathcal{S}(u) = \left \{ x \in u^{-1} \{0\} \, : \, I \left ( x, 0^+ \right ) \geq 1 + \delta_n \right \}.$$
 \item $\mathcal{S}(u)$ is relatively closed in $u^{-1} \{0\}$ and has Hausdorff dimension $(n-2)$ at most.
\end{enumerate}
\end{corollary}

\begin{proof}
The result follows from Lemma \ref{LocalLipschitzStationary}, Proposition \ref{HomogeneousBlowupsGeneral} and Federer's dimension reduction principle, 
cf. \cite[Appendix A]{Simon}. See \cite[Sections 4-6]{TT12} for details.
\end{proof}

A second consequence of Proposition \ref{HomogeneousBlowupsGeneral} is a strengthening of the monotonicity result in Lemma \ref{ClassicalPseudoMonotonicity} at the points in $u^{-1} \{0\}$.
 
\begin{corollary} \label{LinMonotonicity}
Suppose $u \, : \, \Omega \to \Sigma_N$ is a Lipschitz continuous map that is stationary with respect to the variations \eqref{CombinedVariations} of functional \eqref{Functional}, 
$\tilde{r} \leq \sqrt{ \frac{ \log (4/3) (n-1)}{2 \lambda_M}}$ is fixed, and $x \in u^{-1} \{0\}$ satisfies $\mathrm{dist} (x, \partial \Omega ) > \tilde{r}$.
Then there exists a positive constant $\Lambda = \Lambda \left ( n, N, \lambda_M \right )$ such that
$e^{\Lambda r^2} I(x,r)$ is a monotone nondecreasing function of $r$ on the interval $\left ( 0, \tilde{r} \right )$.
\end{corollary}

\begin{proof}
For simplicity, assume that
Combining \eqref{RenormalizedEnergy}, \eqref{RenormalizedHeight} and \eqref{DirichletAlternative}, 
we observe that
\begin{equation*}
I(x,r) = \frac{r \int_{\partial B_r(x)} \frac{\partial u}{\partial \nu} \cdot u \, \mathrm{d}\sigma }{H(x,r)} + \frac{ r \int_{B_r(x)} \sum_{k=1}^N \lambda_k \left | u_k \right |^2 \mathrm{d}y}{H(x,r)},  
\end{equation*}
and for almost every $r \in \left ( 0, \overline{r} \right )$,
\begin{equation*}
\begin{aligned}
\frac{dI}{dr}(x,r) = \frac{2r \int_{\partial B_r(x)} \left | \frac{\partial u}{\partial \nu} \right |^2 \mathrm{d}\sigma}{H(x,r)} 
& - \frac{2 I(x,r) \int_{\partial B_r(x)} \frac{\partial u}{\partial \nu} \cdot u \, \mathrm{d}\sigma}{H(x,r)} \\
& + \frac{r^{n+1}}{H(x,r)} \frac{d}{dr} \left ( \frac{1}{r^n} \int_{B_r(x)} \sum_{k=1}^N \lambda_k \left | u_k \right |^2 \mathrm{d}y \right ),
\end{aligned}
\end{equation*}
where $\overline{r}$ is as defined in Lemma \ref{FirstAlmostMonotonicityFormulas}.
Plugging in the first identity in the second yields for almost every $r \in \left ( 0, \overline{r} \right )$,
\begin{equation*}
\begin{aligned}
\frac{dI}{dr}(x,r) = 
\frac{2r}{H(x,r)^2} \left [ H(x,r) \int_{\partial B_r(x)} \left | \frac{\partial u}{\partial \nu} \right |^2 \, \mathrm{d}\sigma 
- \left ( \int_{\partial B_r(x)} \frac{\partial u}{\partial \nu} \cdot u \, \mathrm{d}\sigma \right )^2 \right ] \\
- 2r \left ( \frac{\int_{B_r(x)} \sum_{k=1}^N \lambda_k \left | u_k \right |^2 \mathrm{d}y}{H(x,r)}  \right ) \frac{F(x,r)}{H(x,r)} 
+ \frac{r^{n+1}}{H(x,r)} \frac{d}{dr} \left ( \frac{1}{r^n} \int_{B_r(x)} \sum_{k=1}^N \lambda_k \left | u_k \right |^2 \mathrm{d}y \right ).
\end{aligned}
\end{equation*}
Since the first term on the right-hand side is nonnegative by the Cauchy-Schwarz inequality, and $F(x,r) \leq D(x,r)$, we have
\begin{equation*}
\frac{dI}{dr}(x,r) \geq - 2 I(x,r) \left ( \frac{ \lambda_M \int_{B_r(x)} |u|^2 \mathrm{d}y}{H(x,r)} \right ) - n \frac{ \lambda_M \int_{B_r(x)} |u|^2 \mathrm{d}y}{H(x,r)} + \frac{r \lambda_m H(x,r)}{H(x,r)},
\end{equation*}
where $\lambda_m = \min_{1 \leq k \leq K} \lambda_k$,
Noting that the last term on the right-hand side is positive and applying \eqref{SpheresControlBalls} with $\alpha = 0$ to the first and second terms on the right-hand side gives 
that for almost every $r \in \left ( 0, \overline{r} \right )$,
\begin{equation}
\frac{dI}{dr}(x,r) \geq  - C_1 \left ( \lambda_M, n, r_0 \right ) r I(x,r) - \frac{C_2 \left ( n, N, \lambda_M \right )}{I(x,r)} r I(x,r).
\label{FirstLowerBound}
\end{equation}

From the monotonicity of $e^{Cr^2} G(x,r)$, we have the lower bound
\begin{equation*}
I(x,r) \geq G(x,r) - 1 \geq e^{- Cr^2} G \left (x, 0^+ \right ) - 1 = e^{- Cr^2} \left ( I \left (x, 0^+ \right ) + 1 \right ) - 1.
\end{equation*}
Since $I \left ( x, 0^+ \right ) \geq 1$ by Proposition \ref{HomogeneousBlowupsGeneral}, 
letting $r' = \min \left \{ \overline{r}, \sqrt{ \log(4/3)/C } \right \}$, 
we have $I(x,r) \geq 1/2$  for $r \in \left ( 0, r' \right )$.
In fact, recalling that $C = \sqrt{2 \lambda_M /(n-1)}$ and observing that $\log (4/3) < 1$, 
we have $r' = \min \left \{ R_x, \sqrt{ \frac{ \log (4/3) (n-1)}{2 \lambda_M}} \right \} = \tilde{r}$.
Combining $I(x,r) \geq 1/2$ with \eqref{FirstLowerBound} yields for almost every $r \in \left ( 0, \tilde{r} \right )$,
\begin{equation*}
 \frac{dI}{dr}(x,r) \geq  - 2 \Lambda r I(x,r),
\end{equation*}
where $\Lambda =  \left ( C_1 + 2 C_2 \right )/2$. Note that $\Lambda$ depends on $\lambda_M$ and $n$ only. 
Hence, we conclude that $e^{\Lambda r^2 } I(x,r)$ is monotone nondecreasing in $r$ on the interval $r \in \left ( 0, \tilde{r} \right )$.
\end{proof}

\begin{remark} \label{UsefulFreqMax}
We note that as for $u$ as in Corollary \ref{LinMonotonicity}, $x \mapsto D(x,r)$ and $x \mapsto H(x,r)$ are continuous functions.
Moreover, observing that $\tilde{r} < \overline{r}$, as defined in Lemma \ref{ClassicalPseudoMonotonicity} and Corollary \ref{LinMonotonicity} respectively,
by Lemma \ref{ClassicalPseudoMonotonicity}, $H \left (x, \tilde{r} \right ) \neq 0$.
Hence, for every compact $K \subset \Omega$ such that the Hausdorff distance between $K$ and $\partial \Omega$, $\mathrm{dist}_{\mathcal{H}} \left ( K, \partial \Omega \right ) > \tilde{r}$,
we have constants $C = C \left ( K, \tilde{r} \right )$, $c = c \left ( K, \tilde{r} \right ) > 0$ such that $0 \leq c \leq H(x,r) \leq C$ for every $x \in K$.
Thus, $I \left ( x , \tilde{r} \right )$ is a continuous function of $x$ over $K$ as well. Furthermore,
\begin{equation}
I_{K,\tilde{r}} = \sup_{x \in K} I \left ( x, \tilde{r} \right ) \leq \sup_{x \in K} e^{\frac{2 \lambda_M}{n-1} \tilde{r}^2} \left [ I \left ( x, \tilde{r} \right ) + 1  \right ].
\label{MaximalFreq}
\end{equation}
By Lemma \ref{FirstAlmostMonotonicityFormulas}, $I_{K,\tilde{r}}$ is majorized by a function monotone increasing in $\tilde{r}$ and continuous in $x \in K$. 
In particular, upper bounds involving constants depending on $I_{K,\tilde{r}}$ are unaffected, when we shrink $\tilde{r}$, which will always be bounded from an above by a constant depending
on $n$ and $\lambda_M$ only and monotone increasing in $\lambda_M$. 
In fact, in Remark \ref{ScaleRemark}, after introducing a new smallness requirement on $\tilde{r}$, we will also introduce a uniform bound $\mathcal{I}$ for the right-hand side of \eqref{MaximalFreq},
and consequently all our estimates will depend on the constant $\mathcal{I}$.
\end{remark}

\subsection{Weiss monotonicity formula} \label{WeissSub}

We will need a version of the classical monotonicity formula due to Weiss \cite{Weiss}.

\begin{proposition} \label{ClassicalWeissMonotonicity}
Suppose $u \, : \, \Omega \to \Sigma_N$ is a locally Lipschitz continuous map that is stationary with respect to the variations \eqref{CombinedVariations} of functional \eqref{Functional},
$\overline{r} \leq \sqrt{\frac{n-1}{2 \lambda_M}}$ is fixed, $x \in \Omega$ satisfies $\mathrm{dist} (x, \partial \Omega ) > \overline{r}$, and $\alpha = I \left ( x, 0^+ \right )$. 
Then there exists a constant $\mathcal{E}$ depending on $\lambda_M$, $N$, $\alpha$, and $n$ only, 
such that for almost every $r \in \left (0, \frac{1}{2} \overline{r} \right )$, 
the following inequality holds:
\begin{equation}
\frac{d}{dr} \left ( \frac{D(x,r)}{r^{n-2+2\alpha}} -  \frac{ \alpha H(x,r)}{r^{n-1+2\alpha}} + \mathcal{E} r^2 \right ) \geq
\frac{2}{r^{n+2\alpha}} \int_{\partial B_r(x)} \left | (y-x) \cdot \nabla u - \alpha u \right |^2 \mathrm{d}\sigma(y).
\label{WeissMonotonicity}
\end{equation}
\end{proposition}

\begin{proof}
For almost every $r \in \left ( 0, \mathrm{dist} (x, \partial \Omega )  \right )$, we have
\begin{equation*}
\begin{aligned}
\frac{d}{dr} \left ( \frac{D(x,r)}{r^{n-2+2\alpha}} \right ) & = 
\frac{2}{r^{n-2+2\alpha}} \int_{\partial B_r(x)} \left | \frac{\partial u}{\partial \nu} \right |^2 \mathrm{d}\sigma
+ \frac{1}{r^{n-2+2\alpha}} \int_{\partial B_r(x)} \sum_{k=1}^N \lambda_k \left | u_k \right |^2 \mathrm{d}\sigma \\
& - \frac{n}{r^{n-1+2\alpha}}\int_{ B_r(x)} \sum_{k=1}^N \lambda_k \left | u_k \right |^2 \mathrm{d}y 
- \frac{2 \alpha D(x,r)  }{r^{n-1 + 2 \alpha}} ,
\end{aligned}
\end{equation*}
and
\begin{equation*}
\begin{aligned}
\frac{d}{dr} \left ( \frac{ H(x,r)}{r^{n-1+2\alpha}} \right ) & = 
\frac{2 D(x,r) }{r^{n-1+2\alpha}} 
- \frac{2}{r^{n-1+2\alpha}} \int_{B_r(x)} \sum_{k=1}^N  \lambda_k \left | u_k \right |^2 \mathrm{d}y
- \frac{2\alpha}{r^{n+2\alpha}} H(x,r). 
\end{aligned}
\end{equation*}
Therefore, for almost every $r \in \left ( 0, \mathrm{dist} (x, \partial \Omega ) \right )$,
\begin{equation*}
\begin{aligned}
\frac{d}{dr} & \left ( \frac{D(x,r)}{r^{n-2+2\alpha}} - \frac{ \alpha H(x,r)}{r^{n-1+2\alpha}} \right ) = \\
& \left [ \frac{2}{r^{n-2+2\alpha}} \int_{\partial B_r(x)} \left | \frac{\partial u}{\partial \nu} \right |^2 \mathrm{d}\sigma 
- \frac{4\alpha D(x,r)}{r^{n-1+2\alpha}} + \frac{2\alpha^2 H(x,r)}{r^{n+2\alpha}} \right ] + \\
&  \frac{1}{r^{n-2+2\alpha}} \int_{\partial B_r(x)} \sum_{k=1}^N \lambda_k \left | u_k \right |^2  \mathrm{d}\sigma
+ \frac{2 \alpha}{r^{n-1+2\alpha}} \int_{B_r(x)} \sum_{k=1}^N \lambda_k \left | u_k \right |^2 \mathrm{d}y \\
& - \frac{n}{r^{n-1+2\alpha}} \int_{B_r(x)} \sum_{k=1}^N \lambda_k \left | u_k \right |^2 \mathrm{d}y.
\end{aligned}
\end{equation*}
Using \eqref{DirichletAlternative}, for almost every $r \in \left ( 0, \mathrm{dist} (x, \partial \Omega ) \right )$, we obtain
\begin{equation}
\begin{aligned}
\frac{d}{dr} \left ( \frac{D(x,r)}{r^{n-2+2\alpha}} - \frac{ \alpha H(x,r)}{r^{n-1+2\alpha}} \right ) = 
\frac{2}{r^{n+2\alpha}} \int_{\partial B_r(x)} \left | (y-x) \cdot \nabla u - \alpha u \right |^2 \mathrm{d}\sigma(y) \\
 + \left [ \frac{1}{r^{n-2+2\alpha}} \int_{\partial B_r(x)} \sum_{k=1}^N \lambda_k \left | u_k \right |^2 \mathrm{d}\sigma 
- \frac{2\alpha + n}{r^{n-1+2\alpha}} \int_{B_r(x)} \sum_{k=1}^N \lambda_k \left | u_k \right |^2 \mathrm{d}y \right ].
\end{aligned}
\label{AlmostWeiss}
\end{equation}
Note that the term in brackets on the right-hand side is precisely
$$
r^2 \frac{d}{dr} \left [ \frac{1}{r^{n+2\alpha}} \int_{B_r(x)} \sum_{k=1}^N \lambda_k \left |u_k \right |^2 \mathrm{d}y \right ]. 
$$
It is easy to observe that the right-hand side is zero, if $u$ is a homogeneous function of $|y-x|$ with degree of homogeneity $\alpha$ in $B_r(x)$.

Next, we observe that
\begin{equation}
\frac{d}{dr} \left ( \frac{H(x,r)}{r^{n-1+2\alpha}} \right ) = \frac{2 H(x,r)}{r^{n+2\alpha}} \left ( G(x,r) - \left ( 1 + \alpha \right ) \right ),
\label{SharpRenormalizedHeight}
\end{equation}
and consequently, for almost every $r \in \left ( 0, \overline{r} \right )$, we have
\begin{equation}
\frac{d}{dr} \log \left ( \frac{H(x,r)}{r^{n-1+2\alpha}} \right ) = \frac{2}{r} \left ( G(x,r) - (1 + \alpha ) \right ). \label{SharpLogHeight}
\end{equation}
Since $e^{\frac{2 \lambda_M}{n-1} r^2} G(x,r)$ is monotone nondecreasing in $r$ on the interval $\left ( 0, \overline{r} \right )$ and 
$$ 1 + \alpha = \lim_{r \to 0} \left [ e^{\frac{2 \lambda_M}{n-1} r^2} G(x,r) \right ], $$
for almost every $r \in \left ( 0, \overline{r} \right )$, we have
\begin{equation}
G(x,r) - \left ( 1 + \alpha \right ) \geq - \left ( 1 + \alpha \right ) \left ( 1 - e^{\frac{2 \lambda_M}{n-1} r^2} \right ) \geq - C \left ( \alpha, n, \lambda_M \right ) r^2.
\label{GeneralizedFreqLB}
\end{equation}
Hence, for almost every $r \in \left ( 0, \overline{r} \right )$, we obtain
\begin{equation*}
\frac{d}{dr} \log \left ( \frac{H(x,r)}{r^{n-1+2\alpha}} \right ) \geq - 2 C r,
\end{equation*}
which implies the monotonicity of $e^{Cr^2} \frac{H(x,r)}{r^{n-1+2\alpha}}$ on $\left ( 0, \overline{r} \right )$.

Using the monotonicity formula for $e^{Cr^2} \frac{H(x,r)}{r^{n-1+2\alpha}}$,
we can bound $\frac{H \left (x, \overline{r}/2 \right ) }{\overline{r}^{n-1+2\alpha}}$ by a constant depending on $n$, $N$, $\lambda_M$, $\alpha$ and $R_x$,
whenever $I \left ( x, 0^+ \right ) = \alpha$.
We use the crude estimate
\begin{equation*}
\begin{aligned}
N = \int_{\Omega} |u|^2 \mathrm{d}y \geq \int_{B_{\overline{r}}(x)} |u|^2 \mathrm{d}y 
& \geq \int_{\overline{r}/2}^{\overline{r}} H(x,r) \mathrm{d}r \\
& \geq e^{- \frac{3}{4}C \overline{r}^2} \frac{H \left ( x, \overline{r}/2 \right )}{ \left ( \overline{r}/ 2 \right )^{n-1+2\alpha}} \int_{\overline{r}/2}^{\overline{r}} r^{n-1+2\alpha} \mathrm{d}r.
\end{aligned}
\end{equation*}
Thus, we obtain
\begin{equation}
\frac{H \left ( x, \overline{r}/2 \right )}{\overline{r}^{n-1+2\alpha}} \leq \frac{ N 2^{n-1+2\alpha} e^{\frac{3}{4} C \overline{r}^2} }{\left ( 2^{n+2\alpha} - 1 \right ) \overline{r}^{n+2\alpha}}.
\label{UniformBoundAv}
\end{equation}
Likewise, using Fubini's theorem and the monotonicity of $e^{Cr^2} \frac{H(x,r)}{r^{n-1+2\alpha}}$, we estimate
\begin{equation}
\begin{aligned}
\int_{B_r(x)} |u|^2 \mathrm{d}y 
& = \int_0^r \left ( e^{Cs^2}\frac{H(x,s)}{s^{n-1+2\alpha}} \right ) e^{-Cs^2} s^{n-1+2\alpha} \mathrm{d}s \\
& \leq e^{Cr^2} \frac{H(x,r)}{r^{n-1+2\alpha}} \int_0^r e^{-Cs^2} s^{n-1+2\alpha} \mathrm{d}s \leq e^{Cr^2} \frac{r H(x,r)}{n+2\alpha}.
\end{aligned}
\label{SpheresControlBalls}
\end{equation}
As a result, we can bound
$$
\left [ \frac{1}{r^{n-2+2\alpha}} \int_{\partial B_r(x)} \sum_{k=1}^N \lambda_k \left | u_k \right |^2 \mathrm{d}\sigma 
- \frac{2\alpha + n}{r^{n-1+2\alpha}} \int_{B_r(x)} \sum_{k=1}^N \lambda_k \left | u_k \right |^2 \mathrm{d}y \right ]
$$
from below by
\begin{equation*}
\begin{aligned}
- \left ( e^{Cr^2} \lambda_M - \lambda_m \right ) \frac{H(x,r)}{r^{n-1+2\alpha}} r 
& \geq  - \left ( e^{Cr^2} \lambda_M - \lambda_m \right ) e^{-Cr^2} e^{\frac{C}{4} \overline{r}^2} 2^{n-1+2\alpha} \frac{H \left (x,\overline{r}/2 \right )}{\overline{r}^{n-1+2\alpha}} r \\
& \geq  - \lambda_M \frac{ N e^{C \overline{r}^2} 2^{n-1+2\alpha} }{\left ( 2^{n+2\alpha} - 1 \right ) \overline{r}^{n+2\alpha}} r,
\end{aligned}
\end{equation*}
where the last inequality is due to \eqref{UniformBoundAv}. 
Recalling that $ C = C \left ( \alpha, n, \lambda_M \right )$, $\overline{r} = \min \left \{ R_x , \sqrt{\frac{n-1}{2 \lambda_M}}, \sqrt{\frac{j_{\frac{n}{2}-1,1}}{\lambda_M}} \right \}$,
setting $\mathcal{E} =  \frac{ N \lambda_M e^{C \overline{r}^2} 2^{n-1+2\alpha} }{2 \left ( 2^{n+2\alpha} - 1 \right ) \overline{r}^{n+2\alpha}}$ and
combining this lower bound with \eqref{AlmostWeiss}, we conclude that for almost every $r \in \left (0, \overline{r} \right )$,
\begin{equation}
\frac{d}{dr} \left ( \frac{D(x,r)}{r^{n-2+2\alpha}} -  \frac{ \alpha H(x,r)}{r^{n-1+2\alpha}} + \mathcal{E} r^2 \right ) 
\geq \frac{2}{r^{n+2\alpha}} \int_{\partial B_r(x)} \left | (y-x) \cdot \nabla u - \alpha u \right |^2 \mathrm{d}\sigma(y).
\end{equation}
\end{proof}

\subsection{Compactness of stationary maps} \label{CompactnessSubsection}

We will need the following lemma for our compactness arguments in the subsequent sections. 

\begin{lemma} \label{ClassicalCompactnessLemma}
Let $u^i \, : \, B_R(0) \to \Sigma_N$ be a sequence of Lipschitz continuous maps that are stationary with respect to the variations \eqref{CombinedVariations} of functional \eqref{Functional},
with corresponding $\lambda_1^{(i)}$, ..., $\lambda_N^{(i)}$, such that $\lambda_M^{(i)} = \max_{1 \leq k \leq N} \lambda_k^{(i)} \leq \lambda_M$,
and $R \in \left ( 0, \tilde{r} \right )$, where $\tilde{r} = \sqrt{ \frac{ \log (4/3) (n-1)}{2 \lambda_M}}$.
If the respective frequency and height functions $I^i(x,R)$ and $H^i(x,R)$ for $u^i$ satisfy the uniform bound, 
\begin{equation}
\sup_{i \geq 1} \left [ I^i(x,R) + H^i(x,R) \right ] \leq C, \label{Equiboundedness}
\end{equation}
for some $C > 0$, then up to subsequences, $\left \{ u^i \right \}$ converges strongly in $L^2 \left (B_R(0) \right )$ to a locally Lipschitz continuous map $u$ in $B_R(0)$.
Moreover, the convergence is strong in $C^{0,\alpha}_{loc} \cap H^{1}_{loc}$ for every $\alpha \in (0,1)$, and $u$ is 
a stationary map with respect to the variations \eqref{CombinedVariations} of functional \eqref{Functional} with corresponding $\lambda_1^{\infty}$, ..., $\lambda_N^\infty$.
\end{lemma}

\begin{proof}
The proof is modeled on \cite[Section 3]{TT12}, where the compactness of blowup sequences for 
Lipschitz continuous maps satisfying \eqref{ExtremalityConditionA}, \eqref{ExtremalityConditionB} and \eqref{RenormalizedEnergy} is proved. 
The key observation is that as long as $R \in \left ( 0, \tilde{r} \right )$, up to minor modifications,
the arguments of \cite[Section 3]{TT12} apply to any sequence $u^i$ of Lipschitz continuous stationary maps satisfying \eqref{Equiboundedness}.
Therefore, we give a brief sketch below.

We note that the assumptions $\lambda_M^{(i)} \leq \lambda_M $, 
$R \in \left ( 0, \tilde{r} \right )$ and \eqref{Equiboundedness} imply the uniform bound $\left \| u^i \right \|_{H^1 \left ( B_R(0) \right )} \leq \tilde{C}$,
as Corollary \ref{LinMonotonicity} and equation \eqref{AlmostMonotoneHeight} apply. 
Hence, up to a subsequence, $u^i$ converge to a limit $u$ weakly in $H^1 \left ( B_R(0) \right )$ and strongly in $L^2 \left ( B_R(0) \right )$,
while $\lambda_k^{(i)} \to \lambda_k^{\infty}$ for every $1 \leq k \leq N$.

Furthermore, $R \leq \tilde{r}$ allows us to use \eqref{SharpLogHeight} with $\alpha =1$, Lemma \ref{FirstAlmostMonotonicityFormulas} and Proposition \ref{HomogeneousBlowupsGeneral}, in order to obtain
\begin{equation}
\frac{d}{dr} \left ( \log \frac{H^i (x,r)}{r^{n+1}} \right ) \geq \frac{4}{r} \left ( e^{-Cr^2 }- 1 \right ), \label{LogHeightLowerBound} 
\end{equation}
as long as $u^i (x) =0$ and $B_r(x) \subset B_R(0)$.
Consequently, integrating \eqref{LogHeightLowerBound} and using the bound \eqref{Equiboundedness}, we conclude that
\begin{equation}
\frac{H^i (x,r)}{r^{n-1}} \leq Cr^2, \label{TowardsLipschitz}
\end{equation}
whenever $u^i (x) = 0$ and $B_r(x) \subset B_R (0)$. 
Using \eqref{TowardsLipschitz} and arguing as in \cite[Lemma 3.10]{TT12}, 
we obtain $\left \| u^i \right \|_{C^{0,1} \left ( B_{R'(0)} \right ) } \leq \tilde{C} \left ( R' \right )$ for every $R' \leq R$.

By the compactness of embedding $C^{0,1} \left ( B_{R'} (0) \right ) \hookrightarrow C^{0,\alpha} \left ( B_{R'} (0) \right )$ for every $\alpha \in (0,1)$, 
we deduce that the convergence $u^i \to u$ holds in $C^{0,\alpha}_{loc}$ for every $\alpha \in (0,1)$.  
Moreover, arguing as in \cite[Lemmas 3.7 and 3.11]{TT12}, we conclude that $u^i \to u$ in $H^1_{loc} \left ( B_R (0) \right )$. 
The stationarity of limit map $u$ follows immediately from strong convergence in $H^1_{loc} \left ( B_R (0) \right )$.
Finally, arguing as in \cite[Section 4]{N2010}, we see that the limit map $u$ is locally Lipschitz continuous in $B_R(0)$.
\end{proof}

The most crucial consequence of Lemma \ref{ClassicalCompactnessLemma} is the following corollary.

\begin{corollary} \label{LimitNodalDimension}
For $u^i$ and $u$ as in Lemma \ref{ClassicalCompactnessLemma}, either $u \equiv 0$ in $B_R(0)$, or the following hold:
\begin{enumerate}[(i)]
 \item $u^{-1} \{0\}$ has Hausdorff dimension $(n-1)$.
 \item $I \left ( x, 0^+ \right ) = 1$ for $x \in u^{-1} \{0\} \backslash \mathcal{S}(u)$, and
 $$ \mathcal{S}(u) = \left \{ x \in u^{-1} \{0\} \, : \, I \left ( x, 0^+ \right ) \geq 1 + \delta_n \right \}.$$
 \item $\mathcal{S}(u)$ is relatively closed in $u^{-1} \{0\}$ and has Hausdorff dimension $(n-2)$ at most.
\end{enumerate}
\end{corollary}

\begin{proof}
Combining Lemma \ref{ClassicalCompactnessLemma} and Corollary \ref{DimensionReduction} yields the claim.
\end{proof}

\section{Smoothed frequency function and related formulas} \label{SmoothSection}

In this section we introduce the smoothed versions of classical quantities defined in the previous section, 
compare the two versions, and prove various estimates and monotonicity formulas involving the smoothed quantities.

\subsection{Smoothed variational quantities} \label{SmoothVarSubsection}
With the purpose of obtaining differential identities that hold at every scale (as opposed to almost every scale), 
as well as directional derivatives of key quantities introduced in Definition \ref{ClassicalQuantities}, we will work with analogous smoothed quantities.
The idea of working with smoothed quantities was originally introduced in \cite{DLS} and also utilized in \cite{DLMSV} and \cite{A18}.

\begin{definition} \label{SmoothedQuantities}
Let $\phi$ be the following Lipschitz function:
$$
\phi(r) =
\begin{cases}
1, \quad & 0 \leq r \leq \frac{1}{2}, \\
2 - 2r, \quad & \frac{1}{2} \leq r \leq 1, \\
0, \quad & 1 \leq r.
\end{cases}
$$
We also let $\nu_x(y) = |y-x|^{-1} \left ( y - x \right )$, a unit vector field.

The smoothed Dirichlet energy is
\begin{equation}
D_\phi (x,r) = \int_{\mathbf{R}^n} | \nabla u(y) |^2 \phi \left ( \frac{|y-x|}{r} \right ) \, \mathrm{d}y.
\label{SmoothedDirichlet}
\end{equation}
The smoothed localized key functional is
\begin{equation}
F_\phi(x,r) = \int_{\mathbf{R}^n} \left ( | \nabla u(y) |^2 - \sum_{k=1}^N \lambda_k \left | u_k (y) \right |^2 \right ) \phi \left ( \frac{|y-x|}{r} \right ) \, \mathrm{d}y.
\label{SmoothedLocalFunctional}
\end{equation}
The smoothed height function is
\begin{equation}
H_\phi(x,r) = - \int_{\mathbf{R}^n} |u (y) |^2 |y-x|^{-1} \phi' \left ( \frac{|y-x|}{r} \right ) \, \mathrm{d}y.
\label{SmoothedHeight}
\end{equation}
The smoothed frequency function is
\begin{equation}
I_\phi(x,r) = \frac{r D_{\phi}(x,r)}{H_{\phi}(x,r)}.
\label{SmoothedFrequency}
\end{equation}
The smoothed, generalized frequency function is
\begin{equation}
G_\phi(x,r) = \frac{r F_\phi(x,r) + H_\phi (x,r)}{H_\phi(x,r)}.
\label{SmoothedGenFrequency}
\end{equation}
Finally, we define the smoothed version of the first term on right-hand side of \eqref{RenormalizedEnergy},
\begin{equation}
E_\phi(x,r) = - \int_{\mathbf{R}^n} \left | \partial_{\nu_x} u(y) \right |^2 |y-x| \phi' \left ( \frac{|x-y|}{r} \right ) \, \mathrm{d}y.
\label{SmoothedRemainder}
\end{equation}
\end{definition}

We state an analogue of Lemma \ref{FirstAlmostMonotonicityFormulas} for the smoothed quantities.

\begin{proposition} \label{SmoothedVariationalFormulas}
Let $u \, : \, \Omega \to \Sigma_N$ be a stationary map with respect to the variations \eqref{CombinedVariations} of functional \eqref{Functional}, $\Omega \subset \mathbf{R}^n$, $x \in \Omega$,
and $r \in \left ( 0, \mathrm{dist}(x, \partial \Omega) \right )$.
Then the corresponding functions $D_\phi$, $F_\phi$, $H_\phi$, $I_\phi$ and $G_\phi$ are $C^1$ in the spatial variable and the positive scales $r \in (0,\infty)$,
and the following identities hold for every $x$, $v \in \mathbf{R}^3$ and $r \in (0, \infty)$:
\begin{equation}
F_\phi(x,r) = - \frac{1}{r} \int_{\mathbf{R}^n} \phi'\left ( \frac{|y-x|}{r} \right ) \left \langle \partial_{\nu_x} u , u \right \rangle \, \mathrm{d}y,
\label{SDirAlt}
\end{equation}
\begin{equation}
\begin{aligned}
\partial_r D_\phi(x,r) & = \frac{n-2}{r} D_\phi(x,r) + \frac{2}{r^2} E_\phi(x,r)
- \frac{n}{r} \int_{\mathbf{R}^n} \sum_{k=1}^N \lambda_k \left | u_k \right |^2 \phi \left ( \frac{|y-x|}{r} \right ) \, \mathrm{d}y \\
& - \frac{1}{r^2} \int_{\mathbf{R}^n} \sum_{k=1}^N \lambda_k \left | u_k \right |^2 |y-x| \phi' \left ( \frac{|y-x|}{r} \right ) \, \mathrm{d}y,
\label{SDirScale}
\end{aligned}
\end{equation}
\begin{equation}
\partial_v F_\phi(x,r)= - \frac{2}{r} \int_{\mathbf{R}^n} \phi' \left ( \frac{|y-x|}{r} \right ) \left \langle \partial_{\nu_x} u , \partial_v u \right \rangle \, \mathrm{d}y,
\label{SDirSpace}
\end{equation}
\begin{equation}
\partial_r H_\phi (x,r) = \frac{n-1}{r} H_\phi (x,r) + 2 F_\phi (x,r) ,
\label{SHeightScale}
\end{equation}
\begin{equation}
\partial_v H_\phi (x,r) = - 2 \int_{\mathbf{R}^n} |y-x|^{-1} \phi' \left ( \frac{|y-x|}{r} \right ) \left \langle \partial_v u , u \right \rangle \, \mathrm{d}y.
\label{SHeightSpace}
\end{equation}
\end{proposition}

\begin{proof}
Firstly, we note that we can approximate $\phi$ by smooth functions $\phi_k$, in which case the smoothed quantities are clearly smooth in $x$ and $r$,
and pass to a limit in $W^{1,p}$ for every $p < \infty$ to obtain uniform convergence to the above differential identities. Hence, we can treat $\phi$ as smooth.

Using \eqref{TargetVariation} with $\psi(y) = \phi \left ( \frac{|x-y|}{r} \right ) u(y)$ gives \eqref{SDirAlt},
while differentiating $D_\phi (x,r)$ with respect to $r$ and testing \eqref{DomainVariation} with $\varphi(y) = \phi \left ( \frac{|x-y|}{r} \right ) (y-x)$ yields \eqref{SDirScale}.
Likewise, differentiating $D_\phi (x,r)$ with respect to $x$ and testing \eqref{DomainVariation} with $\varphi(y) = \phi \left ( \frac{|x-y|}{r} \right ) v$ gives \eqref{SDirSpace}.
A direct calculation analogous to \eqref{RenormalizedHeight}, using \eqref{SDirAlt}, gives \eqref{SHeightScale}. 
Finally, changing variables, $y = x + r \omega$, differentiating in the direction of $v$, and changing back to the variable $y$, we obtain \eqref{SHeightSpace}.
\end{proof}

Using Lemma \ref{SmoothedVariationalFormulas}, we easily obtain a smoothed version of Lemma \ref{ClassicalPseudoMonotonicity}.

\begin{lemma} \label{SmoothedLogLemma}
Let $u$ be as in Lemma \ref{SmoothedVariationalFormulas}. 
Then there exists an $\overline{r} = \overline{r} \left (n , \lambda_M \right )$, where $\lambda_M = \max_{1 \leq k \leq N} \lambda_k$,
such that $H_\phi(x,r) \neq 0$ for $r \in \left ( 0, \overline{r} \right )$, whenever $B_r(x) \subset \Omega$. 
Furthermore, for almost every $r \in \left ( 0, \overline{r} \right )$, we have
\begin{equation}
\frac{d}{dr} \log G_\phi (x,r) = - C \left ( n, \lambda_M \right ) r.
\label{SFreqScalarDeriv}
\end{equation}
whenever $B_r(x) \subset \Omega$.
In particular, $e^{C \left ( n, \lambda_M \right ) r^2} G_\phi(x,r)$ is monotone nondecreasing in $r$ on the interval $ \left ( 0, \overline{r} \right )$, and for almost every $r \in \left ( 0, \overline{r} \right )$,
\begin{equation}
\frac{d}{dr} \left ( \log \frac{ H_\phi (x,r) }{r^{n-1}} \right ) = \frac{2}{r} \left ( G_\phi(x,r) - 1 \right ), 
\label{SLogDerivHeight}
\end{equation}
whenever $B_r(x) \subset \Omega$.
\end{lemma}

\begin{proof}
The proof proceeds analogously to the proof of Lemma \ref{ClassicalPseudoMonotonicity}.
We only point out that by a direct calculation, we have
\begin{equation}
\frac{2}{r} \int_{B_r(x) \backslash B_{r/2}(x)} |u|^2 \, \mathrm{d}y  \leq H_\phi(x,r) \leq \frac{4}{r} \int_{B_r(x) \backslash B_{r/2}(x)} |u|^2 \, \mathrm{d}y.
\label{UsefulComparison}
\end{equation}
Furthermore, using \eqref{RenormalizedHeight} and \eqref{BoundaryPoincareComparison}, for $r \in \left ( 0, \overline{r} \right )$, we observe that \eqref{UsefulComparison} implies
\begin{equation}
\frac{4}{n-1} H \left ( x, r/2 \right ) \leq H_\phi(x,r) \leq \frac{4}{n-1} H(x,r).
\label{UsefulComparison2}
\end{equation}
Using \eqref{UsefulComparison}, \eqref{UsefulComparison2} and arguments in the same spirit, and following the proof of Lemma \ref{ClassicalPseudoMonotonicity}, we obtain \eqref{SFreqScalarDeriv}. 
Finally, $H_\phi(x,r) \neq 0$ on $\left ( 0, \overline{r} \right )$ and \eqref{SLogDerivHeight} for almost every $r \in \left ( 0, \overline{r} \right )$ follow as in Lemma \ref{ClassicalPseudoMonotonicity}.
\end{proof}

\begin{remark} \label{ScaleRemark}
For the rest of the article, we assume that the number $\tilde{r} > 0$ satisfies
\begin{equation}
\tilde{r} \leq \min \left \{ \overline{r}, \sqrt{ \frac{ \log (4/3) (n-1)}{2 \lambda_M}} \right \} = R \left ( n, N, \lambda_M \right ),
\label{ScaleAssumption}
\end{equation}
where $\overline{r} = \overline{r} \left (n , \lambda_M \right )$ is as defined in Lemma \ref{SmoothedLogLemma}, and the second part of the above restriction is due to Corollary \ref{LinMonotonicity}.
Without loss of generality we can assume that $\tilde{r} < 1$. 
In Section \ref{MinkowskiSection} we will shrink $\tilde{r}$ further, if necessary. Therefore, we will make sure that the constants in upper bounds we derive are either uniform
or monotone nondecreasing in $\tilde{r}$.
Finally, for the rest of this article, we assume that maps $u$ for which the right-hand side of \eqref{MaximalFreq} is less than a constant we denote as $\mathcal{I}$.
Consequently all our estimates will depend on the constant $\mathcal{I}$.
\end{remark}

\subsection{Comparison with classical frequency} \label{ComparisonSubsection}
We would like to have an analogue of Corollary \ref{LinMonotonicity} for the smoothed frequency as well. 
However, for the proof we need a lower bound on $I_\phi (x,r)$ whenever $x \in u^{-1} \{0\}$.
Therefore, as a first step we show that the classical and smoothed frequencies are comparable, up to a shift in scale.

\begin{lemma} \label{CSComparisonLemma}
Let the map $u$ be as in Lemma \ref{SmoothedVariationalFormulas}, $\tilde{r} \in \left ( 0, R \left ( n, N, \lambda_M \right ) \right )$, 
and $K \subset \Omega$ compact with $\mathrm{dist}_{\mathcal{H}} \left ( K, \Omega \right ) > \tilde{r}$.
Then there exists a positive constant $C = C \left (n, \lambda_M, \mathcal{I} \right )$ such that for every $x \in u^{-1} \{0\} \cap K$ and every $r \in \left ( 0, \tilde{r} \right )$, we have
\begin{equation}
C I(x,r) \geq I_\phi (x,r) \geq C^{-1} I \left ( x, r/2 \right ),
\label{ComparingCS}
\end{equation}
Furthermore, there exists a positive constant $\epsilon = \epsilon \left (n, \lambda_M, K \right )$ such that for every $x \in u^{-1} \{0\} \cap K$ and $r \in \left ( 0, \tilde{r} \right )$,
\begin{equation}
I_\phi \left ( x, r \right ) > \epsilon.
\label{SFreqLowerBound}
\end{equation}
\end{lemma}

\begin{proof}
Note that \eqref{UsefulComparison2} direclty implies
\begin{equation}
\frac{n-1}{2} \left ( \frac{H(x,r/2)}{H(x,r)} \right ) I(x,r/2) \leq  I_\phi(x,r) \leq \frac{n-1}{4} \left ( \frac{H(x,r)}{H(x,r/2)} \right ) I(x,r).
\label{TowardsComparison}
\end{equation}
Moreover, using \eqref{RenormalizedHeight} and Corollary \ref{LinMonotonicity}, we can estimate
\begin{equation}
\log \frac{H(x,r)}{H(x,r/2)} \leq 2 I \left ( x, \tilde{r} \right ) e^{\Lambda \tilde{r}^2} \log 2,
\label{HeightRatio1}
\end{equation}
and consequently, for $\mathcal{I}$ as defined in \eqref{ScaleRemark}, we have.
\begin{equation}
\frac{H(x,r)}{H(x,r/2)} \leq C \left ( \mathcal{I}, \tilde{r} \right ).
\label{HeightRatio2}
\end{equation}
Therefore, \eqref{TowardsComparison} and \eqref{HeightRatio2} together yield  \eqref{ComparingCS}.

Next, we prove \eqref{SFreqLowerBound}. Using \eqref{ComparingCS} and \eqref{LinMonotonicity}, we observe that
\begin{equation*}
e^{\Lambda \tilde{r}^2 / 4} I_\phi (x,r) >  e^{\Lambda r^2 / 4} I_\phi (x,r) \geq C^{-1} e^{\Lambda r^2 / 4} I(x,r/2) \geq C^{-1} I \left ( x, 0^+ \right ).
\end{equation*}
Since $I \left ( x, 0^+ \right ) \geq 1$ by Proposition \ref{HomogeneousBlowupsGeneral}, letting $\epsilon = C^{-1} e^{- \Lambda \tilde{r}^2 / 4}$ yields \eqref{SFreqLowerBound}. 

Finally, note that $\tilde{r}$ and $\Lambda$ in Corollary \ref{LinMonotonicity} depend only on $n$ and $\lambda_M$, 
as long as $x \in u^{-1} \{0\} \cap K$ and $\mathrm{dist}_{\mathcal{H}} \left ( K, \Omega \right ) > \tilde{r}$.
Hence, $C$ and $\epsilon$ depend on $\mathcal{I}$, $n$, $\lambda_M$ only.
\end{proof}

\subsection{Monotonicity formulas for smooothed frequency at zeros} \label{SmoothLinSubsection}
Modifying the smoothed frequency with a multiplicative exponential factor in order to obtain a monotone quantity will be useful in the following lemmas.

\begin{corollary} \label{LinMonotonicitySmoothed}
Let $u$ be as in Corollary \ref{LinMonotonicity}, $\tilde{r} \in \left ( 0, R \left ( n, N, \lambda_M \right ) \right )$, and $K \subset \Omega$ compact with $\mathrm{dist}_{\mathcal{H}} \left ( K, \Omega \right ) > \tilde{r}$.
Then there exists a $\Lambda' = \Lambda' \left (n, N, \lambda_M, \mathcal{I} \right )$ such that for every $x \in u^{-1} \{0\} \cap K$ and almost every $r \in \left ( 0, \tilde{r} \right )$, we have
\begin{equation}
\frac{d}{dr} I_\phi (x,r) + 2 \Lambda' r I_\phi (x,r) \geq \frac{2}{r H_\phi(x,r)^2} \left [ H_\phi(x,r) E_\phi (x,r) - r^2 F_\phi (x,r)^2 \right ] \geq 0. \label{SLinMonotonicityFormula}
\end{equation}
In particular, for every $x \in u^{-1} \{0\} \cap K$, $e^{\Lambda' r^2 } I_\phi (x,r)$ is a monotone nondecreasing function of $r$ on the interval $\left ( 0, \tilde{r} \right )$. 
\end{corollary}

\begin{proof}
The proof is analogous to the proof of Corollary \ref{LinMonotonicity}. 
We use \eqref{UsefulComparison} and \eqref{UsefulComparison2} in order to estimate the perturbation terms depending on $\lambda_k$ for $1 \leq k \leq N$.
Moreover, we use the lower bound \eqref{SFreqLowerBound} to obtain \eqref{SLinMonotonicityFormula} from the analogue of \eqref{FirstLowerBound}.
Finally, the nonnegativity of $H_\phi(x,r) E_\phi (x,r) - r^2 F_\phi (x,r)^2$ follows from \eqref{SDirAlt} and the Cauchy-Schwarz inequality.
\end{proof}

We also need another modified version of the smoothed frequency with an additive quadratic correction term.

\begin{corollary} \label{LinMonotonicitySmoothedAdditive}
Let $u$ be as in Corollary \ref{LinMonotonicity}, $\tilde{r} \in \left ( 0, R \left ( n, N, \lambda_M \right ) \right )$, and $K \subset \Omega$ compact with $\mathrm{dist}_{\mathcal{H}} \left ( K, \Omega \right ) > \tilde{r}$.
Then there exists an $A = A \left (n, N, \lambda_M, \mathcal{I} \right )$ such that for the frequency with additive correction,
\begin{equation}
I_\phi^A (x,r) = I_\phi (x,r) + Ar^2, 
\label{SmoothedFrequencyAdditive}
\end{equation}
for every $x \in u^{-1} \{0\} \cap K$ and almost every $r \in \left ( 0, \tilde{r} \right )$, we have
\begin{equation}
\frac{d}{dr} I_\phi^A (x,r) \geq \frac{2}{r H_\phi(x,r)^2} \left [ H_\phi(x,r) E_\phi (x,r) - r^2 F_\phi (x,r)^2 \right ] \geq 0. \label{SLinMonotonicityFormulaAdditive}
\end{equation}
\end{corollary}

\begin{proof}
The claim is immediate from \eqref{SLinMonotonicityFormula}, \eqref{ComparingCS} and Corollary \ref{LinMonotonicity} combined.
\end{proof}

\subsection{Elementary upper bounds} \label{ElementaryBoundsSubsection}
We can also compare the smoothed height and frequency functions at nearby points, as long as we adjust the scale. This is proved in the following lemma.

\begin{lemma} \label{ElementaryUpperBounds}
Let $u$ be as in Lemma \ref{SmoothedVariationalFormulas}, $\tilde{r} \in \left ( 0, R \left ( n, N, \lambda_M \right ) \right )$,
and let $K \subset \Omega$ be a compact set satisfying $\mathrm{dist}_{\mathcal{H}} (K, \partial \Omega ) > \tilde{r}$.
Then there exists a constant $C_1 = C_1 \left (n, \lambda_M \right )$  such that whenever $\rho \in \left ( 0, \frac{\tilde{r}}{4} \right )$ and $y \in B_{\rho}(x)$, we have
\begin{equation}
H_\phi (y,\rho) \leq C_1 H_\phi (x,4 \rho).
\label{ShiftedHeightComparison}
\end{equation}
Furthermore, there exist constants $C_2 = C_2 \left ( n, N, \lambda_M, K \right )$ such that whenever, $r \in \left ( 0, \frac{\tilde{r}}{4} \right )$, $x \in K$ and $y \in B_{r/4}(x)$, we have
\begin{equation}
I_\phi(y,r) \leq C_2 I_\phi (x, 4r).
\label{ShiftedFrequencyComparison}
\end{equation}
\end{lemma}

\begin{proof}
We prove \eqref{ShiftedHeightComparison} first. Since $y \in B_\rho(x)$ implies $B_\rho(y) \subset B_{2\rho}(x)$, we observe that
\begin{equation}
H_\phi (y, \rho) = 2 \int_{\rho/2}^\rho \frac{H(y,s)}{s} \, \mathrm{d}s \leq \frac{4}{\rho} \int_0^\rho H(y,s) \, \mathrm{d}s \leq \frac{4}{\rho} \int_0^{2\rho} H(x,s) \, \mathrm{d}s.
\label{SHCStep1}
\end{equation}
Using \eqref{LogHeightLowerBound}, for every $r \in \left ( 2 \rho, \tilde{r} \right )$, we derive the estimate,
\begin{equation}
\frac{4}{\rho} \int_0^{2\rho} H(x,s) \, \mathrm{d}s 
\leq \frac{4}{\rho} \int_0^{2\rho} e^{Cr^2} H(x,r) \left ( \frac{s}{r} \right )^{n+1} \, \mathrm{d}s 
\leq \frac{2^{n+3}}{n+2} e^{Cr^2} \left ( \frac{\rho}{r} \right )^{n+1}  H(x,r).
\label{SHCStep2}
\end{equation}
From \eqref{SHCStep1} and \eqref{SHCStep2} we obtain
\begin{equation}
H_\phi (y, \rho) \leq \frac{2^{n+3}}{n+2} e^{Cr^2} \left ( \frac{\rho}{r} \right )^{n+1}  H(x,r), \label{SHCStep3}
\end{equation}
where $C = C \left ( \lambda_M, n \right )$ as in Lemma \eqref{ClassicalPseudoMonotonicity}.
Now multiplying both sides of \eqref{SHCStep3} by $r^n$ and integrating with respect to $r$ on the interval $\left ( 2 \rho , 4 \rho \right )$, we get
\begin{equation}
\frac{(4 \rho)^{n+1}}{n+1} \left ( 1 - \frac{1}{2^{n+1}} \right ) H_\phi (y, \rho ) \leq \frac{2^{n+2}}{n+2} e^{C \tilde{r}^2} \rho^{n+1} H_\phi (x, 4 \rho). \label{SHCStep4}
\end{equation}
Since $C$ and $\tilde{r}$ both depend on $\lambda_M$ and $n$ only, letting $C_1 = e^{C \tilde{r}^2} \frac{n+1}{n+2} \cdot \frac{2}{2^{n+1}-1}$, \eqref{ShiftedHeightComparison} follows from \eqref{SHCStep4}.

Now we prove \eqref{ShiftedFrequencyComparison}.
Since $y \in B_{r/4}(x)$, we have $B_r(y) \subset B_{2r}(x)$, and therefore,
\begin{equation}
I_\phi(y,r) = \frac{r D_\phi (y,r) }{H_{\phi(y,r)}} \leq \frac{1}{4} \frac{ (4r) D_\phi (x, 4r)}{H_\phi(y,r)}.
\label{SFCStep1} 
\end{equation}
Likewise, $y \in B_{r/4}(x)$ implies $x \in B_{r/4}(y)$, and hence, using \eqref{ShiftedHeightComparison}, we obtain
\begin{equation}
H_{\phi}(y,r) \geq \frac{1}{C_1} H \left ( x, \frac{r}{4} \right ).
\label{SFCStep2}
\end{equation}
By \eqref{SLogDerivHeight}, we also have
\begin{equation}
\frac{H_\phi (x,r/4)}{(r/4)^{n-1}} = \frac{H_\phi (x,4r)}{(4r)^{n-1}} \exp \left [ - \int_{r/4}^{4r} \frac{2}{s} \left ( G_\phi(x,s) - 1\right ) \, \mathrm{d}s \right ].
\label{SFCStep3}
\end{equation}
From \eqref{SFCStep1}, \eqref{SFCStep2} and \eqref{SFCStep3}, we obtain
\begin{equation}
I_\phi(y,r) \leq 4^{2n-3} C \exp \left [ \int_{r/4}^{4r} \frac{2}{s} \left ( G_\phi(x,s) - 1\right ) \, \mathrm{d}s \right ] I_\phi (x,4r).
\label{SFCStep4}
\end{equation}
Since $4r \leq \tilde{r}$, using \eqref{SFreqScalarDeriv}, \eqref{ComparingCS} and Remark \ref{UsefulFreqMax}, we observe that for every $s \in (r/4,4r)$, we have the estimate,
\begin{equation}
\begin{aligned}
& G_\phi(x,s) - 1 
\leq e^{Cs^2} G_\phi (x,s) - 1 
\leq e^{C \tilde{r}^2} G_\phi \left (x, \tilde{r} \right ) - 1 
\leq e^{C \tilde{r}^2} I_\phi \left (x, \tilde{r} \right ) + \left ( e^{C \tilde{r}^2} - 1 \right ) \\
& \leq C e^{C \tilde{r}^2} I \left (x, \tilde{r} \right ) + C' \tilde{r}^2
\leq C \left ( \mathcal{I}, n, \lambda_M \right ).
\label{SFCStep5}
\end{aligned}
\end{equation}
Finally, \eqref{SFCStep4} and \eqref{SFCStep5} together yield \eqref{ShiftedFrequencyComparison}, as $\tilde{r}$ depends on $n$ and $\lambda_M$ only.
\end{proof}

Lemma \ref{ElementaryUpperBounds} allows us to prove an $\epsilon$-clearing result for $u^{-1} \{0\} \cap K$.

\begin{corollary} \label{EpsilonClearing}
Let $u$ be as in Lemma \ref{SmoothedVariationalFormulas}, $\tilde{r} \in \left ( 0, R \left ( n, N, \lambda_M \right ) \right )$,
Also let $K \subset \Omega$ be a compact set satisfying $\mathrm{dist}_{\mathcal{H}} (K, \partial \Omega ) > \tilde{r}$, and $x \in K$.
Then there exists an $\epsilon' = \epsilon' \left ( n, N, \lambda_M, K \right )$ such that if
\begin{equation}
I_\phi (x,r ) < \epsilon', 
\label{EpsilonCondition}
\end{equation}
for some $r \in \left ( 0, \tilde{r} \right )$, then
\begin{equation}
u^{-1} \{0\} \cap B_{r/16}(x) = \emptyset.
\label{Clear}
\end{equation}
\end{corollary}

\begin{proof}
Suppose $y \in u^{-1} \{0\} \cap B_{r/16}(x)$. Then from Lemmas \ref{CSComparisonLemma} and \ref{ElementaryUpperBounds}, we immediately have
\begin{equation}
\epsilon  < I_\phi (y, r/4) \leq C_2 I_\phi (x,r),
\end{equation}
and therefore, $I_\phi(x,r) > \epsilon / C_2$. Letting $\epsilon' = \epsilon / C_2$ yields the claim.
\end{proof}

\section{Main frequency estimates} \label{MainFreqEstimatesSection}

In this section we analyze the oscillations of smoothed frequency function with respect to scale or space.
In particular, we would like to show when and how the former can (almost) control the latter.

\subsection{Frequency pinching} \label{FreqPinchingSubsection}
Alongside $I_\phi^A (x,r)$, the smoothed frequency with an additive quadratic modification defined in \eqref{SmoothedFrequencyAdditive},
a key quantity for our analysis is the corresponding frequency pinching.

\begin{definition} \label{PinchingDefined} 
We define $W_{s,t}^{A}(x)$, the $A$-frequency pinching at $x$ between scales $s$ and $t$ as
\begin{equation}
W_{s,t}^{A} (x) = I_\phi^{A} (x,t) - I_\phi^{A} (x,s),
\label{PinchingCorrected}
\end{equation} 
where $0 \leq s \leq t < \tilde{r}$, $\tilde{r} \leq \sqrt{ \frac{ \log (4/3) (n-1)}{2 \lambda_M}}$, $I_\phi^{A}(r) = I_\phi(x,r) + Ar^2$. 
\end{definition}

Note that by Corollary \ref{LinMonotonicitySmoothedAdditive}, $W_{s,t}^A (x) \geq 0$, 
whenever $x \in u^{-1} \{0\}$ and $A$ is greater than or equal to the constant $A = A \left (n, N, \lambda_M, \mathcal{I} \right )$ in Corollary \ref{LinMonotonicitySmoothedAdditive}.
Hence, for the rest of this article, we fix $A$ as in Corollary \ref{LinMonotonicitySmoothedAdditive} and modify $W_{s,t}^A (x)$, when necessary.

\subsection{A refined Weiss-type monotonicity formula} \label{WeissSubsection}
As a first step, we state and prove an estimate in the spirit of \eqref{WeissMonotonicity}.
This estimate is a generalization of the Weiss-type estimate introduced in \cite{DLMSV}, different versions of which also played analogous roles in \cite{A18}, \cite{FS} and \cite{FS2}.

\begin{lemma} \label{WeissTypeLemma}
Suppose $u \, : \, \Omega \to \Sigma_N$ is a Lipschitz continuous map that is stationary with respect to the variations \eqref{CombinedVariations} of functional \eqref{Functional}, 
$\tilde{r} \in \left ( 0, R \left ( n, N, \lambda_M \right ) \right )$,  and $K \subset \Omega$ is compact with $\mathrm{dist}_{\mathcal{H}} \left ( K, \Omega \right ) > \tilde{r}$. 
Then there exists a constant $C = C \left ( n, N, \lambda_M, \mathcal{I} \right )$ such that 
for every $x \in u^{-1} \{0\} \cap K$
and for every $R$, $r$ such that $r \in \left [ R/8, R \right ]$ and $R \in \left ( 0, \tilde{r}/2 \right )$, 
on the annulus $A_{r,R}(x) = B_R(x) \backslash B_r(x)$, we have
\begin{equation}
\int_{A_{r,R}(x)} \left | (y-x) \cdot \nabla u(y) - I_\phi \left (x, |x-y| \right ) u(y) \right |^2 \, \mathrm{d}y
\leq C R H_\phi(x, 2R) W_{r/2,2R}^{1+A+A^2} (x).
\label{WeissTypeEstimate}
\end{equation}
\end{lemma}

\begin{proof}
By switching to polar coordinates, using the positivity of integrand and invoking Fubini's theorem, we obtain
\begin{equation}
\begin{aligned}
\int_{A_{r,R}(x)} & \left | (y-x) \cdot \nabla u - I_\phi \left (x, |x-y| \right ) u \right |^2 \, \mathrm{d}y \\
& \leq \frac{1}{r} \int_{r}^{2R} \int_{A_{\frac{\tau}{2},\tau}(x)} \left | (y-x) \cdot \nabla u(y) - I_\phi \left (x, |x-y| \right ) u \right |^2 \, \mathrm{d}y \, \mathrm{d}\tau,
\label{WTLStep1}
\end{aligned} 
\end{equation}
and therefore,
\begin{equation}
\begin{aligned}
\int_{A_{r,R}(x)} & \left | (y-x) \cdot \nabla u - I_\phi \left (x, |x-y| \right ) u \right |^2 \, \mathrm{d}y \\
& \leq \frac{2}{r} \int_{r}^{2R} \int_{A_{\frac{\tau}{2},\tau}(x)} \left | (y-x) \cdot \nabla u - I_\phi \left (x, \tau \right ) u \right |^2 \, \mathrm{d}y \, \mathrm{d}\tau \\
& + \frac{2}{r} \int_{r}^{2R} \int_{A_{\frac{\tau}{2},\tau}(x)} \left | I_\phi \left (x, \tau \right ) - I_\phi \left (x, |x-y| \right ) \right |^2 |u|^2  \, \mathrm{d}y \, \mathrm{d}\tau.
\label{WTLStep2}
\end{aligned}
\end{equation}
We estimate the first and second terms on the right-hand side of \eqref{WTLStep2} separately.

We observe that first term on the right-hand side of \eqref{WTLStep2} is bounded from above by
\begin{equation}
\begin{aligned}
\frac{2}{r} \int_r^{2R} \frac{\tau}{2} \int_{\mathbf{R}^n} \frac{-1}{|x-y|} \phi' \left ( \frac{|y-x|}{\tau} \right ) \left | (y-x) \cdot \nabla u (y) - I_\phi(x,\tau) u \right |^2 \, \mathrm{d}y \, \mathrm{d}\tau \\
\leq 16 \int_r^{2R} \left [ E_\phi (x, \tau) - 2 \tau I_\phi (x,\tau) F_\phi(x,\tau) + I_\phi(x,\tau)^2 H_\phi(x,\tau) \right ] \, \mathrm{d} \tau,
\end{aligned}
\label{FirstTermStep1}
\end{equation}
where the inequality follows from $r > R/8$, \eqref{SDirAlt}, and the definitions of $E_\phi(x,\tau)$ and $H_\phi(x,\tau)$.
Furthermore, we have
\begin{equation}
\begin{aligned}
\int_r^{2R} & \left [ E_\phi (x, \tau) - 2 \tau I_\phi (x,\tau) F_\phi(x,\tau) + I_\phi(x,\tau)^2 H_\phi(x,\tau) \right ] \, \mathrm{d} \tau = \\
& (I) + (II) = \int_r^{2R} \left [ E_\phi (x, \tau) - \tau^2 \frac{F_\phi(x,\tau)^2}{H_\phi(x,\tau)} \right ] \, \mathrm{d} \tau + \int_r^{2R} \tau^2 \frac{P_\phi(x,\tau)^2}{H_\phi(x,\tau)} \, \mathrm{d}\tau,
\end{aligned}
\label{FirstTermStep2}
\end{equation}
where
\begin{equation}
P_\phi(x,\tau) =  \int_{\mathbf{R}^n} \sum_{k=1}^N \lambda_k \left | u_k (y) \right |^2 \phi \left ( \frac{|y-x|}{\tau} \right ) \, \mathrm{d}y. \label{SPerturbation}
\end{equation}

We note that by \eqref{SLinMonotonicityFormulaAdditive},
\begin{equation}
\begin{aligned}
E_\phi (x, \tau) - \tau^2 \frac{F_\phi(x,\tau)^2}{H_\phi(x,\tau)} & = \frac{\tau H_\phi (x,\tau)}{2} \cdot \frac{2 \left [ H_\phi(x,\tau) E_\phi (x,\tau) - \tau^2 F_\phi (x,\tau)^2 \right ] }{\tau H_\phi(x,\tau)^2} \\
& \leq \frac{\tau H_\phi (x,\tau)}{2} \partial_\tau I_\phi^A (x,\tau).
\label{FirstTermStep3}
\end{aligned}
\end{equation}
Hence, using \eqref{SLogDerivHeight}, \eqref{SLinMonotonicityFormulaAdditive} and \eqref{FirstTermStep3}, we obtain
\begin{equation}
(I) \leq C_1 R H_\phi(x,2R) \left [ I_\phi^{A} (x,2R) - I_\phi^{A} (x,r/2) \right ],
\label{FirstTermStep4}
\end{equation}
where $C_1$ is an increasing function of $\tilde{r}$, and therefore depends only on $n$ and $\lambda_M$. 

Secondly, from \eqref{ClassicalLogDerivHeight} and \eqref{HeightRatio2},
\begin{equation*}
\frac{P_\phi(x,\tau)}{H_\phi(x,\tau)} 
\leq \lambda_M \frac{\int_0^\tau H(s) \, \mathrm{d}s}{ \frac{2}{\tau} \int_{\tau/2}^\tau H(s) \, \mathrm{d}s} 
\leq \lambda_M \tau  \frac{H(\tau)}{H(\tau/2)}
\leq C \left ( \mathcal{I} \right ) \lambda_M \tau.
\end{equation*}
Hence, we have the estimate,
\begin{equation}
\begin{aligned}
(II) 
& = \int_r^{2R} \tau^2 \left ( \frac{P_\phi(x,\tau)}{H_\phi(x,\tau)} \right )^2 H_\phi(x,\tau) \, \mathrm{d}\tau \\
& \leq C \left ( \mathcal{I} \right )^2 \lambda_M^2  \int_r^{2R} \tau^4 H_\phi(x,\tau) \mathrm{d}\tau \\
& \leq C_2 R H_\phi(x,2R) \left [ (2R)^2 - (r/2)^2 \right ],
\label{FirstTermStep5}
\end{aligned}
\end{equation}
where $C_2$ depends on $n$, $\lambda_M$, and $\mathcal{I}$, as $2R \leq \tilde{r}$.

Letting $C' = \max \left \{ C_1, C_2 \right \}$, from \eqref{FirstTermStep4} and \eqref{FirstTermStep5}, we obtain
\begin{equation}
(I) + (II) \leq C' R H_\phi(x,2R) \left [  I_\phi^{A+1} (x,2R) - I_\phi^{A+1} (x,r/2) \right ].
\label{FirstTermComplete}
\end{equation}

Next we observe that the second term on the right-hand side of \eqref{WTLStep2} is bounded from above by
\begin{equation}
\begin{aligned}
(III) + (IV) = \frac{2}{r} \int_{r}^{2R} \int_{A_{\frac{\tau}{2},\tau}(x)} 
2 \left | I_{\phi}^A \left (x, \tau \right ) - I_{\phi}^A \left (x, |x-y| \right ) \right |^2 |u|^2 \, \mathrm{d}y \, \mathrm{d}\tau \\
+ \frac{2}{r} \int_{r}^{2R} \int_{A_{\frac{\tau}{2},\tau}(x)} 2 \left | A \tau^2 - A|x-y|^2 \right |^2 |u|^2  \, \mathrm{d}y \, \mathrm{d}\tau. 
\end{aligned}
\label{SecondTermStep1}
\end{equation}
By \eqref{SLogDerivHeight}, Corollary \ref{LinMonotonicitySmoothedAdditive} and $r \geq R/8$,
\begin{equation}
\begin{aligned}
(III) 
& \leq 32 \left [ I_\phi^A (x,2R) - I_\phi^A (x,r/2) \right ] \int_r^{2R} H_\phi(x,\tau) \, \mathrm{d}\tau \\
& \leq 64R H_\phi(x,2R) \left [ I_\phi^A (x,2R) - I_\phi^A (x,r/2) \right ].
\label{SecondTermStep2}
\end{aligned}
\end{equation}
Likewise, we have
\begin{equation}
(IV) \leq C \left ( \tilde{r} \right ) R H(x,2R)  \left [ A^2 (2R)^2 - A^2 (r/2)^2 \right ]
\label{SecondTermStep3}
\end{equation}
Hence, letting $C'' = \max \left \{ 64, C \left ( \tilde{r}, A \right ) \right \}$, we get
\begin{equation}
(III) + (IV) \leq C'' R H(x,2R) \left [ I_\phi^{A + A^2} (x,2R) - _\phi^{A + A^2} (x,r/2) \right ].
\label{SecondTermComplete}
\end{equation}

Finally, letting $C = \max \left \{ 16C', C'' \right \}$, which depends on $n$, $\lambda_M$ and $\mathcal{I}$ only, from \eqref{WTLStep2}, \eqref{FirstTermComplete} and \eqref{SecondTermComplete}
we arrive at \eqref{WeissTypeEstimate}.
\end{proof}

\subsection{Frequency oscillation estimate} \label{FOESubsection}
Now we are ready to prove a crucial estimate on the oscillation of frequency. This estimate is a generalization of \cite[Theorem 4.2]{DLMSV}.

\begin{proposition} \label{FrequencyOscillationProposition}
Suppose $u \, : \, \Omega \to \Sigma_N$ is a Lipschitz continuous map that is stationary with respect to the variations \eqref{CombinedVariations} of functional \eqref{Functional}, 
$\tilde{r} \in \left ( 0, R \left ( n, N, \lambda_M \right ) \right )$,  and $K \subset \Omega$ is compact with $\mathrm{dist}_{\mathcal{H}} \left ( K, \Omega \right ) > \tilde{r}$.  
Then there exists a constant $C = C \left ( n, N, \lambda_M, \mathcal{I} \right )$ such that 
for every $r \in \left ( 0, \tilde{r}/4 \right )$ and for every $x_1$, $x_2 \in u^{-1} \{0\} \cap K$ with $\left | x_1 - x_2 \right | \leq r/4$, we have
\begin{equation}
\left | I_\phi(y,r) - I_\phi(z,r) \right | \leq C \left [ \left ( W_{r/8,4r}^{2+A+A^2} \left ( x_1 \right ) \right )^{1/2} + \left ( W_{r/8, 4r}^{2+A+A^2} \left ( x_2 \right ) \right )^{1/2}\right ] \frac{|y-z|}{r}.
\label{FrequencyOscillationEstimate}
\end{equation}
whenever $y$ and $z$ lie on the line segment joining $x_1$ and $x_2$.
\end{proposition}

\begin{proof}
Let $v = x_2 - x_1$, $\mathrm{d}\mu_x = 2|y-x|^{-1} \mathbbm{1}_{[r/2,r]} \left ( |y-x| \right ) \, \mathrm{d}y$. Our goal is to estimate $ \left | \partial_v I_\phi(x,r) \right |$ uniformly,
whenever $x$ lies on the line segment joining $x_1$ and $x_2$.

We have
\begin{equation}
\partial_v I_\phi(x,r) = H_\phi(x,r)^{-1} \left [ r \partial_v F_\phi(x,r) - I_\phi(x,r) \partial_v H_\phi(x,r) \right ] + \frac{r \partial_v P_\phi(x,r)}{H_\phi(x,r)},
\label{FreqDirectionalDeriv}
\end{equation}
where $P_\phi(x,r)$ is as defined in \eqref{SPerturbation}, and
\begin{equation}
\partial_v P_\phi(x,r)  = - \frac{2}{r} \int_{A_{r/2,r}(x)} \sum_{k=1}^N \lambda_k \left | u_k \right |^2 v \cdot \nu_x(y) \, \mathrm{d}y.
\label{PerturbationDirDeriv}
\end{equation}

Firstly, we can estimate the second term on the right-hand side of \eqref{FreqDirectionalDeriv}, as
\begin{equation}
\left | \frac{r \partial_v P_\phi(x,r)}{H_\phi(x,r)} \right | 
= \left | - 2 \frac{\int_{A_{r/2,r}(x)} \sum_{k=1}^N \lambda_k \left | u_k \right |^2 v \cdot \nu_x(y) \, \mathrm{d}y}{\int_{A_{r/2,r}(x)} \frac{2}{|x-y|} \left | u \right |^2 \, \mathrm{d}y}  \right |
\leq \lambda_M r |v| \leq \frac{\lambda_M}{4} r^2.
\label{ErrorTermControlled}
\end{equation}
since $ |v| = \left | x_2 - x_1 \right | \leq r/4$. 

There remains to estimate the first term on the right-hand side of \eqref{FreqDirectionalDeriv}.
Using \eqref{SDirSpace} and \eqref{SHeightSpace}, we express the term inside the brackets as
\begin{equation}
\mathcal{M} = \int_{\mathbf{R}^n} 2 \left [ \left \langle (y-x) \cdot \nabla u, \partial_v u \right \rangle - I_\phi(x,r) \left \langle u, \partial_v u \right \rangle \right ] \, \mathrm{d}\mu_x. 
\label{FOEStep1}
\end{equation}
We write $\partial_v u$ as
\begin{equation}
\partial_v u = \left ( y - x_1 \right ) \cdot \nabla u - \left ( y - x_2 \right ) \cdot \nabla u = \mathcal{E}_1(z) + \mathcal{E}_2(z) + \mathcal{E}_3(z) u,
\label{DirDerivDecomposition}
\end{equation}
where
\begin{equation*}
\begin{aligned}
& \mathcal{E}_i(y) = \left ( y - x_i \right ) \cdot \nabla u - I_\phi \left  ( x_i , \left | y - x_i \right | \right ) u, \quad i=1,2, \\
& \mathcal{E}_3(y) = I_\phi \left ( x_1, \left | y - x_1 \right | \right ) - I_\phi \left ( x_2, \left | y - x_2 \right | \right ). 
\end{aligned}
\end{equation*}
Hence, we have $\mathcal{M} = (I) + (II) + (III)$, where
\begin{equation*}
\begin{aligned}
& (I) = 2 \int_{\mathbf{R}^n} \left ( \mathcal{E}_1 - \mathcal{E}_2 \right ) \cdot (y-x) \cdot \nabla u \, \mathrm{d}\mu_x , \\
& (II) = 2 I_\phi(x,r) \int_{\mathbf{R}^n} \left ( \mathcal{E}_2 - \mathcal{E}_1 \right ) \cdot u \, \mathrm{d}\mu_x , \\
& (III) = 2 \int_{\mathbf{R}^n} \left [ (y-x) \cdot \nabla u \cdot u - I_\phi(x,r) |u|^2 \right ] \mathcal{E}_3 \, \mathrm{d}\mu_x.
\end{aligned}
\end{equation*}

Next, we estimate $(I)$ as follows:
\begin{equation}
\begin{aligned}
(I) & \leq 4 \int_{A_{r/2,r}(x)} \left ( \left | \mathcal{E}_1 \right | + \left | \mathcal{E}_2 \right | \right ) \left | \partial_{\nu_x} u \right | \, \mathrm{d}y \\
    & \leq 4 D_{\phi}(x,2r)^{1/2} \left [ \left ( \int_{A_{r/2,r}(x)} \left | \mathcal{E}_1 \right |^2 \, \mathrm{d}y \right )^{1/2} 
    + \left ( \int_{A_{r/2,r}(x)} \left | \mathcal{E}_2 \right |^2 \, \mathrm{d}y \right )^{1/2} \right ].
\end{aligned}
\label{FOEStep2}
\end{equation}

Similarly, we estimate $(II)$ as follows:
\begin{equation}
(II) \leq 2 I_\phi (x,r) \left [ \left ( \int_{\mathbf{R}^n} \left | \mathcal{E}_1 \right |^2 \, \mathrm{d}\mu_x \right )^{1/2} 
                         + \left ( \int_{\mathbf{R}^n} \left | \mathcal{E}_2 \right |^2 \, \mathrm{d}\mu_x \right )^{1/2} \right ] H_\phi(x,r)^{1/2}.
\label{FOEStep3}
\end{equation}

We observe that $x$ lying on the line segment joining $x_1$ and $x_2$ and $\left | x_1 - x_2 \right | \leq r/4$ together imply that
$$
B_r(x) \backslash B_{r/2}(x) \subset B_{2r} \left ( x_i \right ) \backslash B_{r/4} \left ( x_i \right ), \quad i=1,2.
$$
Consequently,for $i=1,2$, by \eqref{WeissTypeEstimate}, we obtain
\begin{equation}
\int_{A_{r/2,r}(x)} \left | \mathcal{E}_i \right |^2 \, \mathrm{d}y \leq C r H_\phi \left ( x_i, 4r \right ) W_{r/8, 4r}^{1+A+A^2} \left ( x_i \right ), 
\label{FOEStep4a}
\end{equation}
and likewise,
\begin{equation}
\int_{A_{r/2,r}(x)} \left | \mathcal{E}_i \right |^2 \, \mathrm{d}\mu_x \leq C H_\phi \left ( x_i, 4r \right ) W_{r/8, 4r}^{1+A+A^2} \left ( x_i \right ).
\label{FOEStep4b}
\end{equation}
From \eqref{FOEStep2} and \eqref{FOEStep4a}, we obtain
\begin{equation}
(I) \leq 4 \sqrt{Cr}  D_{\phi}(x,2r)^{1/2} \sum_{i=1}^2 \left [ H_\phi \left ( x_i, 4r \right )^{1/2} W_{r/8, 4r}^{1+A+A^2} \left ( x_i \right )^{1/2} \right ],
\label{FOEStep2a} 
\end{equation}
and likewise, from \eqref{FOEStep3} and \eqref{FOEStep4b}, we get
\begin{equation}
(II) \leq C I_\phi (x,r) H_\phi(x,r)^{1/2} \sum_{i=1}^2 \left [  H_\phi \left ( x_i, 4r \right )^{1/2} W_{r/8, 4r}^{1+A+A^2} \left ( x_i \right )^{1/2} \right ]. 
\label{FOEStep3a}
\end{equation}

In order to estimate $(III)$, we decompose $\mathcal{E}_3(z)$ as
\begin{equation}
\mathcal{E}_3(y) =  \mathcal{E}_{3,1} + \mathcal{E}_{3,2}(y) + \mathcal{E}_{3,3}(y),
\label{ThirdTermDecomposition}
\end{equation}
where
\begin{equation*}
\begin{aligned}
& \mathcal{E}_{3,1} = I_\phi \left ( x_1, r \right ) - I_\phi \left ( x_2, r \right ), \\
& \mathcal{E}_{3,2}(y) = I_\phi \left ( x_1, \left | y - x_1 \right | \right ) - I_\phi \left ( x_1, r \right ), \\
& \mathcal{E}_{3,3}(y) = I_\phi \left ( x_2, r \right ) - I_\phi \left ( x_2, \left | y - x_2 \right | \right ). \\
\end{aligned}
\end{equation*}
As a result, we have
\begin{equation}
(III) = 2 \sum_{i=1}^3 \int_{\mathbf{R}^n} \left [ (y-x) \cdot \nabla u \cdot u - I_\phi(x,r) |u|^2 \right ] \mathcal{E}_{3,i} \, \mathrm{d}\mu_x.
\label{ThirdTermSum}
\end{equation}
We note that the first term in the sum is equal to
\begin{equation}
2 \mathcal{E}_{3,1} \int_{\mathbf{R}^n} \left [ (y-x) \cdot \nabla u \cdot u - \frac{r F_\phi(x,r)}{H_\phi(x,r)} |u|^2 \right ] \, \mathrm{d}\mu_x
- 2 \mathcal{E}_{3,1} \frac{r P_\phi(x,r)}{H_\phi(x,r)} \int_{\mathbf{R}^n} |u|^2 \, \mathrm{d} \mu_x,
\label{FOEStep5a} 
\end{equation}
and by \eqref{SDirAlt}, the first term in \eqref{FOEStep5a} is equal to zero, while we we can estimate the second term as
\begin{equation}
2 \left | \mathcal{E}_{3,1} \right | \frac{r P_\phi(x,r)}{H_\phi(x,r)} \int_{\mathbf{R}^n} |u|^2 \, \mathrm{d} \mu_x
\leq 2 \lambda_M \left [ I_\phi \left ( x_1, r \right ) + I_\phi \left ( x_2, r \right ) \right ] r \int_{B_r(x)} |u|^2 \, \mathrm{d}y.
\label{FOEStep5b}
\end{equation}
Hence, using \eqref{ComparingCS} and Remark \ref{UsefulFreqMax}, we conclude that there is a constant $C =C \left ( n, N, \lambda_M, \mathcal{I} \right )$ such that
\begin{equation}
2 \int_{\mathbf{R}^n} \left [ (y-x) \cdot \nabla u \cdot u - I_\phi(x,r) |u|^2 \right ] \mathcal{E}_{3,1} \, \mathrm{d}\mu_x
\leq C r \int_{B_r(x)} |u|^2 \, \mathrm{d}y.
\label{FOEStep5c}
\end{equation}

Next we estimate the second and third terms in the sum in \eqref{ThirdTermSum}. For $i=1,2$ we have
\begin{equation}
\begin{aligned}
2 \int_{\mathbf{R}^n} & \left [ (y-x) \cdot \nabla u \cdot u - I_\phi(x,r) |u|^2 \right ] \mathcal{E}_{3,i} \, \mathrm{d}\mu_x \\
& \leq 2 \int_{\mathbf{R}^n} \left [ |y-x| \left | \nabla u \right | |u| + I_\phi (x,r) |u|^2 \right ] \left | \mathcal{E}_{3,i} \right | \, \mathrm{d} \mu_x.
\end{aligned} 
\label{FOEStep6a}
\end{equation}
Noting that $\mu_x$ is supported on $B_r(x) \backslash B_{r/2}(x)$, where $x$ lies on the line segment joining $x_1$ and $x_2$, where $\left | x_1 - x_2 \right | \leq r/4$,
we observe that for every $y$ in the support of $\mu_x$ and for $i=1,2$,
$$
r/4 \leq \left | y - x_i \right | \leq 2r.
$$
Therefore, for every $y$ in the support of $\mu_x$ and for $i=1,2$,
\begin{equation}
\left | \mathcal{E}_{3,i}(y) \right | \leq W_{r/8, 4r}^{1+A+A^2} \left ( x_i \right ) + Cr^2,
\label{FOEStep6b}
\end{equation}
where $C$ depends on $n$, $N$, $\lambda_M$, $\mathcal{I}$.
Hence, by the H\"older inequality we obtain for $i=1,2$,
\begin{equation}
\begin{aligned}
2 \int_{\mathbf{R}^n} & \left [ (y-x) \cdot \nabla u \cdot u - I_\phi(x,r) |u|^2 \right ] \mathcal{E}_{3,i} \, \mathrm{d}\mu_x \\
& \leq 2 \left [ 2\sqrt{r} D_\phi (x,2r)^{1/2} H_\phi(x,r)^{1/2} + r D_\phi(x,r) \right ] \left ( W_{r/8,4r}^{1+A+A^2} \left ( x_i \right ) + Cr^2 \right ) . 
\end{aligned}
\label{FOEStep6c}
\end{equation}

Summing \eqref{FOEStep5c} and \eqref{FOEStep6c}, we get
\begin{equation}
\begin{aligned}
(III) \leq & \left [ 4\sqrt{r} D_\phi (x,2r)^{1/2} H_\phi(x,r)^{1/2} + 2r D_\phi(x,r) \right ] \sum_{i=1}^2 \left ( W_{r/8, 4r}^{1+A+A^2} \left ( x_i \right ) + Cr^2 \right ) \\
      & + C r \int_{B_r(x)} |u|^2 \, \mathrm{d}y.
\label{ThirdTermCompleted}
\end{aligned} 
\end{equation}

Likewise, summing \eqref{FOEStep2a}, \eqref{FOEStep3a} and \eqref{ThirdTermCompleted}, we obtain
\begin{equation}
\begin{aligned}
\mathcal{M} \leq & C r \int_{B_r(x)} |u|^2 \, \mathrm{d}y + C \sqrt{r} D_{\phi}(x,2r)^{1/2} \sum_{i=1}^2 \left [ H_\phi \left ( x_i, 4r \right )^{1/2} W_{r/8,4r}^{1+A+A^2} \left ( x_i \right )^{1/2} \right ] + \\
                 & C I_\phi (x,r) H_\phi(x,r)^{1/2} \sum_{i=1}^2 \left [  H_\phi \left ( x_i, 4r \right )^{1/2} W_{r/8,4r}^{1+A+A^2} \left ( x_i \right )^{1/2} \right ] + \\
                 & C \left [ \sqrt{r} D_\phi (x,2r)^{1/2} H_\phi(x,r)^{1/2} + r D_\phi(x,r) \right ] \sum_{i=1}^2 \left ( W_{r/8,4r}^{1+A+A^2} \left ( x_i \right ) + Cr^2 \right ),
\end{aligned}
\label{MCompleted} 
\end{equation}
for some $ C = C \left ( n, N, \lambda_M, \mathcal{I} \right )$.

We recall that
\begin{equation}
\left | \partial_v I_\phi (x,r) \right | \leq \frac{\mathcal{M}}{H_\phi(x,r)} + \frac{\lambda_M}{4} r^2.
\label{DirectionalTotal}
\end{equation}
The fact that $x$ lies on the line segment joining $x_1$ and $x_2$ with $\left | x_1 - x_2 \right | \leq r/4$ implies $x_i \in B_{r/4}(x)$, and $4r \leq \tilde{r}$.
Therefore, using \eqref{ShiftedHeightComparison}, \eqref{SLogDerivHeight} and \eqref{ComparingCS}, we observe that
\begin{equation*}
\frac{H_\phi \left ( x_i, 4r \right )}{H_\phi(x,r)} \leq C_1 \frac{H_\phi \left ( x_i, 4r \right )}{H_\phi \left ( x_i ,r/4 \right )} \leq C,
\end{equation*}
where $C_1$ is as in Lemma \ref{ElementaryUpperBounds}, and $C = C \left ( n, N, \lambda_M, \mathcal{I} \right )$.
Utilizing \eqref{SLogDerivHeight} and Remark \ref{UsefulFreqMax} in a similar manner, we estimate each term on the right-hand side of \eqref{DirectionalTotal} to get
\begin{equation}
\begin{aligned}
\left | \partial_v I_\phi (x,r) \right | 
& \leq C \left [ (4r)^2 - (r/8)^2 \right ]^{1/2} + C \sum_{i=1}^2 W_{r/8,4r}^{1+A+A^2} \left ( x_i \right )^{1/2} \\
& \leq C \sum_{i=1}^2 W_{r/8,4r}^{2+A+A^2} \left ( x_i \right )^{1/2},
\label{DirectionalComplete}
\end{aligned}
\end{equation}
for a modified constant $C = \left ( n, N, \lambda_M, \mathcal{I} \right )$.

Finally, we observe that for $y$, $z$ lying on the line segment joining $x_1$ and $x_2$,
\begin{equation}
\left | I_\phi(y,r) - I_\phi (z,r) \right | \leq  \sup_{t \in [0,1]} \left | \partial_v I_\phi \left ( (1-t)z + ty ,r \right ) \right | \cdot |y-z|.
\label{DirDerivExpressed}
\end{equation}
Hence, the claim follows from immediately \eqref{DirDerivExpressed}, \eqref{DirectionalComplete} and Remarks \ref{UsefulFreqMax} and \ref{ScaleRemark}.
\end{proof}

\section{Distortion bound} \label{JonesNumberSection}

In this section we consider nonnegative, finite Radon measures $\mu$ supported on the singular set $\mathcal{S}(u)$ for $u \, : \, \Omega \to \Sigma_N$,
Lipschitz continuous and stationary with respect to the variations \eqref{CombinedVariations} of functional \eqref{Functional}.
Our goal is to understand the geometry of the {\emph{effective}} support of such measures.

\subsection{Mean-flatness} \label{MFSubsection}
An important quantity for studying the size and structure of singular set $\mathcal{S}(u)$ is its mean-flatness, which we define next.

\begin{definition} \label{JonesNumberDefined}
For $\mu$ a nonnegative Radon measure in $\mathbf{R}^n$, $k$ a positive integer less than $n$, $x \in \mathbf{R}^n$, and $r>0$,
we define the $k$-th mean flatness of $\mu$ in $B_r(x)$ as
\begin{equation}
D_\mu^{k} (x,r) = \inf_L \frac{1}{r^{k+2}} \int_{B_r(x)} \mathrm{dist} (y,L)^2 \, \mathrm{d}\mu(y),
\label{JonesNumber}
\end{equation}
where $\mathrm{dist}(y,A) = \inf_{x \in A} |y-x|$, and the infimum in \eqref{JonesNumber} is with respect to all affine $k$-dimensional planes $L \subset \mathbf{R}^n$.
\end{definition}

\begin{remark} \label{JonesBackground}
This quantity is also known as the Jones' $\beta_2$-number, as it was originally introduced in the context of analyst's traveling salesman problem in $\mathbf{R}^2$ by Jones in \cite{Jones}.
We refer to the survey article \cite{Schul} for related problems and various generalizations.
Also see \cite{DT12, AT15, ENV} for more on the Jones's $\beta$-numbers in the context of rectifiability and bi-Lipschitz parametrizations of sets in the Euclidean space. 
\end{remark}

\subsection{Algebraic characterization} \label{LinearAlgebraSubsection}
In order to prove the main estimate of this section, firstly we need the following elementary characterization of $D_\mu^k$.
Let $x \in \mathbf{R}^n$ and $r>0$ be such that $\mu \left ( B_r ( x  )  \right ) > 0$, and define the barycenter of $\mu$ in $B_r \left ( x \right )$ as
$$
\overline{x}_{x,r} = \frac{1}{\mu \left ( B_r ( x ) \right ) } \int_{B_r ( x)} y \, \mathrm{d} \mu(y).
$$
The measure $\mu$ restricted to $B_r(x)$ induces a bilinear form $B_{x,r} \, : \, \mathbf{R}^n \times \mathbf{R}^n \to \mathbf{R}$ given by
$$ \mathcal{B}_{x,r}(v,w) = \int_{B_r(x)} \left ( \left ( y - \overline{x} \right ) \cdot v \right ) \left ( \left ( y - \overline{x} \right ) \cdot w \right ) \, \mathrm{d}\mu(y),$$
where $v$, $w \in \mathbf{R}^n$. Since $\mathcal{B}_{x,r}$ is symmetric and positive semi-definite, there is an orthonormal basis $v_1$, ..., $v_n$ for $\mathbf{R}^n$ such that
$\mathcal{B}_{x,r} \left ( v_i , v_j \right ) = \delta_{ij} \xi_i$, for nonnegative eigenvalues $\xi_i \leq \xi_j$, $n \geq i \geq j \geq 1$, where $\delta_{ij}$ is the Kronecker delta.
We observe that  and for every $i = 1$, $...$, $n$,
\begin{equation}
\int_{B_r(x)} \left ( \left ( y - \overline{x} \right ) \cdot v_i \right ) y \, \mathrm{d}\mu(y) = \xi_i v_i,
\label{Eigendirection}
\end{equation}

The $k$-th mean-flatness of $\mu$ and minimizing affine $k$-planes in \eqref{JonesNumber} can be expressed as
\begin{equation}
D_\mu^k (x,r) = r^{-k-2} \sum_{\ell=k+1}^n \xi_\ell, 
\label{JonesSum}
\end{equation}
and the infimum in the definition of $D_\mu^k$ is achieved by all the affine planes $L = x + \mathrm{Span} \left \{ v_1,..., v_k \right \}$,
for any choice of eigenbasis $v_1$, ...,$v_n$ with nonincreasing eigenvalues $\xi_1 \geq ... \geq \xi_n \geq 0$.

\subsection{Frequency pinching controls mean-flatness} \label{DistBoundSubsection} 
Finally, we state and prove an important estimate on the $(n-2)$-th mean flatness of finite, nonnegative Radon measures supported on the singular set $\mathcal{S}(u)$ of an optimal partition.

\begin{proposition} \label{JonesNumberProposition}
Suppose $u \, : \, \Omega \to \Sigma_N$ is a Lipschitz continuous map that is stationary with respect to the variations \eqref{CombinedVariations} of functional \eqref{Functional}, 
$\tilde{r} \in \left ( 0, R \left ( n, N, \lambda_M \right ) \right )$,  and $K \subset \Omega$ is compact with $\mathrm{dist}_{\mathcal{H}} \left ( K, \Omega \right ) > \tilde{r}$.  
Then there exists a constant $C = C \left ( n, N, \lambda_M, \mathcal{I} \right )$ such that for any finite, nonnegative Radon measure $\mu$ supported on $\mathcal{S}(u)$,
\begin{equation}
D_\mu^{n-2} \left ( x , r \right ) \leq \frac{C}{r^{n-2}} \int_{B_r(x)} W_{r,32r}^{3+A+A^2} (z) \, \mathrm{d}\mu(z) ,
\label{JonesNumberEstimate}
\end{equation}
for every $x \in \mathcal{S}(u) \cap K$ and for every $r \in \left ( 0, \tilde{r}/32 \right )$.
\end{proposition}

\begin{proof}
Note that if $\mu \left ( B_r(x) \right ) = 0$, the claim holds trivially. Hence, we assume that $\mu \left ( B_r(x) \right ) > 0$.
We denote the barycenter of $\mu$ in $B_r(x)$ as $\overline{x}$, and let $v_1, ..., v_n$ be an orthonormal eigenbasis for the bilinear form $B_{x,r}$ 
with respective eigenvalues $\xi_1 \geq ... \geq \xi_n \geq 0$.

{\emph{Step 1:}}
For any constant $\alpha$, \eqref{Eigendirection} implies that for every $j= 1,..,n$ and every $y \in B_{12r}(x) \backslash B_{4r}(x)$,
\begin{equation}
- \xi_j v_j \cdot \nabla u(y) = \int_{B_r(x)} \left ( \left ( z- \overline{x} \right ) \cdot v_j \right ) \left ( (y-z) \cdot \nabla u(y) - \alpha u(y) \right ) \, \mathrm{d}\mu(z),
\label{JonesStep1a}
\end{equation}
as $\overline{x}$ is the barycenter of $\mu$ in $B_r(x)$.
By squaring both sides of \eqref{JonesStep1a} and applying the Cauchy-Schwarz inequality, we get
\begin{equation*}
\begin{aligned}
\xi_j^2 \left | \partial_{v_j} u(y) \right |^2 
& \leq \left ( \int_{B_r(x)} \left | \left ( z- \overline{x} \right ) \cdot v_j \right | \left | (y-z) \cdot \nabla u(y) - \alpha u(y) \right | \, \mathrm{d}\mu(z) \right )^2 \\
& \leq \int_{B_r(x)} \left | \left ( z- \overline{x} \right ) \cdot v_j \right |^2 \, \mathrm{d}\mu(z) \int_{B_r(x)} \left | (y-z) \cdot \nabla u(y) - \alpha u(y) \right |^2 \, \mathrm{d}\mu(z).
\end{aligned}
\end{equation*}
Using $\mathcal{B}_{x,r} \left ( v_i , v_j \right ) = \delta_{ij} \xi_i$ to rewrite the first factor on the right-hand side and dividing both sides by $\xi_j \neq 0$, we have
\begin{equation}
\xi_j \left | \partial_{v_j} u(y) \right |^2 \leq \int_{B_r(x)} \left | (y-z) \cdot \nabla u(y) - \alpha u(y) \right |^2 \, \mathrm{d}\mu(z).
\label{JonesStep1b}
\end{equation}
Using \eqref{JonesSum} and $\lambda_1 \geq ... \geq \lambda_n \geq 0$, for $A(x,r) = B_{10r}(x) \backslash B_{6r}(x)$, we have
\begin{equation*}
\begin{aligned}
r^{n} D_\mu^{n-2}(x,r) \int_{A(x,r)} \sum_{j=1}^{n-1} \left | \partial_{v_j} u(y) \right |^2 \, \mathrm{d}y
& = \int_{A(x,r)} \left ( \xi_{n-1} + \xi_{n} \right ) \sum_{j=1}^{n-1} \left | \partial_{v_j} u(y) \right |^2 \, \mathrm{d}y \\
& \leq 2 \int_{A(x,r)} \sum_{j=1}^n \xi_j \left | \partial_{v_j} u(y) \right |^2 \, \mathrm{d}y.
\end{aligned}
\end{equation*}
Summing in $j=1,...,n$, integrating over $A(x,r)$ both sides of \eqref{JonesStep1b}, and invoking Fubini's theorem, we get
\begin{equation}
\begin{aligned}
r^{n} D_\mu^{n-2}(x,r) & \int_{A(x,r)} \sum_{j=1}^{n-1} \left | \partial_{v_j} u(y) \right |^2 \, \mathrm{d}y 
\leq 2 \int_{A(x,r)} \sum_{j=1}^n \xi_j \left | \partial_{v_j} u(y) \right |^2 \, \mathrm{d}y \\
& \leq 2n \int_{A(x,r)} \int_{B_r(x)} \left | (y-z) \cdot \nabla u(y) - \alpha u(y) \right |^2 \, \mathrm{d}\mu(z) \, \mathrm{d}y \\
& \leq 2n \int_{B_r(x)} \int_{B_{12r}(z) \backslash B_{4r}(z)} \left | (y-z) \cdot \nabla u(y) - \alpha u(y) \right |^2 \, \mathrm{d}y \, \mathrm{d}\mu(z).
\end{aligned}
\label{JonesStep1c}
\end{equation}

{\emph{Step 2:}}
Next we claim that there exists a constant $C_c = C_c \left ( n, N, \lambda_M, \mathcal{I} \right )$ such that
\begin{equation}
D_\phi (x,16r) \leq C_c \int_{A(x,r)} \sum_{j=1}^{n-1} \left | \partial_{v_j} u(y) \right |^2 \, \mathrm{d}y.
\label{JonesStep2Claim}
\end{equation}
If the claim were false, we would have
a sequence of Lipschitz continuous maps $u^{(i)} \, : \, K \to \Sigma_N$,
which are stationary with respect to the variations \eqref{CombinedVariations} of functional \eqref{Functional} with $\lambda_1^{(i)}$, ..., $\lambda_N^{(i)}$, 
where $\lambda_M^{(i)} = \max_{1\leq k \leq N} \lambda_k^{(i)} \leq \lambda_M$, and which satisfy $\max_{y \in K} I \left (y, \tilde{r} \right ) < \mathcal{I}$, as well as
$x_i \in \mathcal{S} \left ( u^{(i)} \right ) \cap B_{16r}(0)$, $0 \leq r_i \leq \tilde{r}/32$, such that 
for suitable scalars $c \left ( r_i \right )$ and rotations $\theta_i \, : \, \mathbf{R}^n \to \mathbf{R}^n$, the maps
$c \left ( r_i \right ) u^{(i)} \left ( x_i + r_i \theta_i(x) \right )$ relabeled as $u^{(i)}$ satisfy the following:
\begin{enumerate}[(i)]
 \item $H_\phi^{(i)} (0, 16r) = 1$, \label{HeightNormalization}
 \item $I^{(i)}_\phi \left ( 0,\tilde{r} \right ) \leq C_0 \left (n, \lambda_M, \mathcal{I} \right )$, \label{FrequencyBoundSet}
 \item $0 \in \mathcal{S} \left ( u^{(i)} \right)$, \label{JunctionPresence}
 \item while for $\nabla_T = \left ( \partial_{x_1}, ..., \partial_{x_{n-1}} \right )$,
\begin{equation}
\int_{A(0,r)} \left | \nabla_T u^{(i)} (y) \right |^2 \, \mathrm{d}y \leq \frac{1}{i} D^{(i)}_\phi (0,16r). \label{ContraBound}
\end{equation}
\end{enumerate}

By \eqref{HeightNormalization} and \eqref{FrequencyBoundSet}, the right-hand side of \eqref{ContraBound} converges to $0$, as $i \to \infty$,
while by $32r \leq \tilde{r} \leq R \left (n, N, \lambda_M \right )$ and Lemma \ref{CSComparisonLemma}, $u^{(i)}$ satisfy the hypothesis of Lemma \ref{ClassicalCompactnessLemma}.
Consequently, by a diagonalization argument, a suitable subsequence of $u^{(i)}$ converges to a map $w \, : \, B_{16r}(0) \to \Sigma_N$ 
in $L^2 \left ( B_{16r}(0) \right )$, $H^1_{loc} \left ( B_{16r}(0) \right )$ and $C^{0,\alpha}_{loc} \left ( B_{16r}(0) \right )$ for every $\alpha \in (0,1)$.
Lemma \ref{ClassicalCompactnessLemma} also ensures that limiting map $w$ is Lipschitz continuous 
and stationary with respect to the variations \eqref{CombinedVariations} of functional \eqref{Functional} 
with $\tilde{\lambda}_1$, ..., $\tilde{\lambda}_N$, and $\tilde{\lambda}_M = \max_{1\leq k \leq N} \tilde{\lambda}_k \leq \lambda_M$. 
Furthermore, $H^{w}_{\phi} (0, 16r) = 1$, and
\begin{equation}
\int_{A(0,r)} \left | \nabla_T w (y) \right |^2 \, \mathrm{d}y = 0. \label{LimitContra}
\end{equation}

The stationarity of $w$ and \eqref{LimitContra} imply that for every $1 \leq k \leq N$, $w_k$ can be expressed as
\begin{equation}
 w_k(y) = a_k \sin \left ( \tilde{\lambda}_k y_n + b_k \right ) \quad \mathrm{in} \quad \left \{ y \in A(0,r) \, : \, w_k(y) > 0 \right \} \neq \emptyset, \label{1DStationary}
\end{equation}
for constants $a_k$, $b_k$, and \eqref{HeightNormalization} and \eqref{SLogDerivHeight} ensure that at least one $a_k \neq 0$.
Note that \eqref{1DStationary} follows from $ - \Delta w_k = \tilde{\lambda}_k w_k$ in $\left \{ y \in B_{16r}(0) \, : \, w_k(y) > 0 \right \}$ for every $1 \leq k \leq N$, 
and the fact that \eqref{LimitContra} implies that in $A(0,r) = B_{10r}(0) \backslash B_{6r}(0)$, $w_k$ depends on the variable $x_n$ only.
Furthermore, each $w_k$ clearly enjoys the unique continuation property in $\left \{ w_k(y) > 0 \right \}$. 
Hence, using the Lipschitz continuity of $w$, we observe that the formula \eqref{1DStationary} is in fact valid in strips $\left \{ y_n^{k_1} \leq y_n \leq y_n^{k_2} \right \}$, the union of which contains
$B_{10r}(0)$. As a result, $w^{-1} \{0\} \cap B_{10r}(0)$ consists of pairwise disjoint hyperplanes, each separating exactly two subdomains $\left \{ w_k(y) > 0 \right \}$.
In particular, we deduce that $0 \in w^{-1} \{0\} \backslash \mathcal{S} \left ( w \right)$ and $I^w \left ( 0, 0^+ \right ) = 1$.
Therefore, for $\rho_0$ small enough $I^w \left ( 0, \rho \right ) < 1 + \delta_n/10$, whenever $\rho \in \left ( 0, \rho_0 \right )$.

On the other hand, using $\lambda_M^{(i)} \leq \lambda_M$, \eqref{FrequencyBoundSet}, Lemma \ref{CSComparisonLemma} and Corollary \ref{LinMonotonicity}, 
we note that there exists a $C = C \left (n, \lambda_M, \mathcal{I} \right )$ such that 
the functions $I^{(i)}(0,\rho) + C\rho^2$ are monotone nondecreasing in $\rho \in \left ( 0, \tilde{r}/2 \right )$.
Furthermore, by \eqref{JunctionPresence} and Corollary \ref{HomogeneousBlowupsGeneral}, for $\rho \in \left (0, \tilde{r}/2 \right )$,
\begin{equation*}
I^{(i)}(0,\rho) + C\rho^2 \geq \lim_{\rho \downarrow 0} \left ( I^{(i)}(0,\rho) + C \rho^2 \right ) \geq 1 + \delta_n.
\end{equation*}
By the strong and uniform convergence of $u^{(i)} \to w$ in $B_{8r}(0)$, for sufficiently small $\rho = C \left ( C, r, \delta_n \right ) \leq \rho_0$,
where uniformity in the index $i$ is guaranteed by the uniform bounds $\lambda_M^{(i)} \leq \lambda_M$ and \eqref{FrequencyBoundSet}, we have
\begin{equation*}
I^w (0,\rho) = \lim_{i \to \infty} I^{(i)}(0,\rho) \geq 1 + \delta_n / 5,
\end{equation*}
which is a contradiction.

{\emph{Step 3:}} Finally, we would like to estimate from above,
\begin{equation*}
\mathcal{R} = \int_{B_r(x)} \int_{B_{12r}(z) \backslash B_{4r}(z)} \left | (y-z) \cdot \nabla u(y) - \alpha u(y) \right |^2 \, \mathrm{d}y \, \mathrm{d}\mu(z),
\end{equation*}
with the choice,
\begin{equation*}
\alpha = \frac{1}{\mu \left ( B_r (x) \right )} \int_{B_r(x)} I_\phi(\zeta,8r) \, \mathrm{d}\mu(\zeta).
\end{equation*}
By the triangle inequality, $\mathcal{R} \leq (I) + (II)+ (III)$, where
\begin{equation*}
\begin{aligned}
& (I) = \int_{B_r(x)} \int_{B_{12r}(z) \backslash B_{4r}(z)} \left | I_\phi(z,8r) - I_\phi \left ( z, |z-y| \right )  \right |^2 |u(y)|^2 \, \mathrm{d}y \, \mathrm{d}\mu(z), \\
& (II) = \int_{B_r(x)} \int_{B_{12r}(z) \backslash B_{4r}(z)} \left | (y-z) \cdot \nabla u(y) - I_\phi \left ( z, |z-y| \right ) \right |^2 \, \mathrm{d}y \, \mathrm{d} \mu(z), \\
& (III) = \int_{B_r(x)} \int_{B_{12r}(z) \backslash B_{4r}(z)} \left ( I_\phi (z,8r) - \alpha \right )^2 \left | u(y) \right |^2 \, \mathrm{d}y \, \mathrm{d}\mu(z).
\end{aligned}
\end{equation*}

\noindent {\emph{Estimate on}} $(I)$: 
For $A$ as in Corollary \ref{LinMonotonicitySmoothedAdditive} and $c(A) = 2 + A + A^2$, expressing $I_\phi(z,8r) - I_\phi \left ( z, |z-y| \right )$ as
\begin{equation*}
\begin{aligned}
I_\phi(z,8r) - I_\phi \left ( z, |z-y| \right ) & = \left (  I_\phi(z,8r) + c(A) (8r)^2 \right ) - \\
                                                & \left ( I_\phi \left ( z, |z-y| \right ) + c(A) |z-y|^2 \right ) + c(A) \left [ |z-y|^2 - (8r)^2 \right ],
\end{aligned}
\end{equation*}
noting that $r \leq |y-z| \leq 32r$ holds, and recalling that $\mu$ is supported on $\mathcal{S}(u)$, by the monotonicity of $I_\phi(z,\rho) + c(A) \rho^2$ for $\rho \leq \tilde{r}$ and $z \in u^{-1} \{0\}$,
and the triangle inequality,
\begin{equation}
\left | I_\phi(z,8r) - I_\phi \left ( z, |z-y| \right ) \right | \leq W_{r,32r}^{c(A)}(z) + Cr^2,
\label{Step3_I_i}
\end{equation}
for some $C = C \left (n, \lambda_M, \mathcal{I}  \right )$.

We also note that since $32r \leq \tilde{r} \leq R \left ( n, N, \lambda_M \right )$, $\mu$ is supported on $\mathcal{S}(u)$, and $z \in B_r(x)$ in the integrand of $(I)$,
using \eqref{RenormalizedHeight}, \eqref{BoundaryPoincareComparison}, \eqref{LinMonotonicity}, \eqref{UsefulComparison2} and \eqref{ShiftedHeightComparison},
for $A_{4r,12r}(x) = B_{12r}(z) \backslash B_{4r}(z)$, we get
\begin{equation}
\begin{aligned}
\int_{A_{4r,12r}(x)} \left | u(y) \right |^2 \, \mathrm{d}y 
\leq 8r H(z, 12r) 
\leq C  r H(z,2r)
\leq C  r H_\phi (z,4r)
\leq C  r H_\phi(x,16r),
\end{aligned}
\label{Step3_I_ii}
\end{equation}
where the constants $C = C \left (n, \lambda_M, \mathcal{I} \right )$ have been updated for each inequality.

Hence, from \eqref{Step3_I_i} and \eqref{Step3_I_ii}, we obtain
\begin{equation*}
\left ( I \right ) \leq C r H_\phi(x,16r) \int_{B_r(x)} \left (  W_{r,32r}^{c(A)}(z) \right )^2 \, \mathrm{d}\mu(z) + C r H_\phi(x,16r) \mu \left ( B_r(x) \right ) r^4.
\end{equation*}
Finally, bounding $W_{r,32r}^{c(A)}(z)$ from above by using Lemma \ref{CSComparisonLemma}, 
the uniform frequency bound $I \left ( z, \tilde{r} \right ) \leq \mathcal{I}$, as $32r \leq \tilde{r} \leq R \left ( n, N, \lambda_M \right )$,
the almost monotonicity of $I(z,r)$, due to $z \in \mathrm{spt}(\mu) \subset \mathcal{S}(u)$,
applying the Cauchy-Schwarz inequality, and updating the constant $C_1 = C_1 \left (n, \lambda_M, \mathcal{I} \right )$ accordingly, we get
\begin{equation}
\left ( I \right ) \leq C_1 r H_\phi(x,16r) \int_{B_r(x)} W_{r,32r}^{c(A)}(z) \, \mathrm{d}\mu(z) + C_1 r H_\phi(x,16r) \mu \left ( B_r(x) \right ) r^4.
\label{Step3_IComplete}
\end{equation}

\noindent {\emph{Estimate on}} $(II)$: Noting that we can apply Proposition \ref{FrequencyOscillationEstimate}, as $\mathrm{spt} \left ( \mu \right ) \subset \mathcal{S}(u)$, 
using Corollary \ref{LinMonotonicitySmoothedAdditive} and arguing as in the proof of \eqref{Step3_I_ii},
we obtain 
\begin{equation*}
\begin{aligned}
\int_{A_{4r,12r}(z)} \left | (y-z) \cdot \nabla u(y) - I_\phi \left ( z, |z-y| \right ) \right |^2 \, \mathrm{d}y \, \mathrm{d} \mu(z)
& \leq C r H_\phi(z,12r) W_{2r,8r}^{c(A)}(z) \\
& \leq C r H_\phi(x, 16r) W_{r,32r}^{c(A)}(z),
\end{aligned}
\end{equation*}
where $C_2 = C_2 \left (n, \lambda_M, \mathcal{I} \right )$ has been updated for the second inequality. Integrating this inequality we obtain
\begin{equation}
(II) \leq C_2 r H_\phi(x, 16r) \int_{B_r(x)} W_{r,32r}^{c(A)}(z) \, \mathrm{d} \mu(z).
\label{Step3_IIComplete}
\end{equation}

\noindent {\emph{Estimate on}} $(III)$: By the Jensen inequality,
\begin{equation}
\left ( I_\phi (z,8r) - \alpha \right )^2 \leq \frac{1}{\mu \left ( B_r (x) \right )} \int_{B_r(x)} \left | I_\phi(z,8r) - I_\phi(\zeta,8r) \right |^2 \, \mathrm{d}\mu(\zeta).
\label{Step3_IIIi}
\end{equation}
By the definition of $\mathcal{R}$, $z,\zeta \in B_r(x)$, and therefore, $|z-y| < 2r$. By this observation and the fact that $\mu$ is supported on $\mathcal{S}(u)$,
we can apply Proposition \ref{FrequencyOscillationProposition} (with $x_1 = z$ and $x_2 = \zeta$) and Corollary \ref{LinMonotonicitySmoothedAdditive}, to obtain the pointwise estimate,
\begin{equation}
\left | I_\phi (z,8r) - I_\phi (\zeta,8r) \right |^2 \leq 4C^2 \left ( W_{r,32r}^{c(A)}(z) + W_{r,32r}^{c(A)}(\zeta) \right ) \frac{|z-\zeta|^2}{r^2} + 2C^2r^4 |z- \zeta|^2.
\label{Step3_IIIii}
\end{equation}
We note that since $32r \leq \tilde{r} \leq R \left (n, N, \lambda_M \right )$, $\mu$ is supported on $\mathcal{S}(u)$, and $z \in B_r(x)$ in the integrand of $(III)$, similarly to \eqref{Step3_I_ii}, we get
\begin{equation}
\int_{A_{4r,12r}(x)} \left | u(y) \right |^2 \, \mathrm{d}y 
\leq C  r H_\phi(x,16r).
\label{Step3_IIIiii}
\end{equation}
where $C = C \left (n, \lambda_M, \mathcal{I} \right )$ is updated at each inequality.
Hence, from \eqref{Step3_IIIi}, \eqref{Step3_IIIii}, \eqref{Step3_IIIiii}, the H\"older inequality and the fact that $\left | z - \zeta \right | < 2r$, as $z,\zeta \in B_r(x)$ in the integrand, we obtain
\begin{equation}
\begin{aligned}
(III) & \leq \frac{C  r H_\phi(x,16r)}{\mu \left ( B_r(x) \right )} \int_{B_r(x)} \int_{B_r(x)} \left [ W_{r,32r}^{c(A)}(z) + W_{r,32r}^{c(A)}(\zeta) + r^6 \right ] \mathrm{d}\mu(\zeta) \, \mathrm{d}\mu(z) \\
     & \leq C_3 r H_\phi(x,16r) \int_{B_r(x)} W_{r,32r}^{c(A)}(z) \, \mathrm{d}\mu(z) + C_3 r H_\phi(x,16r) \mu \left ( B_r(x) \right ) r^6, 
\end{aligned}
\label{Step3_IIIComplete}
\end{equation}
for an updated $C = C_3 \left (n, \lambda_M, \mathcal{I} \right )$.

Thus, from \eqref{Step3_IComplete}, \eqref{Step3_IIComplete}, \eqref{Step3_IIIComplete}, and $32r \leq \tilde{r} \leq R \left (n, N, \lambda_M \right )$, we obtain 
\begin{equation}
\mathcal{R} \leq C r H_\phi(x,16r) \int_{B_r(x)} W_{r,32r}^{c(A)}(z) \, \mathrm{d}\mu(z) + C r H_\phi(x,16r) \mu \left ( B_r(x) \right ) r^4
\label{Step3Completed}
\end{equation}
for $C_0 = C_0 \left (n, \lambda_M, \mathcal{I} \right )$ given by $\max \left \{ C_1, C_2, C_3 \tilde{r}^2 \right \}$.

{\emph{Step 4:}} From \eqref{JonesStep1c}, \eqref{JonesStep2Claim} and \eqref{Step3Completed}, we have
\begin{equation*}
\begin{aligned}
C_c^{-1} D_\phi(x,16r) r^{n} D_\mu^{n-2}(x,r) 
& \leq 2n \cdot C_0 r H_\phi(x,16r) \int_{B_r(x)} W_{r,32r}^{c(A)}(z) \, \mathrm{d}\mu(z) \\
& + C_0 r H_\phi(x,16r) \mu \left ( B_r(x) \right ) r^4.
\end{aligned}
\end{equation*}
Finally, we divide both sides by $C_c^{-1} D_\phi(x,16r) r^{n}$, and use the frequency lower bound \eqref{SFreqLowerBound}, as $x \in \mathcal{S}(u)$, 
as well as $32r \leq \tilde{r} \leq R \left (n, N, \lambda_M \right )$ and $c(A) = 2 + A + A^2$, to get
\begin{equation*}
D_\mu^{n-2} \left ( x , r \right ) \leq \frac{C}{r^{n-2}} \left [ \int_{B_r(x)} W_{r,32r}^{2+A+A^2} (z) \, \mathrm{d}\mu(z) + \mu \left ( B_r(x) \right ) r^2 \right ],
\end{equation*}
where $C$ depends on $n$, $\lambda_M$ and $\mathcal{I}$ only.
In order to obtain \eqref{JonesNumberEstimate}, we absorb into the integrand the second term inside the parantheses, consequently replace $W_{r,32r}^{2+A+A^2}$
with $W_{r,32r}^{3+A+A^2}$, and modify $C$ by multiplying by an absolute constant.
\end{proof}

\section{Approximate spines} \label{SpineSection}

In this section we focus on the singular sets of homogeneous and almost homogeneous $u \, : \, \Omega \to \Sigma_N$, 
Lipschitz continuous and stationary with respect to the variations \eqref{CombinedVariations} of functional \eqref{Functional}.

\subsection{Homogeneous maps of two variables} \label{HomSubsection}
For $u \, : \, \Omega \to \Sigma_N$ Lipschitz continuous and stationary with respect to the variations \eqref{CombinedVariations} of functional \eqref{Functional},
we recall that for every $k=1,2,...,N$, we have
\begin{equation}
- \Delta u_k = \lambda_k u_k \quad \mathrm{in} \quad \left \{ x \in \Omega \, : \, u_k(x) > 0 \right \}.
\label{ELrestated}
\end{equation}
Hence, if $u \, : \, \mathbf{R}^n \to \Sigma_n$ is a homogeneous map depending on the variables $x_{n-1}$, $x_n$ only, we can express it as
$u(x) = u \left ( x', r, \theta \right )$, where $x'= \left ( x_1, ..., x_{n-2} \right )$ and $(r,\theta)$ are the polar coordinates of the $x_{n-1}x_n$-plane. 
In particular, each component $u_k$ is a homogeneous function depending on the variables $x_{n-1}$, $x_n$ only, while simultaneously satisfying \eqref{ELrestated} in its support.
However, these two conditions imply that $\lambda_k = 0$ for $k=1,...,N$. In other words, $u$ is a stationary map with respect to the Dirichlet energy.
We observe that if we identify such a map $u(x)$ with the scalar function given by $\sum_{k=1}^N u_k(x)$, then there is an integer $m \geq 2$ such that, up to a rotation in the $x_{n-1}x_n$-plane,
\begin{equation}
u(x) = u \left ( x' , r, \theta \right ) = r^{m/2} \left | \cos \left ( \frac{m}{2} \theta \right ) \right |.
\label{CylindricalMaps} 
\end{equation}
When $m=2$, $u^{-1} \{0\}$ is the hypersurface $| \theta | = \pi/2$, separating the two nodal domains of $\cos(\theta)$, and $\mathcal{S}(u) = \emptyset$.
When $m \geq 3$, each nodal domain of $\cos \left ( \frac{m}{2} \theta \right )$ corresponds to the support of a component $u_k$,
and $\mathcal{S}(u)$ is given by the $(n-2)$-plane, $\left \{ x \in \mathbf{R}^n \, : \, x_{n-1} = x_n = 0 \right \}$, cf. \cite[Section 9]{CTV05}.
We note that $u$ is invariant with respect to $x' \in \mathbf{R}^{n-2}$. In this case, the linear subspace $\mathbf{R}^{n-2} \times \{0\}$ is often referred to as the {\emph{spine}} of the map $u$.

Finally, we remark that for any $u \, : \, \Omega \to \Sigma_N$ Lipschitz continuous and stationary with respect to the variations \eqref{CombinedVariations} of functional \eqref{Functional}, 
at $\mathcal{H}^{n-2}$-almost every $x \in \mathcal{S}(u)$, for every homogeneous blowup map at $x$  as in Proposition \ref{HomogeneousBlowupsGeneral}, there is a coordinate system of $\mathbf{R}^n$ in which
the blowup map is given \eqref{CylindricalMaps} with an integer $m \geq 3$. In fact, $m = I \left ( x, 0^+ \right )$.
This follows from Proposition \ref{HomogeneousBlowupsGeneral} and a classical Almgren-type stratification, cf. \cite[Section 3.4]{Simon2}.

\subsection{Almost homogeneous maps} \label{AlmostHomSubsection}
By the final remark in the preceding section, it is a natural goal to identify an {\emph{approximate spine}} for maps that are almost homogeneous, with the purpose of studying $\mathcal{S}(u)$.
In particular, it is reasonable to expect that approximate spines of $u$ would serve as good approximations of $\mathcal{S}(u)$ in sufficiently small scales.
Towards this goal we review the definitions of quantitative linear independence and spanning below. 

\begin{definition} \label{QuantLinearAlgebra}
A set of points $x_0, x_1,..., x_k \subset B_r(x)$ is called $\rho r$-linearly independent, if for all $i=1,..,k$,
\begin{equation}
\mathrm{dist} \left ( x_i, x_0 + \mathrm{span} \left \{ x_{i-1} - x_0, ..., x_1 - x_0 \right \} \right ) \geq \rho r,
\label{QuantLinearIndependence}
\end{equation}
where $\mathrm{span}A$ is defined as the linear subspace generated by the elements of $A \subset \mathbf{R}^n$ with the convention $\mathrm{span} \emptyset = \{0\}$.

A set $F \subset B_r(x)$ is said to $\rho r$-span a $k$-dimensional affine subspace $V$, if there are $r\rho$-linearly independent $x_0, x_1,..., x_k \subset F$ such that 
$$V = x_0 + \mathrm{span} \left \{x_1 - x_0,..., x_k - x_0 \right \}.$$
\end{definition}

\begin{remark} \label{UsefulGeometricRemark}
We observe that if a set $F \cap B_r(x)$ fails to $\rho r$-span any $k$-dimensional affine subspace, then there exists a $(k-1)$-dimensional affine subspace $L$ such that
$F \cap B_r(x)$ is contained in $B_{\rho r}(L) = \left \{ x \in \mathbf{R}^n \, : \, \mathrm{dist}(x,L) < \rho r \right \}$.
\end{remark}

Next we prove a technical $\epsilon$-clearing lemma for the singular set $\mathcal{S}(u)$ based on a compactness argument, which will be useful in pinning down $\mathcal{S}(u)$ near suitably chosen affine $(n-2)$-planes.

\begin{lemma} \label{TubularNeighborhoodLemma}
Suppose $u \, : \, \Omega \to \Sigma_N$ is a Lipschitz continuous map that is stationary with respect to the variations \eqref{CombinedVariations} of functional \eqref{Functional}, 
$\tilde{r} \in \left ( 0, R \left ( n, N, \lambda_M \right ) \right )$, $c(A) \geq A$ for $A$ as in Corollary \ref{LinMonotonicitySmoothedAdditive},
$K \subset \Omega$ is compact with $\mathrm{dist}_{\mathcal{H}} \left ( K, \Omega \right ) > \tilde{r}$, $r \in \left ( 0, \tilde{r} \right )$.
Given $\rho, \overline{\rho}, \tilde{\rho} \ll r$, there exists an $\epsilon > 0$, depending on $n$, $\lambda_M$, $\mathcal{I}$, $\rho$, $\overline{\rho}$, $\tilde{\rho}$ only, such that the following holds.
For $x_0$, $x_1$, ..., $x_{n-2}$ in $B_r(x) \cap \mathcal{S}(u) \cap K$, $\rho r$-linearly independent points such that for every $i=0,1,...,n-2$,
\begin{equation}
W_{\tilde{\rho},2r}^{c(A)} \left ( x_i \right ) = \left [ I_\phi \left ( x_i, 2r \right ) + c(A)(2r)^2 \right ] - \left [ I_\phi \left ( x_i, \tilde{\rho} \right ) + c(A) \tilde{\rho}^2 \right ] < \epsilon,
\label{UniformlySmallDrop}
\end{equation}
and $V = x_0 + \mathrm{span} \left \{ x_1 - x_0, ..., x_{n-2} - x_0 \right \}$, we have
\begin{equation}
\mathcal{S}(u) \cap K \cap \left ( B_r(x) \backslash B_{\overline{\rho}} (V) \right ) = \emptyset.  
\label{TubularContainment}
\end{equation}
\end{lemma}

\begin{proof}
Suppose the claim is false.
Then there exist some fixed $\rho$, $\tilde{\rho}$, $\overline{\rho} > 0$, a sequence of maps $u^{(j)}$ satisfying the hypothesis (with corresponding $\lambda_k^{(j)} \leq \lambda_M$),
and a sequence of collections of $\rho r$-linearly independent points $x_0^{(j)}$, $x_1^{(j)}$, ..., $x_{n-2}^{(j)}$ in $B_r(0)$,
where we can assume $x = 0$ is fixed by precomposing each $u^{(j)}$ with a translation by $x_j \in K$, such that for each $i = 0,1,...n-2$,
\begin{equation}
\lim_{j \to \infty} W_{\tilde{\rho},2r, u^{(j)}}^{c(A)} \left ( x_i^{(j)} \right ) = 0,
\label{ContraDrop}
\end{equation}
while for $V^{(j)} = x_0^{(j)} + \mathrm{span} \left \{ x_1^{(j)} - x_0^{(j)}, ..., x_{n-2}^{(j)} - x_0^{(j)} \right \}$, we have
\begin{equation}
z^{(j)} \in \mathcal{S} \left ( u^{(j)} \right ) \cap K \cap \left ( \overline{B_r(0)} \backslash B_{\overline{\rho}} \left ( V^{(j)} \right ) \right ).
\label{ContraPoint}
\end{equation}
By rescaling each map, we can assume that $H_\phi^{(j)} \left ( 0, \tilde{r} \right ) = 1$, 
which together with \eqref{SLogDerivHeight}, the uniform bound $I_{\phi}^{(j)} \left (0, \tilde{r} \right ) \leq \mathcal{I}$, 
and the almost monotonicity of frequency enable us to check the hypothesis of Lemma \ref{ClassicalCompactnessLemma} in $B_{\tilde{r}/2}(0)$. 

Consequently, we have a subsequence $j'$ relabeled as $j$ such that $u^{(j)}$ converge to 
a nontrivial map $w \, : \, B_{\tilde{r}/2}(0) \to \Sigma_N$ in $L^2$, $H^1_{loc}$, $C^{0,\alpha}_{loc}$ for every $\alpha \in (0,1)$,
and $w$ is stationary with respect to the variations \eqref{CombinedVariations} of functional \eqref{Functional} with some limiting $\tilde{\lambda}_k \leq \lambda_M$,
while $x_i^{(j)}$, $i=0,..,n-2$, converge to $\rho r$-linearly independent $x_0$, $x_1$, ..., $x_{n-2}$ in $\overline{B_r(0)}$,
and by \eqref{ContraPoint},
$$ z^{(j)} \to z \in \mathcal{S}(w) \cap K \cap \overline{\left ( B_r(0) \backslash B_{\overline{\rho}} \left ( V \right ) \right )},$$
where $V = x_0 + \mathrm{span} \left \{ x_1 - x_0, ..., x_{n-2}- x_0\right \}$ is an affine $(n-2)$-plane by the stability of $\rho r$-linear independence under convergence.
We remark that $x_0, ..., x_{n-2}, z \in \mathcal{S}(w)$ follows from the argument in {\emph{Step 2}} of the proof of Proposition \ref{JonesNumberProposition}, 
due to the gap condition for frequency at points in the singular set.
Lastly, by \eqref{ContraDrop} and locally strong and uniform convergence, we have
$$W_{\tilde{\rho},2r, w}^{c(A)} \left ( x_i\right ) = 0, \quad i=0,...,n-2.$$
Hence, \eqref{SLinMonotonicityFormulaAdditive} implies that 
$w$ is homogeneous of degree $\alpha_i$ in $\left | y - x_i \right |$ in the annulus $B_{2r} \left ( x_i \right ) \backslash B_{\tilde{\rho}} \left ( x_i \right )$.
In particular, as in Section \ref{HomSubsection}, we observe that $\tilde{\lambda}_k = 0$ for $k=1,...,N$. 
In other words, $w$ is a stationary map with respect to the variations \eqref{CombinedVariations} of the Dirichlet energy.

We note that each component $w_k$ has the unique continuation property in $P_k = \left \{ x \in \mathbf{R}^n \, : \, w_k(x) > 0 \right \}$, as it satisfies $- \Delta w_k = 0$ in $P_k$, $k=1,...,N$.
Since $\alpha_i$ homogeneity of $w$ implies that $P_k$ is a cone in each annulus $B_{2r} \left ( x_i \right ) \backslash B_{\tilde{\rho}} \left ( x_i \right )$, we can extend each $w_k$,
and therefore $w$, from $B_{2r} \left ( x_i \right ) \backslash B_{\tilde{\rho}} \left ( x_i \right )$ to $\mathbf{R}^n$. We denote each extension as $w^{i}$ for $i=0,...,n-2$. 
Since each $x_i \in B_r(0)$, any two annuli $B_{2r} \left ( x_i \right ) \backslash B_{\tilde{\rho}} \left ( x_i \right )$ clearly overlap. 
Hence, $w^0 = w^1 = ... = w^{n-2} = w$ in $B_{\tilde{r}/2}(0)$ by the unique continuation property in each $P_k$, $k=1,...,N$.
Thus, it is easy to check via a blow-down argument (cf. \cite[Lemma 6.8]{DLMSV}) that $\alpha_i = \alpha$ for each $i=0,...,n-2$, 
and in $B_{\tilde{r}/2}(0)$, $w$ is homogeneous of degree $\alpha$ in the variable $y' \in \left ( V - x_0 \right )^{\perp}$. 
Moreover, $w \left ( x_0 + v \right ) = w \left ( x_0 \right )$ for every $v \in \left ( V - x_0 \right )$, 
and therefore, $V \cap B_{\tilde{r}/2}(0) \subset \mathcal{S}(w)$. 
However, $z \in \mathcal{S}(w) \backslash V$ and the homogeneity of $w$ in the variable $y' \in \left ( V - x_0 \right )^{\perp}$ in $B_{\tilde{r}/2}(0)$ 
together imply that $\mathcal{S}(w)$ has Hausdorff dimension $(n-1)$,
contradicting Corollary \ref{LimitNodalDimension}.
\end{proof}

\subsection{Frequency oscillations on approximate spines} \label{FOASsubsection}
We end this section by showing that for $u$, $x_0$, $x_1$, ..., $x_{n-2} \in \mathcal{S}(u)$ satisfying the hypothesis of Lemma \ref{TubularNeighborhoodLemma}, with a sufficiently small $\epsilon > 0$,
on the approximate spine $V$ spanned by these points, we can locally control the oscillations of frequency $I_\phi(y,\rho)$ over space and comparable scales.

\begin{lemma} \label{FOSSLemma}
Suppose $u \, : \, \Omega \to \Sigma_N$ is a Lipschitz continuous map that is stationary with respect to the variations \eqref{CombinedVariations} of functional \eqref{Functional}, 
$\tilde{r} \in \left ( 0, R \left ( n, N, \lambda_M \right ) \right )$, $c(A) \geq A$ for $A$ as in Corollary \ref{LinMonotonicitySmoothedAdditive},
$K \subset \Omega$ is compact with $\mathrm{dist}_{\mathcal{H}} \left ( K, \Omega \right ) > \tilde{r}$, $r \in \left ( 0, \tilde{r} \right )$,
$\rho, \overline{\rho}, \tilde{\rho} \ll r$.
For every $\delta > 0$, there exists an $\epsilon > 0$, depending on $\delta$, $n$, $\lambda_M$, $\mathcal{I}$, $\rho$, $\overline{\rho}$, $\tilde{\rho}$ only, such that the following holds.
For $x_0$, $x_1$, ..., $x_{n-2}$ in $B_r(x) \cap \mathcal{S}(u) \cap K$, $\rho r$-linearly independent points such that for every $i=0,1,...,n-2$,
\begin{equation}
W_{\tilde{\rho},2r}^{c(A)} \left ( x_i \right ) = \left [ I_\phi \left ( x_i, 2r \right ) + c(A)(2r)^2 \right ] - \left [ I_\phi \left ( x_i, \tilde{\rho} \right ) + c(A) \tilde{\rho}^2 \right ] < \epsilon,
\label{ReUniformlySmallDrop}
\end{equation}
and $V = x_0 + \mathrm{span} \left \{ x_1 - x_0, ..., x_{n-2} - x_0 \right \}$, for every $y, \tilde{y} \in B_r(x) \cap V$ and every $s$, $\tilde{s} \in [\rho,r]$, we have
\begin{equation}
\left | I_\phi(y,s) - I_\phi \left ( \tilde{y}, \tilde{s} \right ) \right | < \delta.
\label{FOSSControl}
\end{equation}
\end{lemma}

\begin{proof}
The proof is analogous to the proof of Lemma \ref{TubularNeighborhoodLemma}. Suppose the claim is false for some fixed $\delta$, $\rho$, $\overline{\rho}$, $\tilde{\rho} > 0$.
Then there exist a sequence of maps $u^{(j)}$ satisfying the hypothesis (with corresponding $\lambda_k^{(j)} \leq \lambda_M$),
and a sequence of collections of $\rho r$-linearly independent points $x_0^{(j)}$, $x_1^{(j)}$, ..., $x_{n-2}^{(j)}$ in $B_r(0)$,
where we can assume $x = 0$ is fixed by precomposing each $u^{(j)}$ with a translation by $x_j \in K$, such that for each $i = 0,1,...n-2$,
\begin{equation}
\lim_{j \to \infty} W_{\tilde{\rho},2r, u^{(j)}}^{c(A)} \left ( x_i^{(j)} \right ) = 0,
\label{ContraDropRe}
\end{equation}
while for $V^{(j)} = x_0^{(j)} + \mathrm{span} \left \{ x_1^{(j)} - x_0^{(j)}, ..., x_{n-2}^{(j)} - x_0^{(j)} \right \}$, 
there are points $y^{(j)}$, $\tilde{y}^{(j)} \in V^{(j)} \cap B_r(0)$ and scales $s_j$, $\tilde{s}_j \in [\rho,r]$ such that
\begin{equation}
\left | I_\phi \left ( y^{(j)} ,s_j \right ) - I_\phi \left ( \tilde{y}^{(j)} , \tilde{s}_j \right ) \right | \geq \delta.
\label{ContraFOSS}
\end{equation}
By rescaling each map, we can assume that $H_\phi^{(j)} \left ( 0, \tilde{r} \right ) = 1$, 
which together with \eqref{SLogDerivHeight}, the uniform bound $I_{\phi}^{(j)} \left (0, \tilde{r} \right ) \leq \mathcal{I}$, 
and the almost monotonicity of frequency enable us to check the hypothesis of Lemma \ref{ClassicalCompactnessLemma} in $B_{\tilde{r}/2}(0)$.

As a result, we have a subsequence $j'$ relabeled as $j$ such that $u^{(j)}$ converge to 
a nontrivial map $w \, : \, B_{\tilde{r}/2}(0) \to \Sigma_N$ in $L^2$, $H^1_{loc}$, $C^{0,\alpha}_{loc}$ for every $\alpha \in (0,1)$,
and $w$ is stationary with respect to the variations \eqref{CombinedVariations} of functional \eqref{Functional} with some limiting $\tilde{\lambda}_k \leq \lambda_M$.
Likewise, $x_i^{(j)}$, $i=0,..,n-2$, converge to $\rho r$-linearly independent $x_0$, $x_1$, ..., $x_{n-2}$ in $\overline{B_r(0)}$, 
and by the locally strong and uniform convergence of $u^{(j)}$ in $B_{\tilde{r}/2}(0)$, $I_{\phi}^{(j)} \left ( y^{(j)},s_j \right ) \to I_{\phi,w} (y,s)$ and
$I_{\phi}^{(j)} \left ( \tilde{y}^{(j)},\tilde{s}_j \right ) \to I_{\phi,w} \left ( \tilde{y}, \tilde{s} \right )$, and consequently,
\begin{equation}
\left | I_\phi \left ( y ,s \right ) - I_\phi \left ( \tilde{y} , \tilde{s} \right ) \right | \geq \delta,
\label{ContraFOSSLimit}
\end{equation}
for some $y, \tilde{y} \in V$ and $s, \tilde{s} \in [\rho,r]$,
where $V = x_0 + \mathrm{span} \left \{ x_1 - x_0, ..., x_{n-2}- x_0\right \}$ is an affine $(n-2)$-plane by the stability of $\rho r$-linear independence under convergence.
We remark that $x_0, ..., x_{n-2}, z \in \mathcal{S}(w)$ follows from the argument in {\emph{Step 2}} of the proof of Proposition \ref{JonesNumberProposition}, 
due to the gap condition for frequency at points in the singular set.

As in Section \ref{HomSubsection} and Lemma \ref{TubularNeighborhoodLemma}, using \eqref{ContraDropRe} we observe that $\tilde{\lambda}_k = 0$ for $k=1,...,N$. 
In other words, $w$ is a stationary map with respect to the variations \eqref{CombinedVariations} of the Dirichlet energy.
And repeating the homogeneous extension and unique continuation argument in the proof of Lemma \ref{TubularNeighborhoodLemma}, we conclude that there exists an $\alpha > 0$ such that 
$$
I_{\phi, w} (z,s) = \alpha, \quad \forall z \in V, \forall s > 0,  
$$
which contradicts \eqref{ContraFOSSLimit}.
\end{proof}

\section{Size estimate} \label{MinkowskiSection}

In this section we prove a local Minkowski-type estimate, which implies Theorem \ref{MainTheorem1} immediately through an obvious covering argument.
The proof of this estimate is based on the covering arguments in \cite{NV2, DLMSV} and relies on the discrete Reifenberg theorem of Naber and Valtorta from \cite{NV}.
Note that we prove this result for the singular set $\mathcal{S}(u)$ of arbitrary Lipschitz continuous maps $u \, : \, \Omega \to \Sigma_N$, 
which satisfy stationarity under the variations \eqref{CombinedVariations} of functional \eqref{Functional}.
In particular, our conclusion applies to the singular set $\mathcal{S}(u)$ of optimal partitions in the sense of minimization problem \eqref{EigenvalueMinimization}.

\subsection{Minkowski-type estimate} \label{MTESubsection}
Below we state the main estimate of this section and give a very rough summary of how it will be proved.

\begin{theorem} \label{MTETheorem}
Suppose $u \, : \, \Omega \to \Sigma_N$ is a Lipschitz continuous map that is stationary with respect to the variations \eqref{CombinedVariations} of functional \eqref{Functional}, 
and $K \subset \Omega$ is compact.  
Then there exist an $\tilde{r} = \tilde{r} \left (n, N, \lambda_M, \mathcal{I} \right ) \in \left (0, \min \left \{ R \left (n, N, \lambda_M \right ),  \mathrm{dist} \left ( K , \partial \Omega \right ) \right \}  \right )$  
and a positive constant $C = C \left ( n, N, \lambda_M, \mathcal{I} \right )$ such that for every $x_0 \in K$, $0 < \rho \leq r \leq \tilde{r}$,
\begin{equation}
\mathcal{L}^n \left ( B_\rho \left ( \mathcal{S}(u) \cap K \right ) \cap B_r \left (x_0 \right ) \right ) \leq C r^{n-2} \rho^2. \label{MTE}
\end{equation}
\end{theorem}

Our strategy is to to find for any $\rho \in \left ( 0, \tilde{r} \right )$ a collection of balls $B_r \left ( x_i \right )$ covering $\mathcal{S}(u) \cap K \cap B_{\tilde{r}} \left (x_0 \right )$ such that 
$\mathcal{N}(\rho)$, the number of balls in this collection, is less than or equal to $C r^{n-2} \rho^{2-n}$ for $C$ as in the statement of theorem. 
Once we derive such a bound, observing that
$B_\rho \left ( \mathcal{S}(u) \cap K \right ) \cap B_{r} \left (x_0 \right ) \subset \cup_i B_{2\rho} \left ( x_i \right )$, we arrive at the conclusion,
\begin{equation}
\mathcal{L}^n \left ( B_\rho \left ( \mathcal{S}(u) \cap K \right ) \cap B_{r} \left (x_0 \right ) \right ) \leq 2^n \mathcal{N}(\rho) \rho^n \leq C r^{n-2} \rho^2.
\label{MTEFinalStep}
\end{equation}
The construction of such a collection of balls, and most importantly, the upper bound on its cardinality, will be given in Sections \ref{ECSubsection}, \ref{ICSubsection} and \ref{DRSubsection}.

\begin{proof}[Proof of Theorem \ref{MainTheorem1}]
Combining the estimate \eqref{MTE} with a standard covering argument immediately gives the Minkowski-type estimate \eqref{StrongestSizeEstimate} on $\mathcal{S}(u) \cap K$. 
\end{proof}

\subsection{Efficient covering} \label{ECSubsection}
The collection of balls that yields \eqref{MTEFinalStep} will be a product of the following covering lemma.
While we follow \cite[Section 7.1]{DLMSV} in this section, in the general case $\lambda_M \neq 0$,
we have the additional task of finding a sufficiently small scale at which analogous covering lemmas hold.

\begin{lemma} \label{ECLemma}
Let $u$, $K \subset \Omega$ be as in Theorem \ref{MTETheorem}, the constant $A$ as in Corollary \ref{LinMonotonicitySmoothedAdditive}. Set the constant $c(A) = 3 + A + A^2$.
There exists an $\tilde{r} = \tilde{r} \left (n, N, \lambda_M, \mathcal{I} \right ) \in \left (0, \min \left \{ R \left (n, N, \lambda_M \right ),  \mathrm{dist} \left ( K , \partial \Omega \right ) \right \}  \right )$
such that the following holds.
Given any $x \in K$, $0 < s < r \leq \tilde{r}$, and $D \subset \mathcal{S}(u) \cap K \cap B_r(x)$ with $U = \sup \left \{ I_\phi^{c(A)} (y,r) \, : \, y \in D \right \}$, 
there exist a $\delta = \delta \left ( n, N, \lambda_M, \mathcal{I}, \tilde{r} \right ) > 0$, a constant $C = C(n) \geq 1$,
and a finite collection of balls $\left \{ B_{s_i} \left ( x_i \right ) \right \}$ covering $D$ and a corresponding decomposition of $D$ in sets $A_i \subset D$ with the following properties:
\begin{enumerate}[(a)]
 \item $A_i \subset B_{{s_i}} \left ( x_i \right )$ and $r \geq s_i \geq s$. \label{RadiusLB}
 \item For each $i$, either $s_i = s$, or 
\begin{equation}
 \sup \left \{ I_\phi^{c(A)} \left ( y, s_i \right ) \, : \, y \in D \cap B_{s_i} \left ( x_i \right ) \right \} \leq U - \delta. \label{FreqDropCondition}
\end{equation} \label{RadiusDichotomy}
 \item $\sum_i s_i^{n-2} \leq C r^{n-2}$. \label{DiscreteSum}
\end{enumerate}
\end{lemma}

We postpone the proof of Lemma \ref{ECLemma} to the following sections, and use it to prove Theorem \ref{MTETheorem} instead.

\begin{proof}[Proof of Theorem \ref{MTETheorem}]
We denote $D_0 = \mathcal{S}(u) \cap K \cap B_{\tilde{r}} \left (x_0 \right )$.
From Remark \ref{UsefulFreqMax}, Corollary \ref{LinMonotonicity}, Lemmas \ref{CSComparisonLemma} and \ref{ElementaryUpperBounds}, we deduce the bound
$$
U_0 = \sup \left \{ I_\phi^{c(A)} \left ( y, \tilde{r} \right ) \, : \, y \in D_0 \right \} \leq C_0 \left ( n, N, \lambda_M, \mathcal{I} \right ).
$$
We apply Lemma \ref{ECLemma} to $D = D_0$ with $r \leq \tilde{r}$, $s = \rho$. 
We denote the resulting cover as $\left \{ B_{s_i} \left ( x_i \right ) \right \}_{i \in I_1}$, with the corresponding decomposition $A_i \subset D \cap B_{s_i} \left ( x_i \right )$ for each $i \in I_1$. 
Hence,
$$
\sum_{i \in I_1} s_i^{n-2} \leq C(n) r^{n-2}.
$$
We decompose $I_1 = I_1^g \cup I_1^b$, where $I_1^g = \left \{ i \, : \, s_i = \rho \right \}$.
Hence, for each $i \in I_1^b$, 
$$
\sup \left \{ I_\phi^{c(A)} \left ( y, s_i \right ) \, : \, y \in A_i \right \} \leq U_0 - \delta.
$$
For each $i \in I_1^b$, applying Lemma \ref{ECLemma} with $D = A_i$, $r = s_i$, $s = \rho$,
we obtain collections of balls $ \left \{ B_{s_{i,j}} \left ( x_{i,j} \right ) \right \}$, $j \in I_1^{b,i}$, subsets $A_{i,j} \subset A_i \cap B_{s_{i,j}} \left ( x_{i,j} \right )$, and estimates,
$$
\sum_{j \in I_1^{b,i}} s_{i,j}^{n-2} \leq C(n) \tilde{r}^{n-2} s_i^{n-2}.
$$
Hence, letting $I_2 = I_1^g \cup \left ( \cup_i I_1^{b,i} \right )$, we get
$$
\sum_{i \in I_2} s_i^{n-2} \leq C(n) \tilde{r}^{n-2} \sum_{i \in I_1} s_i^2 \leq C(n)^2 \tilde{r}^{2(n-2)},
$$
after reindexing sets $\left \{ B_{s_i} \left ( x_i \right ) \right \}$, radii $s_i$, and subsets $A_i \subset D \cap B_{s_i} \left ( x_i \right )$ for $i \in I_2$.
But for $s_i > \rho$, $i \in I_2$, we have the improved frequency drop,
$$
\sup \left \{ I_\phi^{c(A)} \left ( y, s_i \right ) \, : \, y \in A_i \right \} \leq U_0 - 2\delta.
$$
By induction, for every $k$, we get $\left \{ A_i \right \}_{i \in I_k}$ corresponding to a collection of balls $\left \{ B_{s_i} \left ( x_i \right ) \right \}_{i \in I_k}$ such that
either $s_i = \rho$, or 
$$ \sup \left \{ I_\phi^{c(A)} \left ( y, s_i \right ) \, : \, y \in A_i \right \} \leq U - k \delta,$$
and
\begin{equation}
\sum_{i \in I_k} s_i^{n-2} \leq C(n)^{k} \tilde{r}^{k(n-2)}. 
\label{UltimateSum} 
\end{equation}
By the positivity of frequency function, this process terminates after $\kappa$ steps, for a positive integer $\kappa \leq 1 + U_0 / \delta$, and yields a collection of balls $B_\rho \left ( x_i \right )$,
which cover $D_0 = \mathcal{S}(u) \cap K \cap B_{\tilde{r}} \left (x_0 \right )$, and by \eqref{UltimateSum},
\begin{equation}
\mathcal{N}(\rho) \rho^{n-2} = \sum_{i \in I_\kappa} s_i^{n-2} \leq C(n)^\kappa \tilde{r}^{\kappa(n-2)}.
\label{CardinalityBound}  
\end{equation}
Observing that $B_\rho \left ( \mathcal{S}(u) \cap K \right ) \cap B_{\tilde{r}} \left (x_0 \right ) \subset \cup_i B_{2\rho} \left ( x_i \right )$, by \eqref{CardinalityBound}, we obtain the desired estimate,
$$
\mathcal{L}^n \left ( B_\rho \left ( \mathcal{S}(u) \cap K \right ) \cap B_{\tilde{r}} \left (x_0 \right ) \right ) \leq 2^n \mathcal{N}(\rho) \rho^n \leq C \tilde{r}^{n-2}  \rho^2.
$$
for $C = 2^n C(n)^\kappa R \left (n, N, \lambda_M \right )^{(\kappa - 1) (n-2)}$. 
Finally, since $\kappa \leq 1 + U_0 / \delta$, where $U_0$ and $\delta$ depend on $n$, $N$, $\lambda_M$ and $\mathcal{I}$ only, $C$ depends exclusively on these parameters as well. 
\end{proof}

\subsection{Intermediate covering} \label{ICSubsection}
The efficient covering in Lemma \ref{ECLemma} will be obtained with the help of an intermediate covering lemma.

\begin{lemma} \label{ICLemma}
Let $u$, $K \subset \Omega$ be as in Theorem \ref{MTETheorem}, the constant $A$ as in Corollary \ref{LinMonotonicitySmoothedAdditive}. Set the constant $c(A) = 3 + A + A^2$.
There exists an $\tilde{r} = \tilde{r} \left (n, N, \lambda_M, \mathcal{I} \right ) \in \left (0, \min \left \{ R \left (n, N, \lambda_M \right ),  \mathrm{dist} \left ( K , \partial \Omega \right ) \right \}  \right )$
such that the following holds.
Given any $x \in K$, $0 < \rho \leq 1/100$, $0 < \sigma < \tau \leq \tilde{r}$ and $D \subset \mathcal{S}(u) \cap K \cap B_{\tilde{r}}(x)$ with 
$U = \sup \left \{ I_\phi^{c(A)} \left (y, \tilde{r} \right ) \, : \, y \in D \right \}$, 
there exist a $\delta = \delta \left ( n, N, \lambda_M, \mathcal{I}, \tilde{r}, \rho \right ) > 0$, a constant $C_R(n) > 0$,
and a finite collection of balls $\left \{ B_{r_i} \left ( x_i \right ) \right \}$ covering $D$ with the following properties:
\begin{enumerate}[(a)]
 \item $10 \rho \sigma \leq r_i \leq \tilde{r}$. \label{RadiusLBR}
 \item For each $i$, either $r_i \leq \sigma$, or there is an $(n-3)$-dimensional affine subspace $L_i \subset \mathbf{R}^n$ such that the set of points, 
 $$F_i = D \cap B_{r_i} \left ( x_i \right ) \cap \left \{ y \, : \, I_\phi^{c(A)} \left (y, \rho r_i \right ) > U - \delta \right \},$$
 is contained in $B_{\rho r_i} \left ( L_i \right ) \cap B_{r_i} \left ( x_i \right )$. \label{RadiusDichotomyR}
 \item $\sum_i r_i^{n-2} \leq C_R(n) \tau^{n-2}$. \label{DiscreteSumR}
\end{enumerate}
\end{lemma}

Firstly, we show that for $\rho > 0$ sufficiently small, Lemma \ref{ICSubsection} implies Lemma \ref{ECLemma}. 
While the proof is identical to the proof of \cite[Lemma 7.3]{DLMSV}, we include it for the reader's convenience. 

\begin{proof}[Proof of Lemma \ref{ECLemma}]
The key issue is to find a sufficiently small $\rho_0 = \rho_(0)(n) > 0$ such that 
choosing a $\rho \in \left ( 0, \rho_0 \right )$, the conclusions of Lemma \ref{ECLemma} can be deduced from Lemma \ref{ICSubsection}.
We treat $\rho \in (0, 1/100)$ as a fixed number for the moment.

Firstly, we apply Lemma \ref{ICLemma} to $B_r(x)$ for $\tau = r \leq \tilde{r}$ and $\sigma = s$, and obtain a first covering $\mathcal{C}(0) = \left \{ B_{r_i} \left ( x_i \right ) \right \}$.
We group the balls in this covering as $\mathcal{C}(0) = \mathcal{G}(0) \cup \mathcal{B}(0)$, where 
$\mathcal{G}(0) = \left \{ B_{r_i} \left ( x_i \right ) \, : \, r_i \leq s \right \}$ and $\mathcal{B}(0) = \left \{ B_{r_i} \left ( x_i \right ) \, : \, r_i > s \right \}$.

For each $B_{r_i} \left ( x_i \right ) \in \mathcal{B}(0)$, we consider $F_i = F \left ( B_{r_i} \left ( x_i \right ) \right )$ and $L_i$ as in Lemma \ref{ICLemma}, \eqref{RadiusDichotomyR}.
Then each $B_{2\rho r_i} \left ( L_i \right ) \cap B_{r_i} \left ( x_i \right )$ can be covered by $\mathcal{N}$ balls of radius $4 \rho r_i$, where $\mathcal{N} \leq C(n) \rho^{3-n}$,
since $L_i$ is an $(n-3)$-dimensional affine subspace.
If $4 \rho r_i < s$, we add these balls in a new collection $\mathcal{C}(1)$. 
We apply Lemma \ref{ICLemma} to those balls with $4 \rho r_i \geq s$ and include the resulting balls in $\mathcal{C}(1)$ as well.
We note that the collection $\mathcal{C}(1)$ has the property that
\begin{equation*}
\sum_{B_{r_i} \left ( x_i \right ) \in \mathcal{C}(1)} r_i^{n-2} 
\leq C(n) \rho^{3-n} \sum_{B_{r_j} \left ( x_j \right ) \in \mathcal{C}(0)} \left ( 4 \rho r_j \right )^{n-2}
= C(n) 4^{n-2} \rho \sum_{B_{r_j} \left ( x_j \right ) \in \mathcal{C}(0)} r_j^{n-2},
\end{equation*}
which, after letting $\rho_0(n) = \left (2 C(n) \right )^{-1} 4^{2-n}$ and fixing a $\rho \in \left ( 0, \rho_0 \right )$, yields
\begin{equation}
\sum_{B_{r_i} \left ( x_i \right ) \in \mathcal{C}(1)} r_i^{n-2} \leq \frac{1}{2} \sum_{B_{r_j} \left ( x_j \right ) \in \mathcal{C}(0)} r_j^{n-2}. \label{SummableBound} 
\end{equation}
Note that fixing $\rho$ also determines $\delta$ in Lemma \ref{ICLemma} fully, and consequently, induces a smallness requirement on $\delta$ in this context as well.

We repeat this procedure finitely many times, until we obtain a collection $\mathcal{C}(k)$ containing balls of radius less than or equal to $s$. 
Setting $\mathcal{C} = \cup_{j \leq k} \mathcal{C}(j)$ and using \eqref{SummableBound} at each step, we get
\begin{equation}
\sum_{ B_{r_i} \left ( x_i \right ) \in \mathcal{C}} r_i^{n-2} \leq \sum_{\ell = 0}^k 2^{- \ell} \sum_{B_{r_j} \left ( x_j \right ) \in \mathcal{C}(0)} r_j^{n-2} \leq 2 C_R(n) r^{n-2}, \label{TowardsDiscreteSum}
\end{equation}
where $C_R(n)$ is as in Lemma \ref{ICLemma}, \eqref{DiscreteSumR}.

Secondly, we set $A_i' = D \cap B_{r_i} \left ( x_i \right )$ for every $B_{r_i} \left ( x_i \right ) \in \mathcal{G}(0)$, that is $r_i \leq s$. 
Otherwise, we set  $A_i' = \left ( D \cap B_{r_i} \left ( x_i \right ) \right ) \backslash F_i$ for $F_i = F_i \left ( B_{r_i} \left ( x_i \right ) \right )$ as in Lemma \ref{ICLemma}, \eqref{RadiusDichotomyR}.
By its construction, $\mathcal{C}(1)$ covers $F_i$, and hence,
$$
D \subset \left ( \bigcup_{B_{r_i} \left ( x_i \right ) \in \mathcal{C}(0)} A_i' \right ) \cup \left ( \bigcup_{B_{r_i} \left ( x_i \right ) \in \mathcal{C}(1)} B_{r_i} \left ( x_i \right )  \right ).
$$
Proceeding with this decomposition inductively for balls in $\mathcal{C}(1)$, ..., $\mathcal{C}(k)$, we obtain a collection of set $\left \{ A_i' \right \}$ exhausting $D$. 
We note that for each $B_{r_i} \left ( x_i \right ) \in \mathcal{C}$, either $r_i \leq s$, due to the inductive covering procedure above, 
or $A_i' = \left ( D \cap B_{r_i} \left ( x_i \right ) \right ) \backslash F_i$, and by the definition of $F_i$ in Lemma \ref{ICLemma}, \eqref{RadiusDichotomyR},
$$
\sup \left \{ I_\phi^{c(A)} \left ( y, \rho r_i \right ) \, : \, y \in A_i' \right \} \leq U - \delta.
$$
In the latter case, we cover $A_i'$ by balls $B_{\rho s_i} \left ( x_{ij} \right )$ replacing $B_{s_i} \left (x_i \right )$,
and replace $A_i$ with $A_{ij} = B_{\rho s_i}  \left ( x_{ij} \right ) \cap A_i'$, on which \eqref{FreqDropCondition} clearly holds.
Noting that $A_i' \subset B_{s_i} \left ( x_i \right )$,  it suffices to choose $C(n) \rho^{-n}$ balls of radius $\rho s_i$ to cover $A_i'$, and since we have already fixed $\rho \in \left ( 0, \rho_0(n) \right )$,
this modification of $\mathcal{C}$ comes at the expense of updating $2C_R(n)$ in \eqref{TowardsDiscreteSum} by another dimensional constant $C(n)$.

Finally, observing that some balls in $\mathcal{C}$ have radii $s_i \in \left [ 10 \rho s, s \right )$, we replace them with concentric balls of radius $s$.
Once again this replacement results in an updated constant $C = C(\rho,n) = C(n)$ on the right-hand side of \eqref{TowardsDiscreteSum}. 
Hence, we conclude that the updated $\mathcal{C} = \left \{ B_{s_i} \left ( x_i \right ) \right \}$ and corresponding $A_i$ satisfy \eqref{RadiusLB}, \eqref{RadiusDichotomy} and \eqref{DiscreteSum} in Lemma \ref{ECLemma}. 
\end{proof}

Finally, we prove Lemma \ref{ICLemma}. The proof is an adaptation of the proof of \cite[Lemma 7.3]{DLMSV}, with the additional task of determining a suitably small scale $\tilde{r} > 0$ so that our construction
satisfies the conditions \eqref{RadiusLBR}, \eqref{RadiusDichotomyR} and \eqref{DiscreteSumR} in the statement.

\begin{proof}[Proof of Lemma \ref{ICLemma}]
We fix $0 < \rho \leq 1/100$.
By translation and scaling, we also fix $x = 0$ and $\tau = \tilde{r}$. 
Thanks to Lemma \ref{ElementaryUpperBounds}, all relevant quantities in our analysis remain uniformly bounded in terms of $\mathcal{I}$ even after such a change of coordinates,
while $\lambda_M$ even shrinks, due to the favorable scaling of \eqref{Functional}.
In addition, we initially treat $\tilde{r} \in \left (0, \min \left \{ R \left (n, N, \lambda_M \right ),  \mathrm{dist} \left ( K , \partial \Omega \right ) \right \}  \right )$, 
(and consequently $0 < \sigma < \tilde{r}$, and $D \subset \mathcal{S}(u) \cap K \cap B_{\tilde{r}}(x)$), as well as $\delta > 0$, as fixed.

In order to construct a finite cover of $D$ satisfying \eqref{RadiusLBR}, \eqref{RadiusDichotomyR} and \eqref{DiscreteSumR} as claimed,
we will update our smallness requirements on $\tilde{r}$ and $\delta$ in the course of the proof.
In particular, for any $\tilde{r}$ and $\delta$ sufficiently small, we will first construct a finite cover of $D$ satisfying \eqref{RadiusLBR} and \eqref{RadiusDichotomyR}. 
Finally, we will derive a possibly more stringent smallness requirement on $\delta$ which will ensure the bound \eqref{DiscreteSumR}.

We denote the smallest integer satisfying $(10\rho)^\kappa \tilde{r} \leq \sigma$ as $\kappa$, 
and begin with inductively constructing a finite cover of $D$ satisfying \eqref{RadiusLBR} and \eqref{RadiusDichotomyR} in $\kappa$ steps.
At each step $k = 0,...,\kappa$, we will cover $D$ by balls $\mathcal{C}(k) = \left \{ B_{\rho_i} \left ( x_i \right ) \, : \, i \in I_k \right \}$,
where $\mathcal{C}(0) = \left \{ B_{\tilde{r}} (0) \right \}$, satisfying the following conditions:
\begin{enumerate}[(i)]
 \item $B_r(x) \in \mathcal{C}(k)$ implies $r = \left ( 10 \rho \right )^{j} \tilde{r}$ for some $j=0,...,k$. \label{RadiusRequirement}
 \item $B_r(x)$, $B_{r'} \left ( x' \right ) \in \mathcal{C}(k)$ implies $B_{r/5}(x) \cap B_{r'/5} \left ( x' \right ) = \emptyset$. \label{VitaliRequirement}
 \item If $j < k$ in \eqref{RadiusRequirement} for $B_r(x) \in \mathcal{C}(k)$, then $B_r(x) \in \mathcal{C}(k+1)$ as well. \label{KeepBadBalls}
\end{enumerate}

{\bf{Inductive procedure.}} 
Consider $B_r(x) \in \mathcal{C}(k)$. By \eqref{RadiusRequirement} and \eqref{KeepBadBalls}, $B_r(x) \in \mathcal{C}(k+1)$, if $r = (10\rho)^j \tilde{r}$ for some $j < k$.
If $j=k$ instead, then we consider the set 
\begin{equation}
F = F \left ( B_r(x) \right ) = D \cap B_r(x) \cap \left \{ y \, : \, I_\phi^{c(A)} (y, \rho r) > U - \delta \right \}. \label{PinchedSet}
\end{equation}

{\emph{Case 1:}} 
$F$ fails to $\rho r$-span an $(n-2)$-dimensional affine subspace. 
We denote such balls of radius $r = (10\rho)^k \tilde{r}$ in $\mathcal{C}(k)$ as $\mathcal{B}(k)$ and add any such $B_r(x)$ to $\mathcal{C}(k+1)$.
Note that in this case there exists an $(n-3)$ dimensional affine subspace $L = L \left (B_r(x) \right )$ such that $F \subset B_{\rho} (L)$ by Remark \ref{UsefulGeometricRemark}.

{\emph{Case 2:}}
$F$ $\rho r$-spans an $(n-2)$-dimensional affine subspace $V$. 
We denote such balls of radius $r = (10\rho)^k \tilde{r}$ in $\mathcal{C}(k)$ as $\mathcal{G}(k)$.
By Lemma \ref{TubularNeighborhoodLemma}, 
choosing firstly $\tilde{r} = \tilde{r} \left (n, N, \lambda_M \right )$ and secondly $\delta = \left ( n, \lambda_M, \mathcal{I}, \rho, \tilde{r} \right )$ small enough, 
we conclude that $D \cap B_r(x) \subset B_{\rho r}(V)$.

Considering all the balls $B^i \in \mathcal{G}(k)$ and corresponding affine spaces $V_i$, we define the set
$$
G(k) = D \cap \bigcup_i B_{2 \rho r} \left ( V_i \right ).
$$
Note that we can cover $G(k)$ with a collection of balls with radius $\left ( 10 \rho \right )^{k+1}$ such that the corresponding concentric balls of radii $2 \rho (10 \rho)^k$ are pairwise disjoint,
and the centers are contained in $D \cap \left ( \cup_i B^i \cap V_i \right )$. We denote this collection as $\mathcal{F}(k+1)$ and note that its cardinality is bounded by a dimensional constant $C(n)$.
We add $B \in \mathcal{F}(k+1)$ to $\mathcal{C}(k+1)$, 
if and only if $B$ does not intersect any ball in $\mathcal{B}_{1/5}(k)$,
where $\mathcal{B}_{1/5}(k)$ is the collection of balls obtained from $\mathcal{B}(k)$ by scaling the radius of each individual ball by $1/5$.

It is not difficult to check that composed with balls added from $\mathcal{B}(k)$ and $\mathcal{F}(k+1)$, the collection $\mathcal{C}(k+1)$ still covers $D$, 
while satisfying \eqref{RadiusRequirement}, \eqref{VitaliRequirement}.
Recall that we have already introduced our first smallness requirements firstly on $\tilde{r} = \tilde{r} \left (n, N, \lambda_M\right )$ 
and secondly on $\delta = \left ( \tilde{r}, \rho, n, \lambda_M, \mathcal{I}, \right )$,
based on our application of Lemma \ref{TubularNeighborhoodLemma}.

{\bf{Frequency pinching requirement.}}
We claim that given any $\eta > 0$, we can update the above-introduced smallness requirements on $\tilde{r}$ and $\delta$ accordingly so that either $\mathcal{C}(k) = \left \{ B_{\tilde{r}}(0) \right \}$, 
or for every $B_s(x) \in \mathcal{C}(k)$,
\begin{equation}
I_\phi^{c(A)} \left ( x, \rho s/5 \right ) \geq U - \eta.  \label{PinchingRequirement}
\end{equation}

We note that unless the inductive procedure leads to the unrefined cover $\mathcal{C}(k) = \left \{ B_{\tilde{r}}(0) \right \}$, 
for every $B_s(x) \in \mathcal{C}(k)$, $s = (10 \rho)^{j+1} \tilde{r}$, and its center $x$ lies in $V \cap B'$ for some
$B' \in \mathcal{C}(j)$ with radius $(10 \rho)^j \tilde{r}$ and some $(n-2)$-dimensional affine subspace $V$ which is $\left [(10 \rho)^{j} \rho \tilde{r} \right ]$-spanned by the set $F \left ( B' \right )$
defined in \eqref{PinchedSet}. We remark that $B' \cap V$ contains at least one point $z \in F \left ( B' \right ) \subset \mathcal{S}(u) \cap K$.
Furthermore, it follows from the definition of bound $U$,the monotonicity of $I_\phi^{c(A)}$ and \eqref{PinchedSet} that at the points in $F \left ( B' \right )$, 
the pinching of additively modified frequency between scales $s$ and $\tilde{r}$ is less than $\delta$.
Therefore, we can apply Lemma \ref{FOSSControl} with $\delta = \delta \left (\eta, \rho, \tilde{r}, n, \lambda_M, \mathcal{I} \right )$ chosen sufficiently small and get
$$
\left | I_\phi \left ( x, \rho s/5 \right ) - I_\phi(z,s) \right | \leq \eta/3.
$$
Combining this estimate with $I_\phi(z,s)^{c(A)} \geq U - \delta$, we can easily verify that
$$
I_\phi^{c(A)} \left ( x, \rho s/5 \right ) \geq U - \delta - \eta/3 - c(A) \left ( 1 - \rho / 5 \right ) s^2 \geq U - \delta - \eta/3 - c(A) \tilde{r}^2.
$$
Since $A = A \left ( n, \lambda_M, N, \mathcal{I} \right )$, firstly updating our smallness requirement on $\tilde{r}$ based on $\eta$, $n$, $\lambda_M$, $N$ and $\mathcal{I}$,
and then updating the smallness requirement on $\delta = \delta \left (\eta, \tilde{r}, \rho, n, \lambda_M, \mathcal{I} \right )$ once again, we obtain
$$
I_\phi^{c(A)} \left ( x, \rho s/5 \right ) \geq U - \eta/3 - \eta/3 - \eta/3 = U - \eta,
$$
as we have claimed.

{\bf{Packing bound.}} 
For every sufficiently small $\tilde{r} = \tilde{r} \left (\eta, \rho, n, N, \lambda_M, \mathcal{I} \right )$ 
and every sufficiently small $\delta = \delta \left ( \eta, \tilde{r}, \rho, n, \lambda_M, \mathcal{I}  \right )$,
we have obtained a covering $\mathcal{C}(k)$ of the set $D$, which itself shrinks as $\tilde{r}$ decreases by definition, 
and this cover satisfies \eqref{RadiusLBR} and \eqref{RadiusDichotomyR} in the statement of lemma.
Our final task is to optimize $\eta$, and update $\tilde{r}$ and $\delta$ respectively, in order to satisfy the condition \eqref{DiscreteSumR} as well.
In other words, we seek to prove the packing bound,
$$
\sum_{B_s(x) \in \mathcal{C}(k)} s^{n-2} \leq C_R(n) \tilde{r}^{n-2}.
$$
Enumerating the balls in $\mathcal{C}(k)$ as $\left \{ B_{5 s_i} \left ( x_i \right ) \right \}_{i \in I}$, we recast this packing bound as the measure estimate,
\begin{equation}
\mu \left ( B_{\tilde{r}} (0) \right ) \leq \frac{C_R(n)}{5^{n-2}} \tilde{r}^{n-2}, \quad \mathrm{where} \quad \mu = \sum_{i \in I} s_i^{n-2} \delta_{x_i}. \label{MeasureEstimate}
\end{equation}
As we will prove \eqref{MeasureEstimate} by an induction argument, we also introduce the truncated version of measure $\mu$,
$$
\mu_s = \sum_{i \in I, s_i \leq s} s_i^{n-2} \delta_{x_i}.
$$
By the construction of $\mathcal{C}(k)$, $\mu_s$ satisfies the following properties:
\begin{enumerate}[(i)]
 \item $\mu_t \leq \mu_\tau$ for $t \leq \tau$. \label{MassMonotonicity}
 \item $\mu = \mu_{\tilde{r}/5}$. \label{MassUpper}
 \item $\mu_s = 0$ for $s < \overline{r}$, where $\overline{r} = (10 \rho)^\kappa \tilde{r}$. \label{MassLower}
\end{enumerate}

We note that for $\chi = \log_2 \left ( \tilde{r} / \overline{r} \right ) - 8$, proving
\begin{equation}
\mu_s \left ( B_s (x) \right ) \leq C_R(n) s^{m-2}, \quad \forall x \in B_{\tilde{r}}(0), \quad \forall s = 2^j \overline{r} \quad \mathrm{with} \quad j = 0,1,...,\chi, \label{StepBound}
\end{equation}
will yield \eqref{MeasureEstimate}. 
Indeed \eqref{StepBound} implies $\mu \left ( B_{\tilde{r}/128} (x) \right ) \leq C_R(n) \tilde{r}^{n-2}$ for every $x \in B_{\tilde{r}}(0)$.
Therefore, covering $B_{\tilde{r}}(0)$ with $C(n)$ balls of radius $\tilde{r}/128$ will give \eqref{MeasureEstimate} with a modified constant $C_R(n)$.
Hence, we will prove \eqref{StepBound} by induction on $j$. 

{\emph{Base case:}} $j=0$, that is $s = \overline{r}$. We have
$$ \mu_{\overline{r}} \left ( B_{\overline{r}}(0) \right ) = \mathcal{N} \left ( x, \overline{r} \right ) \overline{r}^{n-2}, $$
where $\mathcal{N} \left ( x, \overline{r} \right )$ is the number of balls $B_{s_i} \left ( x_i \right )$ with $s_i = \overline{r}$ and $x_i \in B_{\overline{r}}(x)$.
By the Vitali property \eqref{VitaliRequirement} of $\mathcal{C}(k)$, such balls are pairwise disjoint. 
Since they are also contained in $B_{2 \overline{r}}(x)$, clearly $\mathcal{N} \left ( x, \overline{r} \right ) \leq 2^n$, and \eqref{StepBound} holds for $j=0$. 

{\emph{Induction hypothesis:}} The estimate \eqref{StepBound} holds for some $j < \chi$ with a dimensional constant $C_R(n)$.

Thus, proving \eqref{StepBound} for $j+1$ with the same dimensional constant $C_R(n)$ will complete the induction argument and the proof of \eqref{MeasureEstimate}.
Hence, we set $r = 2^j \overline{r}$ and seek to prove \eqref{StepBound} for $s = 2r$.
We remark that decomposing $\mu_{2r}$ as
$$
\mu_{2r} = \mu_r + \mu_{r,2r}, \quad \mathrm{where} \quad \mu_{r,2r} = \sum_{i \in I, r < s_i \leq 2r} s_i^{n-2} \delta_{x_i},
$$
we can estimate $\mu_r \left ( B_{2r}(x) \right )$ by covering $B_{2r}(x)$ by $2^n$ balls of radius $r$ and using the induction hypothesis,
while $\mu_{r,2r} \left (B_{2r}(x) \right ) = \mathcal{N}(x,2r) (2r)^{n-2}$, 
where $\mathcal{N}(x,2r)$ is the number of balls $B_{s_i} \left ( x_i \right )$ with $x_i \in B_{2r}(x)$ and $r < s_i \leq 2r$.
$\mathcal{N}(x,2r)$ bounded by a dimensional constant $C(n)$,
as these balls are pairwise disjoint and contained in $B_{3r}(x)$ and have radii greater than $r$.
Therefore, for every $x \in B_{\tilde{r}}(0)$, we have the {\emph{coarse bound}},
\begin{equation}
\mu_{2r} \left ( B_{2r}(x) \right ) \leq C(n) C_R(n) (2r)^{n-2},
\label{CoarseBound}
\end{equation}
which will be useful on its own, though it does not suffice to complete the induction argument.

\subsection{Improved measure bound} \label{DRSubsection}
We will complete the proof of Lemma \ref{ICLemma} by upgrading \eqref{CoarseBound} for every $x \in B_{\tilde{r}}(0)$ to
\begin{equation}
\mu_{2r} \left ( B_{2r}(x) \right ) \leq C_R(n) (2r)^{n-2}.
\label{ImprovedBound}
\end{equation}
Therefore, fixing $x \in B_{\tilde{r}}(0)$, we denote the restriction of $\mu_{2r}$ to $B_{2r}(x)$ as $\overline{\mu}$,
and seek to prove $\overline{\mu} \left ( B_{2r}(x) \right ) \leq C_R(n) (2r)^{n-2}$.

As naive covering and induction arguments are not enough to obtain the sharp bound \eqref{ImprovedBound},
we will resort to the following Discrete Reifenberg Theorem of Naber and Valtorta from \cite{NV}, which takes into account 
the mean-flatness of $\mathrm{spt} \left ( \overline{\mu} \right )$ averaged over scales and space in a scale-invariant fashion.

\begin{theorem}[Naber-Valtorta,\cite{NV}] \label{DiscreteReifenbergTheorem}
Let $\left \{ B_{s_j} \left ( x_j \right ) \right \}_{j \in J}$ be a collection of pairwise disjoint balls contained in $B_{4r}(x) \subset \mathbf{R}^n$ with $x_j \in B_{2r}(x)$,
and for an integer $k \leq n$, let the measure $\mu$ be defined as
$$
\mu = \sum_{j \in J} s_j^k \delta_{x_j}.
$$
There exist positive dimensional constants $\delta_R = \delta_R(n)$ and $C_R = C_R(n)$ such that if the bound,
\begin{equation}
\int_{B_t(y)} \left [ \int_0^t D_\mu^k (z,s) \frac{\mathrm{d}s}{s} \right ] \, \mathrm{d}\mu(z) < \delta_R^2 t^k,
\label{NVHypothesis}
\end{equation}
holds for every $B_t(y) \subseteq B_{4r}(x)$ with $y \in B_{2r}(x)$, then the following estimate holds:
\begin{equation}
\mu \left ( B_{2r}(x) \right ) = \sum_{j \in J} s_j^k \leq C_R (2r)^k.
\label{NVBound}
\end{equation}
\end{theorem}
Hence, in order to obtain \eqref{ImprovedBound}, it suffices to check \eqref{NVHypothesis} with $k=n-2$ for every $B_t(y) \subseteq B_{4r}(x)$ with $y \in B_{2r}(x)$.
This will be accomplished by combining the induction hypothesis \eqref{StepBound}, the coarse bound \eqref{CoarseBound}, the pointwise bound \eqref{JonesNumberEstimate} on $D_\mu^{n-2}(z,s)$, 
the pinching requirement \eqref{PinchingRequirement},
choosing $\eta$ sufficiently small, and finally updating $\tilde{r}$ and $\delta$ accordingly.
Note that without loss of generality, we can assume \eqref{PinchingRequirement} for $x = x_i$, $s = s_i$, that is
\begin{equation}
I_\phi^{c(A)} \left ( x_i , \rho s_i\right ) \geq U - \eta, \label{PinchingAssumed}
\end{equation}
since otherwise, $\mathcal{C}(k) = \left \{ B_{\tilde{r}}(0) \right \}$, and \eqref{MeasureEstimate} holds trivially.

Recall that $c(A) = 3 + A + A^2$ for $A$ as in Corollary \ref{LinMonotonicitySmoothedAdditive}.
In order to utilize \eqref{JonesNumberEstimate} in verifying \eqref{NVHypothesis}, we define the truncated frequency pinching,
\begin{equation*}
\overline{W}_s \left ( x_i \right ) =  
\begin{cases}
I_\phi^{c(A)} \left ( x_i , 32s \right ) - I_\phi^{c(A)} \left ( x_i ,s \right ), \quad & \mathrm{if} \; s > s_i, \\
0 \quad & \mathrm{otherwise}.
\end{cases}
\end{equation*}
Observing that for $s < s_i$, $\mathrm{spt} \left ( \mu \right ) \cap B_s \left ( x_i \right ) = \left \{ x_i \right \}$, and applying \eqref{JonesNumberEstimate}, 
for every $s \in \left (0, 4r \right )$, we obtain
\begin{equation}
D_{\overline{\mu}}^{n-2} \left ( x_i , s \right ) \leq \frac{C}{s^{n-2}} \int_{B_s \left ( x_i \right )} \overline{W}_s (y) \, \mathrm{d}\overline{\mu}(y) , \label{TruncatedJonesEstimate}
\end{equation}
where $C = C \left ( n, N, \lambda_M, \mathcal{I} \right )$.
Fixing $t \leq \mathrm{dist} \left ( y,\partial B_{4r}(x) \right )$, from \eqref{TruncatedJonesEstimate} and Fubini's theorem, we get
\begin{equation}
\begin{aligned}
\int_{B_t(y)} \left [ \int_0^t D_{\overline{\mu}}^{n-2} (z,s) \frac{\mathrm{d}s}{s} \right ] \mathrm{d}\overline{\mu}(z)
& \leq C \int_{B_t(y)} \int_0^t s^{1-n} \int_{B_s(z)} \overline{W}_s (\zeta) \mathrm{d} \overline{\mu}(\zeta) \mathrm{d}s \mathrm{d} \overline{\mu}(z) \\
& = C \int_0^t s^{1-n} \int_{B_t(y)} \int_{B_s(z)} \overline{W}_s (\zeta) \mathrm{d} \overline{\mu}(\zeta) \mathrm{d} \overline{\mu}(z) \mathrm{d}s.
\end{aligned} \label{TruncatedJonesApplied}
\end{equation}
Taking into consideration the supports of $\overline{\mu}$, $\mu_s$, the pairwise disjointness of $B_{s_i} \left ( x_i \right )$ in $\mathcal{C}(k)$, 
and the definition of $\overline{W}_s \left ( x_i \right )$,
we can first restrict the domains of integrations to $B_{2r}(x)$ and then substitute $\overline{\mu}$ with $\mu_s$ in \eqref{TruncatedJonesApplied}.
Doing so and applying Fubini's theorem again, for $I$ denoting the left-hand side of \eqref{TruncatedJonesApplied}, we have
$$
I \leq C \int_0^t s^{1-n} \int_{B_{t+s}(y) \cap B_{2r}(x)} \overline{W}_s ( \zeta ) \int_{B_s(\zeta) \cap B_{2r}(x)} \mathrm{d}\mu_s(z) \mathrm{d}\mu_s(\zeta)\mathrm{d}s.
$$
In order to bound the integral with respect to $z$ from above by $C(n) s^{n-2}$, we use the induction hypothesis \eqref{StepBound} in the case $s \leq r$, 
and the covering argument in the derivation of coarse bound \eqref{CoarseBound} combined with \eqref{StepBound} again in the case $r < s \leq 4r$. 
As a result,
\begin{equation}
\begin{aligned}
I \leq C \cdot C(n) \int_0^t \int_{B_{t+s}(y) \cap B_{2r}(x)} \overline{W}_s ( \zeta ) \mathrm{d}\mu_s(\zeta)\frac{\mathrm{d}s}{s}
& \leq \tilde{C} \int_0^t \int_{B_{2t}(y)} \overline{W}_s ( \zeta ) \mathrm{d}\mu_t \frac{\mathrm{d}s}{s} \\
& \leq \tilde{C} \int_{B_{2t}(y)} \int_0^t \overline{W}_s ( \zeta ) \frac{\mathrm{d}s}{s} \mathrm{d}\mu_t,
\end{aligned} \label{PenultimateEst}
\end{equation}
where the second inequality follows from $s \leq t$, inclusion and the monotonicity of $\mu_s$ in $s$, 
the third inequality is due to Fubini's theorem,
and $\tilde{C} = \tilde{C} \left ( n, N, \lambda_M, \mathcal{I} \right )$.

Next we will estimate the inner integral on the right-hand side of \eqref{PenultimateEst} for fixed $\zeta \in \mathrm{spt} \left ( \mu_t \right )$,
that is $\zeta = z_i$ for some $i$ such that $s_i > s$, as $\overline{W}_s \left ( z_i \right ) = 0$ otherwise. Consequently,
\begin{equation}
\int_0^t \overline{W}_s ( \zeta ) \frac{\mathrm{d}s}{s} 
= \int_{s_i}^t \overline{W}_s \left ( z_i \right ) \frac{\mathrm{d}s}{s} 
= \int_{s_i}^t \left [ I_\phi^{c(A)} \left ( z_i, 32s \right ) - I_\phi^{c(A)} \left ( z_i, s \right ) \right ] \frac{\mathrm{d}s}{s}.
\label{ScaleIntegral1}
\end{equation}
Having fixed $s_i$ and choosing $\kappa$ to be the smallest integer such that $2^\kappa s_i \geq t$, we note that
$t \leq 4r$, $r = 2^j \overline{r}$, $j < \chi = \log_2 \left ( \tilde{r} / \overline{r} \right ) - 8$ together imply
\begin{equation}
32 \cdot 2^{\kappa+1} s_i = 128 \cdot 2^{\kappa-1} s_i < 128 t \leq 512 r \leq \frac{512}{2} \overline{r} 2^{\chi} = \tilde{r}. \label{ScaleUB}
\end{equation}
Then using \eqref{ScaleIntegral1} and the monotonicity of $I_\phi^{c(A)} \left (z_i, \tau \right )$, we can estimate
\begin{equation}
\begin{aligned}
\int_0^t \overline{W}_s ( \zeta ) \frac{\mathrm{d}s}{s}
& \leq \sum_{k=0}^\kappa \int_{2^k s_i}^{2^{k+1} s_i} \left [ I_\phi^{c(A)} \left ( z_i, 32s \right ) - I_\phi^{c(A)} \left ( z_i, s \right ) \right ] \frac{\mathrm{d}s}{s} \\
& \leq \sum_{k=0}^\kappa \left [ I_\phi^{c(A)} \left ( z_i, 2^{k+6} s_i \right ) - I_\phi^{c(A)} \left ( z_i, 2^k s_i \right ) \right ] \int_{2^k s_i}^{2^{k+1} s_i} \frac{\mathrm{d}s}{s} \\
& = \log 2 \sum_{k=0}^\kappa \left [ I_\phi^{c(A)} \left ( z_i, 2^{k+6} s_i \right ) - I_\phi^{c(A)} \left ( z_i, 2^k s_i \right ) \right ].
\end{aligned}
\label{ScaleIntegral2}
\end{equation}
Note that we can write
\begin{equation}
\begin{aligned}
\sum_{k=0}^\kappa & \left [ I_\phi^{c(A)} \left ( z_i, 2^{k+6} s_i \right ) - I_\phi^{c(A)} \left ( z_i, 2^k s_i \right ) \right ] \\
& = \sum_{\ell = 0}^5 \sum_{k=0}^\kappa \left [ I_\phi^{c(A)} \left ( z_i, 2^{k+\ell+1} s_i \right ) - I_\phi^{c(A)} \left ( z_i, 2^{k+\ell} s_i \right ) \right ] \\
& = \sum_{\ell = 0}^5 \left [ I_\phi^{c(A)} \left ( z_i, 2^{\kappa + \ell + 1} s_i \right ) - I_\phi^{c(A)} \left ( z_i, 2^{\ell} s_i \right ) \right ],
\end{aligned}
\label{ScaleIntegral3}
\end{equation}
By \eqref{ScaleIntegral2}, \eqref{ScaleIntegral3}, the monotonicity of $I_\phi^{c(A)} \left (z_i, \tau \right )$, \eqref{ScaleUB}, the definition of $U$, and \eqref{PinchingAssumed}, we obtain the estimate
\begin{equation}
\int_0^t \overline{W}_s ( \zeta ) \frac{\mathrm{d}s}{s}
\leq 6 \log 2 \left [ I_\phi^{c(A)} \left ( z_i, \tilde{r} \right ) - I_\phi^{c(A)} \left ( z_i, \rho s_i \right ) \right ] \leq \left ( 6 \log 2 \right ) \eta.
\label{ScaleIntegralComplete}
\end{equation}

From \eqref{PenultimateEst} and \eqref{ScaleIntegralComplete} we get
\begin{equation*}
\int_{B_t(y)} \left [ \int_0^t D_{\overline{\mu}}^{n-2} (z,s) \frac{\mathrm{d}s}{s} \right ] \mathrm{d}\overline{\mu}(z) \leq \tilde{C} \left ( 6 \log 2 \right ) \cdot \eta \cdot \mu_t \left ( B_{2t}(y) \right ).
\end{equation*}
In addition, by covering $B_{2t}(y)$ with balls of radius $t$, and using the induction hypothesis \eqref{StepBound} in the case $t \leq r$, while using the coarse estimate \eqref{CoarseBound} in the case $ r < t \leq 4r$,
we can estimate
$$
\mu_t \left ( B_{2t}(y) \right ) \leq C(n) t^{n-2},
$$
and conclude that
\begin{equation}
\int_{B_t(y)} \left [ \int_0^t D_{\overline{\mu}}^{n-2} (z,s) \frac{\mathrm{d}s}{s} \right ] \mathrm{d}\overline{\mu}(z) \leq C \left ( n, N, \lambda_M, \mathcal{I} \right ) \eta t^{n-2}.
\label{UltimateEstimate}
\end{equation}
Once we choose $\eta = \eta \left ( n, N, \lambda_M, \mathcal{I} \right )$ small enough, and accordingly update our smallness requirements on $\tilde{r}$ and $\delta$ respectively,
\eqref{UltimateEstimate} implies \eqref{NVHypothesis} for the measure $\overline{\mu} = \mu_{2r} \rest B_{2r}(x)$. 
Hence, Theorem \ref{DiscreteReifenbergTheorem} gives \eqref{ImprovedBound}, and the proof of this covering lemma is complete.
\end{proof}

\section{Rectifiability} \label{RectifiabilitySection}

In this section, following the strategy of \cite[Section 8]{DLMSV}, we prove Theorem \ref{MainTheorem2}.
The main ingredients are Theorem \ref{MainTheorem1} and a characterization of rectifiable measures by Azzam and Tolsa \cite{AT15}.
Below we denote the $k$-dimensional Hausdorff measure as $\mathcal{H}^k$.

\begin{theorem}[Azzam-Tolsa, \cite{AT15}] \label{ATTheorem}
Let $\mathcal{S} \subset \mathbf{R}^n$ be  $\mathcal{H}^k$-measurable with $\mathcal{H}^k ( \mathcal{S} ) < \infty$ and $\mu = \mathcal{H}^k \rest \mathcal{S}$.
Then $\mathcal{S}$ is countably $k$-rectifiable if and only if
\begin{equation}
\int_0^1 D_\mu^k (x,s) \frac{\mathrm{d}s}{s} < \infty \quad \mathrm{for} \quad \mu-\mathrm{almost} \; \mathrm{every} \; x \in \mathcal{S}.
\label{FinitenessHypothesis} 
\end{equation}
\end{theorem}

We remark that if the upper limit of integral in \eqref{FinitenessHypothesis} is replaced with any finite number, the claim is still valid, as evidently only the fine scales matter.

\begin{proof}[Proof of Theorem \ref{MainTheorem2}]
We fix a compact subset $K \subset \Omega$, $x \in K$, and the scale $r \leq \tilde{r} \left (n, N, \lambda_M, \mathcal{I} \right )$ as in Theorem \ref{MTETheorem}.
We also denote $\mathcal{S} = \mathcal{S}(u) \cap K \cap B_{r}(x)$.

We observe that for any $\delta \in \left ( 0, r \right )$ fixed
and $\mathcal{C} = \left \{ B_{\delta} \left ( x_i \right ) \right \}$, a maximal, piecewise disjoint collection of open balls with centers in $\mathcal{S}$,
by Theorem \ref{MTETheorem},
\begin{equation}
\mathcal{N}_\delta \omega_n  \delta^n = \mathcal{L} \left ( B_\delta \left ( \mathcal{S} \right ) \right ) \leq C r^{n-2} \rho^2, \label{FiniteMass1}
\end{equation}
where $\mathcal{N}_\delta$ is the number of balls in $\mathcal{C}$ and the constant $C = C \left ( n, N, \lambda_M, \mathcal{I} \right )$ is as in Theorem \ref{MTETheorem}.
By the maximality of $\mathcal{C}$, $\mathcal{S} \subset \cup_{i=1}^{N_\delta} B_{3\delta} \left ( x_i \right )$. 
Thus, by the definition of $\mathcal{H}^{n-2}_\delta$ and \eqref{FiniteMass1},
\begin{equation}
\mathcal{H}_\delta^{n-2} \left ( \mathcal{S} \right ) \leq  \mathcal{N}_\delta \omega_{n-2} \cdot ( 3 \delta )^n \leq C 3^{n-2} \omega_{n-2}\omega_n^{-1} r^{n-2}. \label{FiniteMass2}
\end{equation}
Letting $\delta \downarrow 0$ in \eqref{FiniteMass1}, and setting $\tilde{C} = C 3^{n-2} \omega_{n-2}\omega_n^{-1}$, 
by the definition of $\mu$ and $\mathcal{H}^{n-2}$, we have for every $x \in K$ and $r \leq \tilde{r}$,
\begin{equation}
\mu \left ( B_r (x) \right ) = \mathcal{H}^{n-2} \left ( \mathcal{S} \right ) \leq \tilde{C} r^{n-2}. \label{FiniteMassComplete}
\end{equation}
Hence, it suffices to check \eqref{FinitenessHypothesis} to conclude that $\mathcal{S}$ is countably $(n-2)$-rectifiable, 
which implies the countable $(n-2)$-rectifiability of $\mathcal{S}(u)$ by a simple covering argument.

For $t \leq \tilde{r}/128$, using the pointwise bound \eqref{JonesNumberEstimate} on $D_\mu^{n-2}(z,s)$  and Fubini's theorem as in the estimate \eqref{TruncatedJonesApplied}, we get
\begin{equation}
\begin{aligned}
\int_{B_t(y)} \int_0^t D_{\mu}^{n-2} (z,s) \frac{\mathrm{d}s}{s} \mathrm{d} \mu(z)
& \leq C \int_{B_t(y)} \int_0^t s^{1-n} \int_{B_s(z)} W_{s,32s}^{3+A+A^2} (\zeta) \mathrm{d} \mu (\zeta) \mathrm{d}s \mathrm{d} \mu(z) \\
& = C \int_0^t s^{1-n} \int_{B_t(y)} \int_{B_s(z)} W_{s,32s}^{3+A+A^2} (\zeta) \mathrm{d} \mu(\zeta) \mathrm{d} \mu(z) \mathrm{d}s \\
& \leq C \int_0^t \int_{B_{t+s}(y)} W_{s,32s}^{3+A+A^2} (\zeta) \int_{B_s (\zeta)}  \mathrm{d} \mu(z) \mathrm{d}\mu(\zeta) \frac{\mathrm{d}s}{s^{n-1}}.
\end{aligned} 
\label{RectEst1}
\end{equation}
Using \eqref{FiniteMassComplete} with $r=s$ to estimate the innermost integral on the third line of \eqref{RectEst1}, that is $\mu \left (B_s (\zeta) \right )$, and applying Fubini's Theorem, we have
\begin{equation}
\begin{aligned}
\int_{B_t(y)} \int_0^t D_{\mu}^{n-2} (z,s) \frac{\mathrm{d}s}{s} \mathrm{d} \mu(z)
& \leq C \int_0^t \int_{B_{t+s}(y)} W_{s,32s}^{3+A+A^2} (\zeta) \mathrm{d} \mu (\zeta) \frac{\mathrm{d}s}{s} \\
& \leq C \int_{B_{2t}(y)} \int_0^t W_{s,32s}^{3+A+A^2} (\zeta) \frac{\mathrm{d}s}{s} \mathrm{d} \mu (\zeta).
\end{aligned}
\label{RectEst2}
\end{equation}

Next we claim that for every $t \leq 2^{-7} \tilde{r}$ and $\epsilon \in (0,t)$,
\begin{equation}
\int_\epsilon^t W_{s,32s}^{3+A+A^2} (\zeta) \frac{\mathrm{d}s}{s} 
\leq 6 \log 2 \left [ I_\phi^{3+A+A^2} \left ( \zeta, \tilde{r} \right ) - I_\phi^{3+A+A^2} \left ( \zeta, \epsilon \right ) \right ].
\label{EpsilonEstimate}
\end{equation}
As in the proof of \eqref{ScaleIntegralComplete}, we set $\kappa$ to be the smallest integer such that $2^\kappa \epsilon \geq t$.
Then $t \leq \tilde{r}/128$ implies $32 \cdot 2^{\kappa+1} \epsilon \leq \tilde{r}$.
Therefore, $\left [ \epsilon, t \right ] \subset \cup_{k=0}^\kappa \left [ 2^k \epsilon, 2^{k+1} \epsilon \right ]$, and
we can break up the integral with respect to $s$ into the corresponding dyadic pieces. 
Then arguing exactly as in \eqref{ScaleIntegral2}, \eqref{ScaleIntegral3} and exploiting the monotonicity of $I_\phi^{3+A+A^2}(\zeta,s)$, we get the estimate \eqref{EpsilonEstimate}.
Lastly, letting $\epsilon \downarrow 0$ in \eqref{EpsilonEstimate}, we have
\begin{equation}
\int_0^t W_{s,32s}^{3+A+A^2} (\zeta) \frac{\mathrm{d}s}{s} 
\leq 6 \log 2 \left [ I_\phi^{3+A+A^2} \left ( \zeta, \tilde{r} \right ) - I_\phi^{3+A+A^2} \left ( \zeta, 0^+ \right ) \right ] \leq C,
\label{ScaleAvEstimate}
\end{equation}
where $C = C \left ( n, N, \lambda_M , \mathcal{I} \right )$, since $\tilde{r}$ and $A$ are bounded uniformly from above by constants depending on $n$, $N$, $\lambda_M$, $\mathcal{I}$,
while by Remark \ref{ScaleRemark}, $I_\phi \left ( \zeta, \tilde{r} \right ) \leq \mathcal{I}$ for every $\zeta \in K$.

Combining \eqref{RectEst2}, \eqref{ScaleAvEstimate}, and \eqref{FiniteMassComplete} with $r = 2t$, we conclude that whenever $t \leq 2^{-7} \tilde{r}$,
$$
\int_{B_t(y)} \int_0^t D_{\mu}^{n-2} (z,s) \frac{\mathrm{d}s}{s} \mathrm{d} \mu(z) \leq C,
$$
for a new constant $C = C \left ( n, N, \lambda_M , \mathcal{I} \right )$.
Thus, by Fubini's Theorem for $\mu$-almost every $z \in \mathcal{S}$, we have
$$
\int_0^t D_\mu^{n-2} (z,s) \frac{\mathrm{d}s}{s} < \infty.
$$
Therefore, by Theorem \ref{ATTheorem}, $\mathcal{S}$ is countably $(n-2)$-rectifiable.
Finally, exhausting $\mathcal{S}(u)$ by its compact subsets and covering each compact subsets by balls in which we have verified rectifiability, 
we conclude that $\mathcal{S}(u)$ is countably $(n-2)$-rectifiable.
\end{proof}

\section{Partitions on manifolds and maps into homogeneous trees} \label{FinalRemarksSection}

In this section we briefly discuss how the proofs of Theorems \ref{MainTheorem1} and \ref{MainTheorem2} imply Theorem \ref{MainTheorem3} as well.

\begin{proof}[Proof of Theorem \ref{MainTheorem3}]
Firstly, we discuss the case of problem \eqref{Manifolds}, that is the partitioning of a smooth, bounded $\Omega \subseteq M$ 
into a collection of $N$ open, connected and pairwise disjoint subsets, $\left \{ \Omega_k \right \}$,
that minimizes \eqref{EigenvalueMinimization} with respect to $\lambda_k = \lambda_1 \left ( \Omega_k \right )$, the first Dirichlet eigenvalues of $\Omega_k$ with zero boundary data and 
with respect to the Laplace-Beltrami operator $\Delta_g$ on $(M,g)$. We note that the case $\Omega = M$ is possible and rather interesting, whether $\partial M = \emptyset$ or not.

The only modification we need to introduce in this case concerns the variational setup. In particular, the variational quantities such as Dirichlet energy, height, $L^2$-norm, and frequency
in Definition \ref{ClassicalQuantities} should be replaced with the analogous counterparts for a Riemannian manifold $M$ with metric $g$.
Consequently, the constants appearing in the variational formulas in Section \ref{PrelimSection} will depend on this metric as well. We refer to \cite[Section 2]{GarofaloLin} for the details.
Note that as an additional techicality, we have to restrict ourselves to the radius of injectivity of $M$, as we use the exponential map $\exp_p \, : \, T_p M \to M$ to define the domain variations in \eqref{CombinedVariations}.
Once we adapt the classical variational setting to $(M,g)$, it is straightforward to define the analogous smoothed variational quantities in Section \ref{SmoothSection}.
Hence, the remaining sections carry over veribatim to the manifold setting, while the constants appearing in the estimates depend on the metric $g$.

Finally, we point out that problem \eqref{Trees} is indeed a simpler special case. 
Instead of considering the constrained minimization problem \eqref{EquivalentProblem}, we consider $u \, : \, \Omega \to \Sigma_N$
minimizing the Dirichlet energy in $\Omega$ with respect to some boundary data $g \in H^{1/2} \left ( \partial \Omega, \Sigma_N \right )$.
Many aspects of the problem such as regularity and compactness simplify, and all the variational formulas in Sections \ref{PrelimSection}, \ref{SmoothSection} hold at every point of $\Omega$, 
as opposed to the restriction to $u^{-1} \{0\}$ for formulas involving logarithmic derivatives in the optimal partition setting. 
We refer to \cite{GS92, CL08} for the details such as regularity theory, local Lipschitz estimates and classical variational formulas in this special case.
In fact, the careful reader will easily notice that since $\lambda_M = 0$ in this setting, the error terms and additive or multiplicative constants appearing in the generalized frequencies are all zero,
and the proofs carry through in a simpler fashion.
\end{proof}

\section*{Acknowledgments}

The author would like to thank Professor Fanghua Lin and Professor Changyou Wang for many helpful discussions and their interest in this work.

\bibliography{bibext}
\bibliographystyle{amsplain}

\end{document}